\documentclass[smallextended,referee,envcountsect,]{environment/svjour3}
\smartqed

\usepackage{environment/basic}
\usepackage{environment/local_convergence_specific}

\isextendedversion{
\journalname{}
}{
\journalname{JOTA}
}

\begin{document}
    
    
    \maketitle

    \begin{abstract}
        We consider Riemannian optimization problems with inequality and equality constraints and analyze a class of Riemannian interior point methods for solving them.
        The algorithm of interest consists of outer and inner iterations.
        We show that, under standard assumptions, the algorithm achieves local superlinear convergence by solving a linear system at each outer iteration, removing the need for further computations in the inner iterations.
        We also provide a specific update for the barrier parameter that achieves local near-quadratic convergence of the algorithm.
        We apply our results to the method proposed by Obara, Okuno, and Takeda (2026) and show its local superlinear and near-quadratic convergence with an analysis of the second-order stationarity.
        To our knowledge, this is the first algorithm for constrained optimization on Riemannian manifolds that achieves both local convergence and global convergence to a second-order stationary point.
        Numerical results support the theoretical analyses of the proposed methods.
    \end{abstract}

    \vskip\baselineskip

    \noindent Communicated by Shiqian Ma.
    
    \keywords{Riemannian optimization \and Nonlinear optimization \and Interior point method \and Local convergence}
    \subclass{65K05 \and 90C30}

    \section{Introduction}
        In this paper, we consider the following optimization problem:
        \isextendedversion{
        \begin{mini}
            {\pt \in \mani}{\objfun\paren*{\pt}}
            {\label{prob:RNLO}}{}
            \addConstraint{\ineqfun[\ineqidx]\paren*{\pt}}{\geq 0,}{\quad \ineqidx \in \ineqset\coloneqq\brc*{1,\ldots,\ineqdime}}
            \addConstraint{\eqfun[\eqidx]\paren*{\pt}}{= 0,}{\quad \eqidx \in \eqset\coloneqq\brc*{1,\ldots,\eqdime},}
        \end{mini}}
        {\begin{align}
            \begin{split}\label{prob:RNLO}
                \minimize[\pt \in \mani] ~ \objfun\paren*{\pt} \quad \subjectto \quad &\ineqfun[\ineqidx]\paren*{\pt}\geq 0, \text{ for all } \ineqidx \in \ineqset\coloneqq\brc*{1,\ldots,\ineqdime},\\ 
                &\eqfun[\eqidx]\paren*{\pt} = 0, \text{ for all } \eqidx \in \eqset\coloneqq\brc*{1,\ldots,\eqdime},
            \end{split}
        \end{align}}
        where $\mani$ is a $\dime$-dimensional, connected, complete Riemannian manifold, and $\objfun$, $\brc*{\ineqfun[\ineqidx]}_{\ineqidx\in\ineqset}$, and $\brc*{\eqfun[\eqidx]}_{\eqidx\in\eqset}\colon\mani\to\setR[]$ are twice continuously differentiable functions.
        If $\mani$ is not connected, we can apply our results to each connected component of $\mani$ separately.
        We call problem \cref{prob:RNLO} the Riemannian nonlinear optimization problem, abbreviated as \RNLO{}. 
        \RNLO{}~\eqref{prob:RNLO} is a natural extension of the standard nonlinear optimization problem from a Euclidean space to a Riemannian manifold. 

        Riemannian optimization, namely, optimization on Riemannian manifolds, has several strengths over classical optimization on the Euclidean space.
        First, Riemannian modeling guarantees the feasibility of the conditions defining the manifold, a property that does not hold in Euclidean models. 
        When $\mani$ is an embedded manifold, traditional Euclidean solvers may compute solutions faster than Riemannian methods, since Riemannian optimization entails additional computations for retractions. 
        Nevertheless, it is often beneficial to use Riemannian optimization when we compute more accurate solutions, as in \cite{Andreanietal2024GlobConvALMforNLPviaRiemOptim,LiuBoumal2020SimpleAlgoforOptimonRiemManiwithCstr}.
        Second, Riemannian modeling can naturally capture specific structures in applications.
        \isextendedversion{
        Examples include the positive definiteness of matrices in $H^{2}$ optimal model reduction~\cite{SatoSato17StructurePreservByRTR}, the Grassmann manifolds in low rank matrix completion~\cite{DaiMilenkovicKerman2011SubspcEvolTransfSETLowRankMatCompl,DaiKermanMilenkovic2012GeomApptoLowRankMatCompl}, the special orthogonal group in robot posture computation~\cite{BrossettEscandeKheddar2018MCPostureComputonMani} and structure-from-motion problems~\cite{Boumal23IntroOptimSmthMani}, and the flag manifold, a nested sequence of subspaces, in statistics~\cite{Yeetal2022OptimonFlagMani}. 
        Such structures can be properly addressed within the Riemannian optimization framework, whereas Euclidean approaches require ad hoc treatments and may fail to preserve these structures. 
        }{
        Examples include the special orthogonal group, which arises in robot posture computation~\cite{BrossettEscandeKheddar2018MCPostureComputonMani} and structure-from-motion problems~\cite{Boumal23IntroOptimSmthMani}, and the flag manifold, a nested sequence of subspaces, in statistics~\cite{Yeetal2022OptimonFlagMani}. 
        Such structures can be properly addressed within the Riemannian optimization framework, whereas Euclidean approaches require ad hoc treatments and may fail to preserve these structures. 
        }
        Third, exploiting the geometry of the problem can lead to high computational performance. 
        For example, geodesically convex optimization problems, which arise in various applications~\cite{SraVishnoiYildiz2018OnGeodesicConvFormBrascampLieqConst,SraHosseini2015ConicGeomOptimonManiPosDefiMat,Nguyenetal2019CalcOptimLUsingGeodsicConvOptim}, guarantee that any local optimum is also a global optimum~\cite{Boumal23IntroOptimSmthMani}. 
        Furthermore, preconditioning on a metric has been investigated to improve numerical performance~\cite{MishraSepulchre2016RiemPrecon,Gaoetal2023OptimonProdManiunderAPreconMetr}. 
    
        Riemannian optimization has been extensively studied for unconstrained settings over the last two decades~\cite{Absiletal08OptimAlgoonMatMani,Boumal23IntroOptimSmthMani}, whereas research on constrained Riemannian optimization has emerged more recently and has been rapidly developing. 
        Constrained Riemannian optimization inherits the advantages of the unconstrained framework while further expanding its applicability.
        Yang et al.~\cite{YangZhangSong14OptimCondforNLPonRiemMani} provided optimality conditions for \RNLO{}~\cref{prob:RNLO}. 
        Constraint qualifications have been investigated in \cite{BergmannHerzog2019IntrinsicFormulationKKTcondsConstrQualifSmthMani,YamakawaSato2022SeqOptimCondforNLOonRiemManiandGlobConvALM,Andreanietal2024ConstrQualifStrongGlobConvPropALMonRiemMani,Andreanietal2024GlobConvALMforNLPviaRiemOptim}. 
        Liu and Boumal~\cite{LiuBoumal2020SimpleAlgoforOptimonRiemManiwithCstr} proposed an exact penalty method combined with smoothing, as well as an augmented Lagrangian method, and established their global convergence properties. 
        These methods have recently been further developed to achieve stronger global convergence properties~\cite{YamakawaSato2022SeqOptimCondforNLOonRiemManiandGlobConvALM,Andreanietal2024ConstrQualifStrongGlobConvPropALMonRiemMani,Andreanietal2024GlobConvALMforNLPviaRiemOptim,SmthL1ExactPenaMethforIntrinsicConstRiemOptimProb}. 
        Obara et al.~\cite{ObaraOkunoTakeda2021SQOforNLOonRiemMani} proposed sequential quadratic optimization (\SQO{}) for \RNLO{}~\cref{prob:RNLO} and proved both global convergence and local quadratic convergence. 
        Schiela and Ortiz~\cite{SchielaOrtiz2021SQPMethforEqCstrOptimonHilbertMani} proposed \SQO{} for optimization problems on Hilbert manifolds with equality constraints and proved local quadratic convergence. 
        We note that these \SQO{} methods are designed to use the Riemannian Hessian operator. 
        Although computing the Hessian is generally expensive, employing second-order information enables the computation of highly accurate solutions.

        \isextendedversion{Interior point methods (\IPM{}s), which are the focus of the paper, were among the most prominent second-order algorithms for nonlinear optimization~\cite{ForsgrenGillWright2002IntMethforNLO,NocedalWright2006NumerOptim}.
        In Euclidean optimization, there are various \IPM{}s that possess local superlinear or quadratic convergence~\cite{MartiacutenezParadaTapia1995OntheCharactofQSuplinConvofQuasiNewtonIPMsforNonlinProgram,YamashitaYabe1996SupLinQuadConvSomePrimeDualIntPtMethforCstrOptim,ElBakryetal1996FormulationTheoryNewtonIntPtMethforNLP,Armandetal2007DynUpdofBarrParaminPDMethforNonlinProgram,ArmandBenoist2008ALocConvPropofPDMethforNonlinProgram,Armandetal2012FromGlobtoLocConvofIntMethsforNonlinOptim}.  
        Some algorithms have both local and global convergence properties~\cite{YamashitaYabe2003AnIPMwithPDQuadBarrPenaFunforNonlinOptim,Titsetal2003PrimalDualIntPtMethforNLPwithStrongGlobalLocalConv,WumlautachterBiegler2005LnSearchMethforNLPMotivGlobalConv,WumlautachterBiegler2005LnSearchMethforNLPLocalConv,WumlautachterBiegler2006ImplIntPtFltrLnSearchAlgoforLgeScaleNLP,LiuYuan2010NullSpcPrimalDualIntPtAlgoforNLOwithNiceConvProp,DaiLiuSun2020PrimalDualIntPtMethCpblofDetectInfeasiforNLP,WrightOrban2002PropofLogBarrFunonDegenNLP,VicenteWright2002LocalConvPDMethforDegenNonlinProg,YamashitaYabe2005QuadConvofPDIPMforDegenNonlinOptimProb}.
        Interior point trust region methods (\IPTRM{}s), which are variants of \IPM{}s, have also been designed to achieve both global and local superlinear convergence~\cite{ByrdGilbertNocedal2000TRMethBsdIntPtTechniqueforNLP,ByrdLiuNocedal1997LocalBehavIntPtMethforNLP,Connetal2000PrimalDualTRAlgoforNonconvexNLP,Gouldetal2001SuperlinConvofPDIntptAlgoforNonlinProgram,YamashitaYabeTanabe2005GlobalSuperlinConvPrimalDualIntPtTRNethforLgeScaleCstrOptim}.}{Interior point methods (\IPM{}s), which are the focus of the paper, were among the most prominent second-order algorithms for nonlinear optimization~\cite{ForsgrenGillWright2002IntMethforNLO}.
        In Euclidean optimization, there are various \IPM{}s that possess local superlinear or quadratic convergence such as~\cite{YamashitaYabe1996SupLinQuadConvSomePrimeDualIntPtMethforCstrOptim,ElBakryetal1996FormulationTheoryNewtonIntPtMethforNLP}.  
        Some algorithms have both local and global convergence properties, e.g., \cite{WumlautachterBiegler2005LnSearchMethforNLPMotivGlobalConv,WumlautachterBiegler2005LnSearchMethforNLPLocalConv,WumlautachterBiegler2006ImplIntPtFltrLnSearchAlgoforLgeScaleNLP}.
        Interior point trust region methods (\IPTRM{}s), which are variants of \IPM{}s, have also been designed to achieve both global and local superlinear convergence, for example, \cite{ByrdGilbertNocedal2000TRMethBsdIntPtTechniqueforNLP,ByrdLiuNocedal1997LocalBehavIntPtMethforNLP,YamashitaYabeTanabe2005GlobalSuperlinConvPrimalDualIntPtTRNethforLgeScaleCstrOptim}.}
        On the other hand, \IPM{}s for Riemannian optimization, abbreviated as \RIPM{}s, are relatively scarce.
        Lai and Yoshise~\cite{LaiYoshise2024RiemIntPtMethforCstrOptimonMani} proposed two types of \IPM{}s for \RNLO{}~\cref{prob:RNLO}, one possessing global convergence and the other local quadratic convergence.
        Hirai~et~al.~\cite{HiraiNieuwboerWalter2023IntptMethonManiTheoandAppl} analyzed an \IPM{} for convex optimization problems on Riemannian manifolds.
        In our companion paper~\cite{Obaraetal2025APrimalDualIPTRMfor2ndOrdStnryPtofRiemIneqCstrOptimProbs}, we proposed a Riemannian \IPTRM{} (\RIPTRM{}) for the inequality-constrained case, i.e., \RNLO{}~\cref{prob:RNLO} with $\eqset=\emptyset$, and established its global convergence to a second-order stationary point (\SOSP{}).

        Based on the above discussion and practical performance in the literature, developing interior point methods in the Riemannian setting is expected to be advantageous, as such methods can deliver highly accurate solutions while accommodating special structures with nonlinear constraints. 
        Illustrative applications include stable system identification~\cite{Obaraetal2024StabLinSysIdentifwithPriorKnwlbyRiemSQO} and $H^{2}$ optimal model reduction~\cite{MisawaSato2022H2OptimReductofPosNetusingRiemALM} in control, where the optimization involves positive definite matrices subject to nonlinear constraints.
        
    \subsection{Our Contribution}
        In this paper, we analyze a class of \RIPM{}s for solving \RNLO{}~\cref{prob:RNLO}.
        The algorithm of interest consists of a double-loop structure: an outer iteration and an inner iteration.
        In the outer iteration, the algorithm carefully selects an initial point and stopping conditions for the inner iteration.
        Then, it calls the inner iteration to compute the next iterate, generating a sequence of outer iterates to reach a local optimum under assumptions.
        Our contributions are summarized as follows: 
        \begin{enumerate}
            \item We prove the local superlinear and near-quadratic convergence of the algorithm.
            We first extend the extrapolation in \cite{Gouldetal2001SuperlinConvofPDIntptAlgoforNonlinProgram}, a methodology for computing an appropriate initial point for the inner iteration, to our \RNLO{}~\cref{prob:RNLO} setting.
            Using the extrapolation, we show that the algorithm works satisfactorily without inner iterations and then prove the local convergence.
            Note that, by the extrapolation, our method is compatible with any algorithm for the inner iteration since it is irrelevant to our convergence result.
            In addition, since the previous research~\cite{Gouldetal2001SuperlinConvofPDIntptAlgoforNonlinProgram} deals with optimization problems in the Euclidean space with linear equality and nonlinear inequality constraints, our work complements their result for the nonlinear equality- and nonlinear inequality-constrained setting even in the Euclidean case.            
            \item We establish the local superlinear and near-quadratic convergence of \RIPTRM{} for \RNLO{}~\cref{prob:RNLO} with $\eqset=\emptyset$, which was introduced by the same authors in a companion paper~\cite{Obaraetal2025APrimalDualIPTRMfor2ndOrdStnryPtofRiemIneqCstrOptimProbs}. 
            The algorithm employs a trust-region strategy to achieve global convergence to \SOSP{}s. 
            However, its local convergence had not yet been established. 
            In this paper, we analyze the second-order stationarity around a solution and prove the local convergence by applying the aforementioned results to the algorithm. 
            Combined with the global convergence analysis~\cite{Obaraetal2025APrimalDualIPTRMfor2ndOrdStnryPtofRiemIneqCstrOptimProbs}, this constitutes, to the best of our knowledge, the first algorithm for constrained optimization on Riemannian manifolds that achieves both local convergence and global convergence to an \SOSP{}.
            \item
            Our numerical results support the theoretical analyses of the proposed methods.
            In the experiments, our methods consistently exhibit the local convergence properties and yield solutions with high accuracy.
        \end{enumerate}

        A theoretical comparison between our algorithms and existing methods is provided in \cref{table:RIPMcomparison}.

        \begin{table}[t]
        \centering
        \caption{
        Comparison of algorithms for Riemannian optimization with nonlinear constraints.
        The symbols $\ineqset$ and $\eqset$ indicate the ability to handle inequality and equality constraints, respectively.
        The column ``Global'' summarizes global convergence guarantees: \KKT{}, \AKKT{}, and \SOSP{} correspond, respectively, to convergence to a Karush–Kuhn–Tucker point, an approximate-\KKT{} point, and a second-order stationary point.
        }
        \begin{tabular}{ccccc}
        \hline
        Reference & Method & Constraints        & Global      & Local          \\ \hline
        \cite{LiuBoumal2020SimpleAlgoforOptimonRiemManiwithCstr,YamakawaSato2022SeqOptimCondforNLOonRiemManiandGlobConvALM,Andreanietal2024ConstrQualifStrongGlobConvPropALMonRiemMani,Andreanietal2024GlobConvALMforNLPviaRiemOptim} & \ALM{}    & $\ineqset, \eqset$ & \AKKT{}, \SOSP{}$^{\dagger}$ &                \\
        \cite{LiuBoumal2020SimpleAlgoforOptimonRiemManiwithCstr,SmthL1ExactPenaMethforIntrinsicConstRiemOptimProb} & \EPM{}    & $\ineqset, \eqset$ & \AKKT{}        &                \\
        \cite{ObaraOkunoTakeda2021SQOforNLOonRiemMani} & \SQO{}    & $\ineqset, \eqset$ & \KKT{}         & quadratic      \\
        \cite{SchielaOrtiz2021SQPMethforEqCstrOptimonHilbertMani} & \SQO{}    & $\eqset$           &         & quadratic      \\
        \cite[Algorithm~2]{LaiYoshise2024RiemIntPtMethforCstrOptimonMani} & \IPM{}    & $\ineqset, \eqset$ &             & quadratic      \\
        \cite[Algorithm~5]{LaiYoshise2024RiemIntPtMethforCstrOptimonMani} & \IPM{}    & $\ineqset, \eqset$ & \KKT{}         &                \\
        \textbf{Algorithm~1} & \IPM{}    & $\ineqset, \eqset$ &             & near-quadratic \\
        \textbf{Algorithm~2} & \IPTRM{}  & $\ineqset$         & \SOSP{}$^{\ddagger}$        & near-quadratic \\ \hline
        \end{tabular}
        
        \footnotesize{
            $\dagger$~Global convergence to several variants of \AKKT{} points is investigated in~\cite{YamakawaSato2022SeqOptimCondforNLOonRiemManiandGlobConvALM,Andreanietal2024ConstrQualifStrongGlobConvPropALMonRiemMani,Andreanietal2024GlobConvALMforNLPviaRiemOptim}. 
            The \ALM{} in~\cite{LiuBoumal2020SimpleAlgoforOptimonRiemManiwithCstr} achieves global convergence to a \KKT{} point and, when $\ineqset=\emptyset$, an \SOSP{}.\\
            $\ddagger$~The global convergence analysis can be found in the companion paper~\cite{Obaraetal2025APrimalDualIPTRMfor2ndOrdStnryPtofRiemIneqCstrOptimProbs}.}
        \label{table:RIPMcomparison}
        \end{table}

    \subsection{Paper Organization}
        The rest of this paper is organized as follows.
        In \cref{sec:preliminaries}, we review fundamental concepts from Riemannian geometry and Riemannian optimization.
        In \cref{sec:propmeth}, we state a class of \RIPM{}s of our interest.
        In \cref{sec:localconv}, we prove the local convergence property. 
        In \cref{sec:applicationtoRIPTRM}, we also provide an application of the aforementioned results to the \RIPTRM{} in \cite{Obaraetal2025APrimalDualIPTRMfor2ndOrdStnryPtofRiemIneqCstrOptimProbs} and show its local convergence.
        In \cref{sec:experiments}, we conduct numerical experiments.
        In \cref{sec:conclusion}, we summarize our research and discuss future work.
        \isextendedversion{}{All omitted proofs can be found in the extended version available on arXiv~\cite{Obaraetal2025LocalConvofRiemIPMs}.}

    \section{Preliminaries}\label{sec:preliminaries}
        Define $\setNz\coloneq\brc*{0, 1, 2, \ldots}$ and 
        let $\setR[\dime], \setRp[\dime]$, and $\setRpp[\dime]$ be $\dime$-dimensional Euclidean space, its nonnegative orthant, and its positive orthant, respectively. 
        We denote by $\onevec[\ineqdime]$ the $\ineqdime$-dimensional vector of ones. 
        We omit the subscript $\ineqdime$ when the context is clear. 
        A continuous function $\forcingfun\colon\setRp[]\to\setRp[]$ is said to be a forcing function if $\forcingfun\paren*{\barrparam[]}=0$ if and only if $\barrparam[]=0$.
        For related positive quantities $\orderconstone$ and $\orderconsttwo$, we write $\orderconstone=\bigO[\orderconsttwo]$ if there exists a constant $\ordercoeff > 0$ such that $\orderconstone\leq\ordercoeff\orderconsttwo$ for all $\orderconsttwo$ sufficiently small.
        We write $\orderconstone=\smallO[\orderconsttwo]$ if $\lim_{\orderconsttwo\to 0}\frac{\orderconstone}{\orderconsttwo}=0$ holds.
        We also write $\orderconstone=\bigomega[\orderconsttwo]$ if $\orderconsttwo=\bigO[\orderconstone]$, and $\orderconstone=\bigtheta[\orderconsttwo]$ if $\orderconstone=\bigO[\orderconsttwo]$ and $\orderconsttwo=\bigO[\orderconstone]$.
        Given two normed vector spaces $\vecspc, \vecspctwo$ and a linear operator $\opr\colon\vecspc\to\vecspctwo$, we define the operator norm as $\opnorm{\opr}\coloneqq\sup\brc*{\norm{\opr\vecone}_{\vecspctwo} \colon \vecone \in \vecspc \text{ and } \norm{\vecone}_{\vecspc} \leq 1}$, where $\norm{\plchold}_{\vecspc}$ and $\norm{\plchold}_{\vecspctwo}$ are the norms on $\vecspc$ and $\vecspctwo$, respectively.
    
    \subsection{Notation and Terminology from Riemannian Geometry}
        We briefly review some concepts from Riemannian geometry, following the notation of \cite{Absiletal08OptimAlgoonMatMani,Boumal23IntroOptimSmthMani}.
        Let $\ptone[] \in \mani$, and let $\tanspc[{\ptone[]}]\mani$ be the tangent space to $\mani$ at $\ptone[]$.
        We denote by $\paren*{\opensubset, \chart}$ a chart of $\mani$.
        Here, $\opensubset\subseteq\mani$ is an open set, and $\chart\colon\opensubset\to\chart\paren*{\opensubset}\subseteq\setR[\dime]$ is a homeomorphism.
        Given a chart $\paren*{\opensubset, \chart}$ with $\pt\in\opensubset$, we write the coordinate expressions as
        \isextendedversion{
        \begin{align}
            \loccoord{\pt} \coloneqq \chart\paren*{\pt}, \loccoord{\tanvecone[]} \coloneqq \D\chart\paren*{\pt}\sbra*{\tanvecone[]}, \text{ and } \loccoord{\funtwo} \coloneqq \funtwo \circ \inv{\chart}
        \end{align}}{$\loccoord{\pt} \coloneqq \chart\paren*{\pt}, \loccoord{\tanvecone[]} \coloneqq \D\chart\paren*{\pt}\sbra*{\tanvecone[]}$, and  $\loccoord{\funtwo} \coloneqq \funtwo \circ \inv{\chart}$}
        for any $\tanvecone[]\in\tanspc[]\mani$ and $\funtwo\colon\mani\to\setR[]$.
        A vector field on $\mani$ is a map $\funthr\colon\mani\to\tanspc[]\mani$ with $\funthr\paren*{\ptone}\in\tanspc[\ptone]\mani$, where $\tanspc[]\mani$ is the tangent bundle. 
        \isextendedversion{Let $\intvl\subseteq\setR[]$ be an open interval, and let $\curve\colon\intvl\to\mani$ be a smooth curve.}{}
        A Riemannian metric on $\mani$ is a choice of inner product $\metr[{\ptone[]}]{\plchold}{\plchold}\colon\tanspc[{\ptone[]}]\mani\times\tanspc[{\ptone[]}]\mani\to\setR[]$ for every $\ptone\in\mani$ satisfying that, for all smooth vector fields $\funthr, \funfou$ on $\mani$, the function $\ptone\mapsto\metr[{\ptone[]}]{\funthr\paren*{\ptone}}{\funfou\paren*{\ptone}}$ is smooth from $\mani$ to $\setR[]$.
        A Riemannian manifold is a smooth manifold endowed with a Riemannian metric.
        The Riemannian metric induces the norm $\Riemnorm[\ptone]{\tanvecone[\ptone]} \coloneqq \sqrt{\metr[{\ptone[]}]{\tanvecone[\ptone]}{\tanvecone[\ptone]}}$ for $\tanvecone[\ptone]\in\tanspc[\ptone]\mani$ and the Riemannian distance $\Riemdist{\plchold}{\plchold}\colon\mani\times\mani\to\setR$.
        It follows from \cite[Theorem 13.29]{Lee12IntrotoSmthManibook2ndedn} that $\mani$ is a metric space under the Riemannian distance.
        According to the Hopf-Rinow theorem, every closed bounded subset of $\mani$ is compact for any finite-dimensional, connected, complete Riemannian manifold by regarding $\mani$ as a metric space~\cite[Chapter 7, Theorem 2.8]{doCarmo92RiemGeo}. 
        \isextendedversion{
        The following lemma shows the local equivalence of the Riemannian distance and the Euclidean distance in $\chart\paren*{\opensubset}$.
        \begin{lemma}[{\hspace{-0.04em}\cite[Lemma~14.1]{GallivanQiAbsil2012RiemDennisMoreCond}}]\label{lemm:EucliRiemdistequiv}
            Let $\mani$ be a Riemannian manifold, $\paren*{\opensubset, \chart}$ be a chart, and $\compactsubset\subseteq\opensubset$ be a compact subset, respectively.
            Then, there exist $\constsix[1], \constsix[2] \in \setRpp[]$ such that, for all $\pt[1], \pt[2] \in \compactsubset$, 
            \isextendedversion{
            \begin{align}
                \constsix[1]\norm*{\loccoord{\pt[1]} - \loccoord{\pt[2]}} \leq \Riemdist{\pt[1]}{\pt[2]} \leq \constsix[2]\norm*{\loccoord{\pt[1]} - \loccoord{\pt[2]}},
            \end{align}}
            {$\constsix[1]\norm*{\loccoord{\pt[1]} - \loccoord{\pt[2]}} \leq \Riemdist{\pt[1]}{\pt[2]} \leq \constsix[2]\norm*{\loccoord{\pt[1]} - \loccoord{\pt[2]}}$,}
            where $\norm*{\cdot}$ denotes the Euclidean norm.
        \end{lemma}}{}

        For two smooth manifolds $\mani_{1}, \mani_{2}$ and a differentiable map $\funone\colon\mani_{1}\to\mani_{2}$, we denote by $\D\funone\paren*{\ptone[]}\colon\tanspc[\ptone]\mani_{1}\to\tanspc[\funone\paren*{\ptone}]\mani_{2}$ the differential of $\funone$ at $\ptone\in\mani_{1}$.
        We use the canonical identification $\tanspc[\ptone]\vecspc\cid\vecspc$ for a vector space $\vecspc$ and $\ptone\in\vecspc$.
        Let $\funset\paren*{\mani}$ denote the set of sufficiently differentiable real-valued functions.
        Given $\funtwo\in\funset\paren*{\mani}$, $\D\funtwo\paren*{\ptone}\sbra*{\tanvecone[\ptone]} \in \tanspc[\funtwo\paren*{\ptone}]\setR[]\cid\setR$ is the differential of $\funtwo$ at $\ptone\in\mani$ along $\tanvecone[\ptone]\in\tanspc[\ptone]\mani$.
        The Riemannian gradient of $\funtwo$ at $\ptone$, denoted by $\gradstr\funtwo\paren*{\ptone}$, is defined as a unique element of $\tanspc[\ptone]\mani$ that satisfies
        \begin{align}\label{eq:riemgrad}
            \metr[\ptone]{\gradstr\funtwo\paren*{\ptone}}{\tanvecone[\ptone]} = \D\funtwo\paren*{\ptone}\sbra*{\tanvecone[\ptone]}    
        \end{align}
        for any $\tanvecone[\ptone]\in\tanspc[\ptone]\mani$.
        Here, $\gradstr\funtwo\colon\mani\to\tanspc[]\mani\colon\ptone\mapsto\gradstr\funtwo\paren*{\ptone}$ is the gradient vector field.
        Let $\Riemcxt{}{}$ be the Levi-Civita connection and $\vecfldset\paren*{\mani}$ be the set of sufficiently smooth vector fields. 
        For any $\funthr\in\vecfldset\paren*{\mani}$, we define the Jacobian of $\funthr$ at $\ptone$ as $\Jacobian[\funthr]\paren*{\ptone}\colon\tanspc[\ptone]\mani\to\tanspc[\ptone]\mani\colon\tanvecone[\ptone]\mapsto\Riemcxt{\tanvecone[\ptone]}{\funthr}$.
        \isextendedversion{
        The nonsingularity of the Jacobian and the boundedness of its inverse around a nonsingular point are shown in the following lemma:
        \begin{lemma}[{\hspace{-0.04em}\cite[Lemma~3.2]{FernandesFerreiraYuan2017SuplinConvofNewtonMethRiemMani}, \cite[Lemma~4.6]{LaiYoshise2024RiemIntPtMethforCstrOptimonMani} \unskip}]\label{lemm:Jacobnonsingbouundardaccum}
        Let $\funthr\in\vecfldset\paren*{\mani}$.
        If $\Jacobian[\funthr]$ is continuous at $\ptaccum\in\mani$ and $\Jacobian[\funthr]\paren*{\ptaccum}$ is nonsingular, then there exist $\consteig > 0$ and a neighborhood $\subsetmani\subseteq\mani$ of $\ptaccum$ such that, for all $\pt[]\in\subsetmani$, $\Jacobian[\funthr]\paren*{\pt[]}$ is nonsingular and $\opnorm{\inv{\Jacobian[\funthr]\paren*{\pt}}}\leq\consteig$.
        \end{lemma}}{}
        In particular, for the case $\funthr=\gradstr\funtwo$, the operator $\Hess\funtwo\paren*{\ptone}\colon\tanspc[\ptone]\mani\to\tanspc[\ptone]\mani$ denotes the Riemannian Hessian of $\funtwo$ at $\ptone$; that is, $\Hess\funtwo\paren*{\ptone}\sbra*{\tanvecone[\ptone]}\coloneqq\Riemcxt{\tanvecone[\ptone]}{\gradstr\funtwo}$ for all $\tanvecone[\ptone]\in\tanspc[\ptone]\mani$.
        When $\mani$ is a Euclidean space, we have $\Hess\funtwo\paren*{\ptone}\sbra*{\tanvecone[\ptone]} = \D\paren*{\gradstr\funtwo}\paren*{\ptone}\sbra*{\tanvecone[\ptone]}$
        and $\metr[\ptone]{\Hess\funtwo\paren*{\ptone}\sbra*{\tanvecone[\ptone]^{1}}}{\tanvecone[\ptone]^{2}} = \D[2]\funtwo\paren*{\ptone}\sbra*{\tanvecone[\ptone]^{1}, \tanvecone[\ptone]^{2}}$ for all $\tanvecone[\ptone], \tanvecone[\ptone]^{1}, \tanvecone[\ptone]^{2} \in \tanspc[\ptone]\mani$.  

        \isextendedversion{For each $\pt\in\mani$, let $\expmap[\pt]\colon\tanspc[\pt]\mani\to\mani$ denote the exponential map at $\pt$, the map satisfying that $\tmefiv\mapsto\expmap[\pt]\paren*{\tmefiv\tanvecone[\pt]}$ is the unique geodesic passing through $\pt$ with velocity $\tanvecone[\pt]\in\tanspc[\pt]\mani$ at $\tmefiv=0$.
        The injectivity radius at $\pt$ is defined as 
        \begin{align}
            \injradius\paren*{\pt}\coloneqq\sup\brc*{\injradiusval > 0 \colon \restrfun{\expmap[\pt]}{\brc*{\tanvecone[\pt]\in\tanspc[\pt]\mani\colon\Riemnorm[\pt]{\tanvecone[\pt]} < \injradiusval}} \text{ is a diffeomorphism}}.            
        \end{align}
        Note that $\injradius\paren*{\pt} > 0$ for any $\pt\in\mani$.
        For any $\pt[1], \pt[2]\in\mani$ with $\Riemdist{\pt[1]}{\pt[2]} < \injradius\paren*{\pt[1]}$, there is a unique minimizing geodesic connecting $\pt[1]$ and $\pt[2]$.
        Given a smooth curve $\curve$ connecting $\pt[1], \pt[2] \in \mani$,  we denote by $\partxp[\curve]{\pt[2]}{\pt[1]}\colon\tanspc[{\pt[1]}]\mani\to\tanspc[{\pt[2]}]\mani$ the parallel transport of $\tanspc[{\pt[1]}]\mani$ to $\tanspc[{\pt[2]}]\mani$.
        We omit the superscript $\curve$ when $\Riemdist{\pt[1]}{\pt[2]} < \injradius\paren*{\pt[1]}$ holds, and we use the unique minimizing geodesic as the curve $\curve$.}
        {Let $\injradius\paren*{\pt} > 0$ be the injectivity radius at $\pt\in\mani$.
        For any $\pt[1], \pt[2]\in\mani$ with $\Riemdist{\pt[1]}{\pt[2]} < \injradius\paren*{\pt[1]}$, we denote by $\partxp[]{\pt[2]}{\pt[1]}\colon\tanspc[{\pt[1]}]\mani\to\tanspc[{\pt[2]}]\mani$ the parallel transport of $\tanspc[{\pt[1]}]\mani$ to $\tanspc[{\pt[2]}]\mani$ along a unique minimizing geodesic connecting $\pt[1]$ and $\pt[2]$.}
        Note that the parallel transport is isometric, i.e., $\Riemnorm[{\pt[2]}]{\partxp[]{\pt[2]}{\pt[1]}\sbra*{\tanvecone[{\pt[1]}]}} = \Riemnorm[{\pt[1]}]{\tanvecone[{\pt[1]}]}$ for any $\tanvecone[{\pt[1]}]\in\tanspc[{\pt[1]}]\mani$ and that $\partxp[]{\pt}{\pt}\colon\tanspc[\pt]\mani\to\tanspc[\pt]\mani$ is the identity map.
        The adjoint of the parallel transport coincides with its inverse, that is, $\metr[{\pt[2]}]{\partxp[]{\pt[2]}{\pt[1]}\sbra*{\tanvecone[{\pt[1]}]}}{\tanvectwo[{\pt[2]}]} = \metr[{\pt[1]}]{\tanvecone[{\pt[1]}]}{\partxp[]{\pt[1]}{\pt[2]}\sbra*{\tanvectwo[{\pt[2]}]}}$ for any $\tanvecone[{\pt[1]}]\in\tanspc[{\pt[1]}]\mani$ and any $\tanvectwo[{\pt[2]}]\in\tanspc[{\pt[2]}]\mani$.
        
        A retraction $\retr[]\colon\tanspc[]\mani\to\mani$ is a smooth map with the following properties: by letting $\retr[\ptone]\colon\tanspc[\ptone]\mani\to\mani$ be the restriction of $\retr[]$ to $\tanspc[\ptone]\mani$, it satisfies
        \isextendedversion{
        \begin{align}
            &\retr[\ptone]\paren*{\zerovec[\ptone]} = \ptone,\label{eq:retrzero}\\
            &\D\retr[\ptone]\paren*{\zerovec[\ptone]} = \id[{\tanspc[\ptone]\mani}]\label{eq:retrdiffzero}
        \end{align}
        }{$\retr[\ptone]\paren*{\zerovec[\ptone]} = \ptone$ and $\D\retr[\ptone]\paren*{\zerovec[\ptone]} = \id[{\tanspc[\ptone]\mani}]$}
        under $\tanspc[{\zerovec[\ptone]}]\paren*{\tanspc[\ptone]\mani}\cid\tanspc[\ptone]\mani$, where $\zerovec[\ptone]$ is the zero vector of $\tanspc[\ptone]\mani$ and $\id[{\tanspc[\ptone]\mani}]$ denotes the identity map on $\tanspc[\ptone]\mani$.
        Let 
        \isextendedversion{
        \begin{align}\label{def:pullbackfun}
            \pullback[]{\funtwo}\coloneqq\funtwo\circ\retr[] \text{ and } \pullback[\ptone]{\funtwo}\coloneqq\funtwo\circ\retr[\ptone]
        \end{align}}
        {$\pullback[]{\funtwo}\coloneqq\funtwo\circ\retr[]$ and $\pullback[\ptone]{\funtwo}\coloneqq\funtwo\circ\retr[\ptone]$}
        denote the pullback of the function $\funtwo\colon\mani\to\setR[]$ and the restriction of $\pullback{\funtwo}$ to $\tanspc[\ptone]\mani$, respectively.
        Note that it follows from \isextendedversion{\cref{eq:retrzero,eq:retrdiffzero}}{the definition of the retraction} that
        \begin{align}\label{eq:pullbackgradzero}
            \gradstr\pullback[]{\funtwo}\paren*{\zerovec[\ptone]} = \gradstr\funtwo\paren*{\ptone}.
        \end{align}
        \isextendedversion{
        The equivalence of the Riemannian norm between two tangent vectors and the Riemannian distance between their retracted points is shown as follows:
        \begin{lemma}[{\hspace{-0.04em}\cite[Lemma~2]{HuangAbsilGallivan2015RiemSymRankOneTRM}}]\label{lemm:retrtanvecequiv}
            Let $\mani$ be a Riemannian manifold endowed with a retraction $\retr[]$ and let $\pt[1]\in\mani$.
            There exist positive scalars $\constsix[1], \constsix[2], \tholdvalsev[{\constsix[1], \constsix[2]}] \in \setRpp[]$ such that, for all $\pt[2]$ in a sufficiently small neighborhood of $\pt[1]$ and for all $\tanvecone[{\pt[2]}], \tanvectwo[{\pt[2]}] \in \tanspc[{\pt[2]}]\mani$ with $\Riemnorm[{\pt[2]}]{\tanvecone[{\pt[2]}]} \leq \tholdvalsev[{\constsix[1], \constsix[2]}]$ and $\Riemnorm[{\pt[2]}]{\tanvectwo[{\pt[2]}]} \leq \tholdvalsev[{\constsix[1], \constsix[2]}]$,
            \isextendedversion{
            \begin{align}
                \constsix[1] \Riemnorm[{\pt[2]}]{\tanvecone[{\pt[2]}] - \tanvectwo[{\pt[2]}]} \leq \Riemdist{\retr[{\pt[2]}]\paren*{\tanvecone[{\pt[2]}]}}{\retr[{\pt[2]}]\paren*{\tanvectwo[{\pt[2]}]}} \leq \constsix[2] \Riemnorm[{\pt[2]}]{\tanvecone[{\pt[2]}] - \tanvectwo[{\pt[2]}]}. 
            \end{align}}
            {$\constsix[1] \Riemnorm[{\pt[2]}]{\tanvecone[{\pt[2]}] - \tanvectwo[{\pt[2]}]} \leq \Riemdist{\retr[{\pt[2]}]\paren*{\tanvecone[{\pt[2]}]}}{\retr[{\pt[2]}]\paren*{\tanvectwo[{\pt[2]}]}} \leq \constsix[2] \Riemnorm[{\pt[2]}]{\tanvecone[{\pt[2]}] - \tanvectwo[{\pt[2]}]}$ holds.}
        \end{lemma}}{}
        A second-order retraction is a retraction satisfying that, for all $\ptone[]\in\mani$ and all $\tanvecone[\ptone]\in\tanspc[\ptone]\mani$, the curve $\curve\paren*{\tmefiv}=\retr[\ptone]\paren*{\tmefiv\tanvecone[\ptone]}$ has zero acceleration at $\tmefiv=0$.
        Note that second-order retractions are neither restrictive nor extremely expensive in practice.
        For instance, metric projection, which is commonly used in first-order methods and computationally efficient in many applications, is itself a second-order retraction~\cite[Proposition~5.55]{Boumal23IntroOptimSmthMani}. 
        Another example of a second-order retraction is the exponential map~\cite[Section~5.5]{Absiletal08OptimAlgoonMatMani}, which can be less expensive to compute depending on the manifold.
        Second-order retractions have the following property:
        \begin{proposition}[
        {\hspace{-0.5em}\cite[Proposition~5.5.5]{Absiletal08OptimAlgoonMatMani}}\unskip
        ]\label{prop:secodrretrHess}
            If the retraction $\retr[]$ is second order, then, for any $\funtwo\in\funset\paren*{\mani}$,
            \isextendedversion{\begin{align}\label{eq:secondretrHesszero}
                \Hess\pullback[\ptone]{\funtwo}\paren*{\zerovec[\ptone]} = \Hess\funtwo\paren*{\ptone},
            \end{align}}{$\Hess\pullback[\ptone]{\funtwo}\paren*{\zerovec[\ptone]} = \Hess\funtwo\paren*{\ptone}$ holds,}
            where the left-hand side is the Hessian of $\pullback[\ptone]{\funtwo}\colon\tanspc[\ptone]\mani\to\setR$ at $\zerovec[\ptone]\in\tanspc[\ptone]\mani$.
        \end{proposition}

    \subsection{Optimality Conditions for \RNLO{}~\cref{prob:RNLO}}
        Let $\ineqfun[]\paren*{\pt}\coloneqq\trsp{\paren*{\ineqfun[1]\paren*{\pt},\ldots,\ineqfun[\ineqdime]\paren*{\pt}}}\in\setR[\ineqdime]$ and $\eqfun[]\paren*{\pt}\coloneqq\trsp{\paren*{\eqfun[1]\paren*{\pt},\ldots,\eqfun[\eqdime]\paren*{\pt}}}\in\setR[\eqdime]$ for $\pt\in\mani$.
        We define the Lagrangian of \RNLO{}~\cref{prob:RNLO} as
        \isextendedversion{
        \begin{align}
            \Lagfun\paren*{\allvar[]} \coloneqq \objfun\paren*{\pt[]} - \trsp{\ineqLagmult[]}\ineqfun[]\paren*{\pt} + \trsp{\eqLagmult[]}\eqfun[]\paren*{\pt},
        \end{align}}
        {$\Lagfun\paren*{\allvar[]} \coloneqq \objfun\paren*{\pt[]} - \trsp{\ineqLagmult[]}\ineqfun[]\paren*{\pt} + \trsp{\eqLagmult[]}\eqfun[]\paren*{\pt}$,}
        where $\allvar[] \coloneqq \paren*{\pt[], \ineqLagmult[], \eqLagmult[]} \in \mani \times \setR[\ineqdime] \times\setR[\eqdime]$ and $\ineqLagmult[]\in\setR[\ineqdime]$ and $\eqLagmult[]\in\setR[\eqdime]$  are the vectors of Lagrange multipliers for the inequality and equality constraints, respectively.
        The Riemannian gradient and the Riemannian Hessian of the Lagrangian with respect to $\pt[]\in\mani$ are represented as 
        \isextendedversion{
        \begin{align}
            &\gradstr[\pt]\Lagfun\paren*{\allvar[]} = \gradstr\objfun\paren*{\pt[]} - \sum_{\ineqidx\in\ineqset} \ineqLagmult[\ineqidx]\gradstr\ineqfun[\ineqidx]\paren*{\pt[]} + \sum_{\eqidx\in\eqset} \eqLagmult[\eqidx]\gradstr\eqfun[\eqidx]\paren*{\pt[]}, \label{eq:gradLagfundef}\\
            &\Hess[\pt]\Lagfun\paren*{\allvar[]} = \Hess\objfun\paren*{\pt[]} - \sum_{\ineqidx\in\ineqset}\ineqLagmult[\ineqidx]\Hess\ineqfun[\ineqidx]\paren*{\pt[]} + \sum_{\eqidx\in\eqset}\eqLagmult[\eqidx]\Hess\eqfun[\eqidx]\paren*{\pt[]},
        \end{align}}
        {$\gradstr[\pt]\Lagfun\paren*{\allvar[]} = \gradstr\objfun\paren*{\pt[]} - \sum_{\ineqidx\in\ineqset} \ineqLagmult[\ineqidx]\gradstr\ineqfun[\ineqidx]\paren*{\pt[]} + \sum_{\eqidx\in\eqset} \eqLagmult[\eqidx]\gradstr\eqfun[\eqidx]\paren*{\pt[]}$ and $\Hess[\pt]\Lagfun\paren*{\allvar[]} = \Hess\objfun\paren*{\pt[]} - \sum_{\ineqidx\in\ineqset}\ineqLagmult[\ineqidx]$ $\Hess\ineqfun[\ineqidx]\paren*{\pt[]} + \sum_{\eqidx\in\eqset}\eqLagmult[\eqidx]\Hess\eqfun[\eqidx]\paren*{\pt[]}$,}
        respectively.
        We often write $\gradstr[\pt]\Lagfun\paren*{\pt,\ineqLagmult[], \eqLagmult[]}$ for $\gradstr[\pt]\Lagfun\paren*{\allvar[]}$.
        Given $\paren*{\pt,\ineqLagmult[]}\in\mani\times\setR[\ineqdime]$, for any $\vecone\in\setR[\ineqdime]$, any $\eqvecone\in\setR[\eqdime]$, and all $\tanvecone[\ptone]\in\tanspc[\ptone]\mani$, 
        \isextendedversion{
        we define the following maps:
        \begin{align}
            &\Ineqfunmat[]\paren*{\pt} \coloneqq \diag\paren*{\ineqfun[]\paren*{\pt}}\in\setR[\ineqdime\times\ineqdime], \\
            &\IneqLagmultmat[] \coloneqq \diag\paren*{\ineqLagmult[]}\in\setR[\ineqdime\times\ineqdime],\\
            &\ineqgradopr[\pt]\sbra*{\vecone} \coloneqq \sum_{\ineqidx\in\ineqset} \vecone[\ineqidx] \gradstr\ineqfun[\ineqidx]\paren*{\pt}\in\tanspc[\pt]\mani,\\ 
            &\eqgradopr[\pt]\sbra*{\eqvecone} \coloneqq \sum_{\eqidx\in\eqset} \eqvecone[\eqidx] \gradstr\eqfun[\eqidx]\paren*{\pt}\in\tanspc[\pt]\mani,\\ 
            &\coineqgradopr[\pt]\sbra*{\tanvecone[\ptone]} \coloneqq \trsp{\paren*{\metr[\pt]{\gradstr\ineqfun[1]\paren*{\pt}}{\tanvecone[\ptone]}, \ldots,\metr[\pt]{\gradstr\ineqfun[\ineqdime]\paren*{\pt}}{\tanvecone[\ptone]}}}\in\setR[\ineqdime],\\
            &\coeqgradopr[\pt]\sbra*{\tanvecone[\ptone]} \coloneqq \trsp{\paren*{\metr[\pt]{\gradstr\eqfun[1]\paren*{\pt}}{\tanvecone[\ptone]}, \ldots,\metr[\pt]{\gradstr\eqfun[\eqdime]\paren*{\pt}}{\tanvecone[\ptone]}}}\in\setR[\eqdime],
        \end{align}}
        {
        we define $\Ineqfunmat[]\paren*{\pt} \coloneqq \diag\paren*{\ineqfun[]\paren*{\pt}}\in\setR[\ineqdime\times\ineqdime]$, $\IneqLagmultmat[] \coloneqq \diag\paren*{\ineqLagmult[]}\in\setR[\ineqdime\times\ineqdime]$, $\ineqgradopr[\pt]\sbra*{\vecone} \coloneqq \sum_{\ineqidx\in\ineqset} \vecone[\ineqidx] \gradstr\ineqfun[\ineqidx]\paren*{\pt}\in\tanspc[\pt]\mani$, $\eqgradopr[\pt]\sbra*{\eqvecone} \coloneqq \sum_{\eqidx\in\eqset} \eqvecone[\eqidx] \gradstr\eqfun[\eqidx]\paren*{\pt}\in\tanspc[\pt]\mani$, $\coineqgradopr[\pt]\sbra*{\tanvecone[\ptone]} \coloneqq \trsp{\paren*{\metr[\pt]{\gradstr\ineqfun[1]\paren*{\pt}}{\tanvecone[\ptone]}, \ldots,\metr[\pt]{\gradstr\ineqfun[\ineqdime]\paren*{\pt}}{\tanvecone[\ptone]}}}$$\in\setR[\ineqdime]$, and  $\coeqgradopr[\pt]\sbra*{\tanvecone[\ptone]} \coloneqq$\\
        $\trsp{\paren*{\metr[\pt]{\gradstr\eqfun[1]\paren*{\pt}}{\tanvecone[\ptone]}, \ldots,\metr[\pt]{\gradstr\eqfun[\eqdime]\paren*{\pt}}{\tanvecone[\ptone]}}}$$\in\setR[\eqdime]$,}
        where $\diag\colon\setR[\ineqdime]\to\setR[\ineqdime\times\ineqdime]$ is the diagonal operator.
        Let 
        \isextendedversion{
        \begin{align}
            \feasirgn \coloneqq \brc*{\pt[] \in \mani \relmiddle{|} \ineqfun[\ineqidx]\paren*{\pt[]} \geq 0 \text{ for all } \ineqidx\in\ineqset} \text{ and } \strictfeasirgn \coloneqq \brc*{\pt[] \in \mani \relmiddle{|} \ineqfun[\ineqidx]\paren*{\pt[]} > 0 \text{ for all } \ineqidx\in\ineqset}
        \end{align}}
        {$\feasirgn \coloneqq \brc*{\pt[] \in \mani \relmiddle{|} \ineqfun[\ineqidx]\paren*{\pt[]} \geq 0 \text{ for all } \ineqidx\in\ineqset}$ and  $\strictfeasirgn \coloneqq \brc*{\pt[] \in \mani \relmiddle{|} \ineqfun[\ineqidx]\paren*{\pt[]} > 0 \text{ for all } \ineqidx\in\ineqset}$} be the feasible region associated with the inequality constraints
        and the strictly feasible region associated with the inequality constraints, respectively.
        We define the index set of active inequalities at $\pt[] \in \mani$ as 
        \isextendedversion{
        \begin{align}\label{def:activeineqset}
            \activeineqset\paren*{\pt[]} \coloneqq \brc*{\ineqidx \in \ineqset \relmiddle{|} \ineqfun[\ineqidx]\paren*{\pt[]} = 0}.
        \end{align}}{$\activeineqset\paren*{\pt[]} \coloneqq \brc*{\ineqidx \in \ineqset \relmiddle{|} \ineqfun[\ineqidx]\paren*{\pt[]} = 0}$.}
        
        We now introduce the optimality conditions and related concepts.       
        \begin{definition}[{\hspace{-0.04em}\cite[Equation~(4.3)]{YangZhangSong14OptimCondforNLPonRiemMani}}]
            The linear independence constraint qualification (\LICQ{}) holds at $\ptaccum\in\mani$ if $\brc*{\gradstr\ineqfun[\ineqidx]\paren*{\ptaccum}, \gradstr\eqfun[\eqidx]\paren*{\ptaccum}}_{\ineqidx\in\activeineqset\paren*{\ptaccum}, \eqidx\in\eqset}$ are linearly independent on $\tanspc[\ptaccum]\mani$.
        \end{definition}

        \begin{theorem}[{\hspace{-0.04em}\cite[Theorem~4.1]{YangZhangSong14OptimCondforNLPonRiemMani}}]
            Suppose that $\ptaccum \in \mani$ is a local minimum of \RNLO{}~\cref{prob:RNLO} and the \LICQ{} holds at $\ptaccum$.
            Then, there exist vectors of Lagrange multipliers for the inequality and equality constraints $\paren*{\ineqLagmultaccum[], \eqLagmultaccum[]}\in\setR[\ineqdime]\times\setR[\eqdime]$ such that the following hold:
            \begin{align}\label{eq:KKTconditions}
                \barrKKTvecfld\paren*{\allvaraccum}\coloneqq
                \begin{bmatrix}
                    \gradstr[\pt]\Lagfun\paren*{\allvaraccum}\\
                    \Ineqfunmat[]\paren*{\ptaccum}\ineqLagmultaccum[]\\
                    \eqfun[]\paren*{\ptaccum}
                \end{bmatrix}
                =
                \begin{bmatrix}
                    \zerovec[\ptaccum]\\
                    0\\
                    0
                \end{bmatrix}, \, 
                \ineqLagmultaccum[] \geq 0, \text{ and } \ineqfun\paren*{\ptaccum} \geq 0.
            \end{align}
            We call \cref{eq:KKTconditions} the \KKT{} conditions of \RNLO{}~\cref{prob:RNLO} and $\ptaccum$ a \KKT{} point of \RNLO{}~\cref{prob:RNLO}.
        \end{theorem}

        \begin{definition}[{\hspace{-0.04em}\cite[Definition 2.5]{ObaraOkunoTakeda2021SQOforNLOonRiemMani}}]
            Given $\ptaccum \in \feasirgn$ satisfying the \KKT{} conditions with associated Lagrange multipliers $\paren*{\ineqLagmultaccum[], \eqLagmultaccum[]}\in \setRp[\ineqdime]\times\setR[\eqdime]$, we say that the strict complementarity condition (\SC{}) holds if exactly one of $\ineqLagmultaccum$ and $\ineqfun[\ineqidx]\paren*{\ptaccum}$ is zero for each index $\ineqidx \in \ineqset$. 
            Hence, we have $\ineqLagmultaccum > 0$ for every $\ineqidx \in \activeineqset\paren*{\ptaccum}$.
        \end{definition}
        We define the critical cone associated with $\allvaraccum\in\mani\times\setRp[\ineqdime]\times\setR[\eqdime]$ as
        \begin{align}\label{def:criticalcone}
            &\criticalcone\paren*{\allvaraccum} \coloneqq \left\{ \tanvecone[\ptaccum] \in \tanspc[\ptaccum]\mani \middle|
            \begin{aligned}
                &\metr[\ptaccum]{\gradstr\ineqfun[\ineqidx]\paren*{\ptaccum}}{\tanvecone[\ptaccum]} = 0 \text{ for all } \ineqidx \in \activeineqset\paren*{\ptaccum} \text{ with } \ineqLagmultaccum > 0, \\
                &\metr[\ptaccum]{\gradstr\ineqfun[\ineqidx]\paren*{\ptaccum}}{\tanvecone[\ptaccum]} \geq 0 \text{ for all } \ineqidx \in \activeineqset\paren*{\ptaccum} \text{ with } \ineqLagmultaccum = 0, \\
                &\metr[\ptaccum]{\gradstr\eqfun[\eqidx]\paren*{\ptaccum}}{\tanvecone[\ptaccum]} = 0 \text{ for all } \eqidx \in \eqset
            \end{aligned}
            \right\}.
        \end{align}

        \begin{theorem}[{\hspace{-0.04em}\cite[Theorem~4.3]{YangZhangSong14OptimCondforNLPonRiemMani}}]
            Let $\ptaccum\in\feasirgn$ be a \KKT{} point with associated Lagrange multipliers $\paren*{\ineqLagmultaccum[], \eqLagmultaccum[]}\in\setRp[\ineqdime]\times\setR[\eqdime]$.
            Suppose
            \isextendedversion{\begin{align}\label{eq:secondordersufficientcond}
                \metr[\ptaccum]{\Hess[{\pt}]\Lagfun\paren*{\allvaraccum}\sbra*{\tanvecone[\ptaccum]}}{\tanvecone[\ptaccum]} > 0 \text{ for all } \tanvecone[\ptaccum] \in \criticalcone\paren*{\allvaraccum}\backslash\brc*{\zerovec[\ptaccum]}.
            \end{align}}
            {
            $\metr[\ptaccum]{\Hess[{\pt}]\Lagfun\paren*{\allvaraccum}\sbra*{\tanvecone[\ptaccum]}}{\tanvecone[\ptaccum]}$\\ $> 0 \text{ for all } \tanvecone[\ptaccum] \in \criticalcone\paren*{\allvaraccum}\backslash\brc*{\zerovec[\ptaccum]}$.
            }
            Then, $\ptaccum$ is a strict local minimum of \RNLO{}~\cref{prob:RNLO}.
            We call \isextendedversion{\cref{eq:secondordersufficientcond}}{the above inequality} the second-order sufficient condition (\SOSC{}).
        \end{theorem}

    \section{Proposed Method}\label{sec:propmeth}
        In this section, we propose a class of Riemannian \IPM{}s for solving \RNLO{}~\cref{prob:RNLO}.
        Let $\barrparam[] > 0$ be the barrier parameter, and let $\tildestd\paren*{\pt}\coloneqq\sbra*{\zerovec[\pt]; \onevec[\ineqdime]; \zerovec[]}\in\tanspc[\pt]\mani\times\setR[\ineqdime]\times\setR[\eqdime]$ for $\pt\in\mani$.
        We define the barrier \KKT{} vector field~\cite[Section~3.1]{LaiYoshise2024RiemIntPtMethforCstrOptimonMani} as 
        \isextendedversion{
        \begin{align}\label{eq:barrKKTvecfld}
            &\barrKKTvecfld\paren*{\allvar[];\barrparam[]}\coloneqq 
            \barrKKTvecfld\paren*{\allvar[]} - 
            \barrparam[]\tildestd\paren*{\pt}.
        \end{align}}{$\barrKKTvecfld\paren*{\allvar[];\barrparam[]}\coloneqq 
            \barrKKTvecfld\paren*{\allvar[]} - 
            \barrparam[]\tildestd\paren*{\pt}$.}
        With $\barrparam[]$ set to zero, equation $\barrKKTvecfld\paren*{\allvar[];0}=0$ coincides with the \KKT{} conditions~\cref{eq:KKTconditions} for $\allvar[]\in\feasirgn\times\setRp[\ineqdime]\times\setR[\eqdime]$.
        The algorithm aims to find a solution of $\barrKKTvecfld\paren*{\allvar[]; 0} = 0$ with $\allvar[]\in\feasirgn\times\setRp[\ineqdime]\times\setR[\eqdime]$
        by approximately solving
        \begin{align}\label{eq:barrKKTequation}
            \barrKKTvecfld\paren*{\allvar[];\barrparamotriter}=0, \quad \allvar[]\in\strictfeasirgn\times\setRpp[\ineqdime]\times\setR[\eqdime]
        \end{align}
        for a sequence of barrier parameters $\brc*{\barrparamotriter}_{\otriteridx}\subset \setRpp[]$ that converges to zero.
        Let $\allvarotriter\coloneqq\paren*{\ptotriter, \ineqLagmultotriter[], \eqLagmultotriter[]} \in \strictfeasirgn \times \setRpp[\ineqdime]\times\setR[\eqdime]$ denote an approximate solution of \cref{eq:barrKKTequation}.
        We call the sequence $\brc*{\allvarotriter}_{\otriteridx}$ the outer iterates. 
        Adjusting the barrier parameter and the tolerances for residuals defines the outer iteration, whose index is the subscript $\otriteridx\in\setNz$.
        At the $\otriteridx$-th outer iteration, the algorithm first sets $\barrparamotriter > 0$ so that $\lim_{\otriteridx \to \infty} \barrparamotriter = 0$ and then 
        computes the initial point $\allvarotriterpinit\in\strictfeasirgn\times\setRpp[\ineqdime]\times\setR[\eqdime]$ using the information from the previous iterates.
        Then, starting from this initial point, it calls the inner iteration to find the next iterate $\allvarotriterp\in\strictfeasirgn\times\setRpp[\ineqdime]\times\setR[\eqdime]$ that satisfies the following stopping conditions with $\barrparam > 0$:
        \begin{align}
            &\Riemnorm[{\ptotriterp}]{\gradstr[{\pt[]}]\Lagfun\paren*{\allvarotriterp}} \leq \forcingfungradLag\paren*{\barrparam},\label{eq:stopcondKKT}\\
            &\norm*{\Ineqfunmat[]\paren*{\ptotriterp}\ineqLagmultotriterp[] - \barrparam \onevec} \leq \forcingfuncompl\paren*{\barrparam},\label{eq:stopcondbarrcompl}\\
            &\norm*{\eqfun[]\paren*{\ptotriterp}}\leq\forcingfuneqvio\paren*{\barrparam},\label{eq:stopcondeqvio}\\
            &\ineqfun[]\paren*{\ptotriterp} > 0, \text{ and } \ineqLagmultotriterp[] > 0. \label{eq:stopcondstrictfeasi}
        \end{align}
        Any algorithm can be employed for the inner iteration as long as it computes a point satisfying \cref{eq:stopcondKKT,eq:stopcondbarrcompl,eq:stopcondeqvio,eq:stopcondstrictfeasi}; see \cref{sec:conclusion} for more discussion.
        \cref{algo:RIPMOuter} formally states the procedure of the outer iteration.

        For the initial point computation, we extend the extrapolation—originally introduced in the Euclidean setting with linear equality and nonlinear inequality constraints~\cite{Gouldetal2001SuperlinConvofPDIntptAlgoforNonlinProgram}—to \RNLO{}~\cref{prob:RNLO}.
        Applying the Newton method to \cref{eq:barrKKTequation} with $\barrparam[] = \barrparamotriter > 0$ fixed, we obtain
        the Newton equation
        \begin{align}\label{eq:RiemNewtoneq}
            \Jacobian[\barrKKTvecfld]\paren*{\allvar[]}\sbra*{\dirallvar} = - \barrKKTvecfld\paren*{\allvar[];\barrparam[]},
        \end{align}
        where, 
        under $\tanspc[{\ineqLagmult[]}]\setR[\ineqdime]\cid\setR[\ineqdime]$ and $\tanspc[{\eqLagmult[]}]\setR[\eqdime]\cid\setR[\eqdime]$,
        \begin{align}
            \begin{split}\label{def:JacobbarrKKTvecfld}
                \Jacobian[\barrKKTvecfld]\paren*{\allvar[]}\colon\tanspc[\pt]\mani\times\setR[\ineqdime]\times\setR[\eqdime]&\to\tanspc[\pt]\mani\times\setR[\ineqdime]\times\setR[\eqdime]\\
                \dirallvar\coloneqq\paren*{\dirpt, \dirineqLagmult[], \direqLagmult[]} &\mapsto 
                \begin{bmatrix}
                    \Hess[\pt]\Lagfun\paren*{\allvar}\sbra*{\dirpt} - \ineqgradopr[\pt]\sbra*{\dirineqLagmult[]}+\eqgradopr[\pt]\sbra*{\direqLagmult[]}\\
                    \IneqLagmultmat[]\coineqgradopr[\pt]\sbra*{\dirpt} + \Ineqfunmat[]\paren*{\pt}\dirineqLagmult[]\\
                    \coeqgradopr[\pt]\sbra*{\dirpt}
                \end{bmatrix}
            \end{split}
        \end{align}
        is the Jacobian of $\barrKKTvecfld\paren*{\cdot;\barrparam[]}$ at $\allvar$.
        Here, we omit $\barrparam[]$ in the notation since the Jacobian does not depend on $\barrparam[]$.
        Then, we define the extrapolation as 
        \begin{align}
            \allvarotriterpinit\leftarrow\paren*{\retr[\ptotriter]\paren*{\dirNewtonptbarrparamotriter}, \ineqLagmultotriter[] + \dirNewtonineqLagmultbarrparamotriter[], \eqLagmultotriter[] + \dirNewtoneqLagmultbarrparamotriter[]},\label{def:extrapolation}
        \end{align}
        where        $\dirNewtonallvarbarrparamotriter\in\tanspc[\ptotriter]\mani\times\setR[\ineqdime]\times\setR[\eqdime]$ is the Newton step defined as
        \begin{align}
            \dirNewtonallvarbarrparamotriter = \paren*{\dirNewtonptbarrparamotriter, \dirNewtonineqLagmultbarrparamotriter[], \dirNewtoneqLagmultbarrparamotriter[]} \coloneqq - \inv{\Jacobian[\barrKKTvecfld]\paren*{\allvarotriter}}\barrKKTvecfld\paren*{\allvarotriter;\barrparamotriter}\label{def:Newtondirection}
        \end{align}
        at any $\otriteridx$-th iteration where the Jacobian is nonsingular and the point $\allvarotriterpinit$ lies in $\strictfeasirgn\times\setRpp[\ineqdime]\times\setR[\eqdime]$; that is, $\allvarotriterpinit$ satisfies \cref{eq:stopcondstrictfeasi}.        

        \begin{algorithm2e}[t]
            \SetKwInOut{Require}{Require}
            \SetKwInOut{Input}{Input}
            \SetKwInOut{Output}{Output}
            \Require{Riemannian manifold $\mani$, twice continuously differentiable functions $\objfun$, $\brc*{\ineqfun[\ineqidx]}_{\ineqidx\in\ineqset}$, and $\brc*{\eqfun[\eqidx]}_{\eqidx\in\eqset}\colon\mani\to\setR[]$, 
            forcing functions $\forcingfungradLag,  \forcingfuncompl, \forcingfuneqvio\colon\setRp[]\to\setRp[]$.}
            \Input{Initial point $\allvar[0] = \paren*{\pt[0], \ineqLagmult[0], \eqLagmult[0]} \in \strictfeasirgn \times \setRpp[\ineqdime]\times\setR[\eqdime]$,
            initial barrier parameter $\barrparam[-1] > 0$}.
            \For{$\otriteridx = 0, 1, \ldots$}{
                Set $\barrparamotriter > 0$ so that $\lim_{\otriteridx \to \infty} \barrparamotriter = 0$.
                
                Compute the initial point $\allvarotriterpinit\in\strictfeasirgn\times\setRpp[\ineqdime]\times\setR[\eqdime]$ by \cref{def:extrapolation}.
                
                Compute $\allvarotriterp\in \strictfeasirgn \times \setRpp[\ineqdime]\times\setR[\eqdime]$ 
                that satisfies the stopping conditions \cref{eq:stopcondKKT,eq:stopcondbarrcompl,eq:stopcondeqvio,eq:stopcondstrictfeasi} using an inner iteration starting from $\allvarotriterpinit$.
            }
            \caption{Outer iteration of \RIPM{} 
            }\label{algo:RIPMOuter}
        \end{algorithm2e}

    \section{Local Convergence Analysis}\label{sec:localconv}
        In this section, we establish the local convergence property of \cref{algo:RIPMOuter}.
        Let $\allvaraccum=\paren*{\ptaccum, \ineqLagmultaccum[], \eqLagmultaccum[]}\in\feasirgn\times\setRp[\ineqdime]\times\setR[\eqdime]$ be any point satisfying $\barrKKTvecfld\paren*{\allvaraccum}=0$.
        In \cref{subsec:exactsolstopcond}, we prove that the initial point defined by \cref{def:extrapolation} satisfies the stopping conditions \cref{eq:stopcondKKT,eq:stopcondbarrcompl,eq:stopcondeqvio,eq:stopcondstrictfeasi} in a sufficiently small neighborhood of $\allvaraccum$ under certain assumptions, 
        which implies that \cref{algo:RIPMOuter} will require no further inner iterations.
        In \cref{subsec:locnearquadconv}, we prove the local convergence of \cref{algo:RIPMOuter}.
        Some complicated proofs of lemmas, a proposition, and a theorem are deferred to the appendices for readability.

    \stepcounter{assumptionsection}

    \subsection{Analysis of the Extrapolation}\label{subsec:exactsolstopcond}
        In this subsection, we first show that the point produced by the extrapolation \cref{def:extrapolation} is well-defined and satisfies \cref{eq:stopcondstrictfeasi} within a sufficiently small neighborhood of $\allvaraccum$. 
        This ensures that the point defined by~\cref{def:extrapolation} is eligible as the initial point for the inner iteration.
        Then, we prove that the initial point $\allvarotriterpinit$ by \cref{def:extrapolation} also satisfies the other conditions \cref{eq:stopcondKKT,eq:stopcondbarrcompl,eq:stopcondeqvio} if the current iterate lies in a sufficiently small neighborhood of $\allvaraccum$ and the current barrier parameter is sufficiently small.
        This implies that \cref{algo:RIPMOuter} can work without an inner iteration.
        These results will form the basis of the local convergence analysis in \cref{subsec:locnearquadconv}.

        We assume the following:
        \begin{assumption}\label{assu:SCLICQOSOSC}
            The point $\allvaraccum$ satisfies the \SC{}, the \LICQ{}, and the \SOSC{}.
        \end{assumption}
        \begin{assumption}\label{assu:forcingfuncomplgradLagbounded}
            There exist $\constgradLaglower, \constgradLagupper, \constcompllower, \constcomplupper, \consteqviolower, \consteqvioupper \in \setR[]$ with $0 < \constgradLaglower < 1 < \constgradLagupper$, $0 < \constcompllower < 1 < \constcomplupper$, and $0 < \consteqviolower < 1 < \consteqvioupper$ such that, for all $\otriteridx\in\setNz$,
            \isextendedversion{\begin{align}
                &\constgradLaglower\barrparamotriter \leq \forcingfungradLag\paren*{\barrparamotriter} \leq \constgradLagupper\barrparamotriter,\label{ineq:forcingfungradLagbounded}\\
                &\constcompllower\barrparamotriter \leq \forcingfuncompl\paren*{\barrparamotriter} \leq \constcomplupper\barrparamotriter,\label{ineq:forcingfuncomplbounded} \\
                &\consteqviolower\barrparamotriter \leq \forcingfuneqvio\paren*{\barrparamotriter} \leq \consteqvioupper\barrparamotriter.\label{ineq:forcingfuneqviobounded} 
            \end{align}}
            {$\constgradLaglower\barrparamotriter \leq \forcingfungradLag\paren*{\barrparamotriter} \leq \constgradLagupper\barrparamotriter$, $\constcompllower\barrparamotriter \leq \forcingfuncompl\paren*{\barrparamotriter} \leq \constcomplupper\barrparamotriter$, and $\consteqviolower\barrparamotriter \leq \forcingfuneqvio\paren*{\barrparamotriter} \leq \consteqvioupper\barrparamotriter$.}
        \end{assumption}
        \begin{assumption}\label{assu:barrparamdecreaseupodr}
            $\barrparamotriter^{2} = \smallO[\barrparamotriterp]$ holds.
        \end{assumption}
        \begin{assumption}\label{assu:barrparammonononincr}
            $\barrparamotriterp \leq \barrparamotriter$ holds for all $\otriteridx\in\setNz$ sufficiently large.
        \end{assumption}
        \begin{assumption}\label{assu:retrsecondord}
            The retraction $\retr[]$ is second order.
        \end{assumption}
        \begin{assumption}\label{assu:pullbackHessobjineqfunLipschitz}
            There exist positive scalars $\Lipschitzconstfiv[\objfun], \brc*{\Lipschitzconstfiv[{\ineqfun[\ineqidx]}]}_{\ineqidx\in\ineqset}, \brc*{\Lipschitzconstfiv[{\eqfun[\eqidx]}]}_{\eqidx\in\eqset} \in \setRpp[]$ such that, for all $\pt[]\in\mani$ sufficiently close to $\ptaccum\in\feasirgn$ and all $\tanvecone[\pt]\in\tanspc[\pt]\mani$ sufficiently small, 
            \isextendedversion{\begin{align}
                &\opnorm{\Hess\pullback[\pt]{\objfun}\paren*{\tanvecone[\pt]} - \Hess\pullback[\pt]{\objfun}\paren*{\zerovec[\pt]}} \leq \Lipschitzconstfiv[\objfun] \Riemnorm[\pt]{\tanvecone[\pt]},\label{ineq:pullbackHessobjfunLipschitz}\\
                &\opnorm{\Hess\pullback[\pt]{\ineqfun[\ineqidx]}\paren*{\tanvecone[\pt]} - \Hess\pullback[\pt]{\ineqfun[\ineqidx]}\paren*{\zerovec[\pt]}} \leq \Lipschitzconstfiv[{\ineqfun[\ineqidx]}] \Riemnorm[\pt]{\tanvecone[\pt]} \text{ for all } \ineqidx\in\ineqset,\label{ineq:pullbackHessineqfunLipschitz}\\
                &\opnorm{\Hess\pullback[\pt]{\eqfun[\eqidx]}\paren*{\tanvecone[\pt]} - \Hess\pullback[\pt]{\eqfun[\eqidx]}\paren*{\zerovec[\pt]}} \leq \Lipschitzconstfiv[{\eqfun[\eqidx]}] \Riemnorm[\pt]{\tanvecone[\pt]} \text{ for all } \eqidx\in\eqset.\label{ineq:pullbackHesseqfunLipschitz}
            \end{align}}{$\opnorm{\Hess\pullback[\pt]{\objfun}\paren*{\tanvecone[\pt]} - \Hess\pullback[\pt]{\objfun}\paren*{\zerovec[\pt]}} \leq \Lipschitzconstfiv[\objfun] \Riemnorm[\pt]{\tanvecone[\pt]}$, $\opnorm{\Hess\pullback[\pt]{\ineqfun[\ineqidx]}\paren*{\tanvecone[\pt]} - \Hess\pullback[\pt]{\ineqfun[\ineqidx]}\paren*{\zerovec[\pt]}} \leq \Lipschitzconstfiv[{\ineqfun[\ineqidx]}] \Riemnorm[\pt]{\tanvecone[\pt]}$ for all  $\ineqidx\in\ineqset$, and $\opnorm{\Hess\pullback[\pt]{\eqfun[\eqidx]}\paren*{\tanvecone[\pt]} - \Hess\pullback[\pt]{\eqfun[\eqidx]}\paren*{\zerovec[\pt]}} \leq \Lipschitzconstfiv[{\eqfun[\eqidx]}] \Riemnorm[\pt]{\tanvecone[\pt]}$ for all $\eqidx\in\eqset$.}
        \end{assumption}
        Note that \cref{assu:SCLICQOSOSC} is standard in the literature~\cite[Assumptions~(A2)--(A4)]{LaiYoshise2024RiemIntPtMethforCstrOptimonMani} and \cite[Assumption~B1]{ObaraOkunoTakeda2021SQOforNLOonRiemMani}.
        We will discuss the relaxation of \cref{assu:SCLICQOSOSC} in \cref{sec:conclusion}.
        \cref{assu:forcingfuncomplgradLagbounded} is fulfilled if we employ updates $\forcingfungradLag\paren*{\barrparamotriter} \leftarrow \constforcingfungradLag\barrparamotriter$, $\forcingfuncompl\paren*{\barrparamotriter} \leftarrow \constforcingfuncompl\barrparamotriter$, and $\forcingfuneqvio\paren*{\barrparamotriter} \leftarrow \constforcingfuneqvio\barrparamotriter$ for some $\constforcingfungradLag, \constforcingfuncompl, \constforcingfuneqvio \in \setRpp[]$.
        We will provide a specific update for $\brc*{\barrparamotriter}_{\otriteridx}$ that satisfies \cref{assu:barrparamdecreaseupodr,assu:barrparammonononincr} in \cref{subsec:locnearquadconv}.
        \cref{assu:pullbackHessobjineqfunLipschitz} holds if the functions $\objfun$, $\brc*{\ineqfun[\ineqidx]}_{\ineqidx\in\ineqset}$, and $\brc*{\eqfun[\eqidx]}_{\eqidx\in\eqset}$ are of class $C^{3}$~\cite[Lemma~10.57]{Boumal23IntroOptimSmthMani}.

        We first prove that the extrapolation is well-defined around $\allvaraccum$.
        The next lemma shows that the Jacobian at $\allvaraccum$ is nonsingular, following an argument similar to \cite[Theorem~3.1]{LaiYoshise2024RiemIntPtMethforCstrOptimonMani}.
        \isextendedversion{Recall the definitions of the critical cone~\cref{def:criticalcone}, the barrier \KKT{} vector field~\cref{eq:barrKKTvecfld}, and its Jacobian \cref{def:JacobbarrKKTvecfld}.}{}
        \newcommand{\lemmaJacobnonsingaccum}{\begin{lemma}\label{lemm:Jacobnonsingaccum}
            Under \cref{assu:SCLICQOSOSC}, 
            $\Jacobian[\barrKKTvecfld]\paren*{\allvaraccum}$ is nonsingular. 
        \end{lemma}
        }
        \isextendedversion{
        \lemmaJacobnonsingaccum{}
        \begin{proof}
            See \cref{appx:proofJacobnonsingaccum}.
            \qed
        \end{proof}}{        \lemmaJacobnonsingaccum{}
        }
        Since the Jacobian $\Jacobian[\barrKKTvecfld]$ is continuous, we obtain the following corollary:
        \begin{corollary}\label{coro:Jacobnonsingularardaccum}
            Under \cref{assu:SCLICQOSOSC}, 
            $\Jacobian[\barrKKTvecfld]\paren*{\allvar}$ is nonsingular for any $\allvar\in\mani\times\setR[\ineqdime]\times\setR[\eqdime]$ sufficiently close to $\allvaraccum$.
        \end{corollary}

        Next, we prove that the point produced by \cref{def:extrapolation} belongs to the interior of the feasible region for $\pt$ and $\ineqLagmult[]$.
        To this end, we first provide two auxiliary lemmas, which can be proved in similar manners to \cite[Theorem 5.1]{Gouldetal2001SuperlinConvofPDIntptAlgoforNonlinProgram} and  \cite[(6.11)]{Gouldetal2001SuperlinConvofPDIntptAlgoforNonlinProgram}, respectively.
        Note that  $\ineqLagmultotriter$ denotes the $\ineqidx$-th component of $\ineqLagmultotriter[]\in\setRpp[\ineqdime]$.
        \newcommand{\lemmaineqLagmultineqfunbarrparambounds}{\begin{lemma}\label{lemm:ineqLagmultineqfunbarrparambounds}
            Suppose the \SC{} and \cref{assu:forcingfuncomplgradLagbounded}.
            For any $\otriteridx$-th iteration of \cref{algo:RIPMOuter} satisfying that $\allvarotriter$ is sufficiently close to $\allvaraccum$, the following hold:
            \begin{enumerate}
                \item There exist $\constineqLagmultlower, \constineqLagmultupper \in \setRpp[]$ such that $\constineqLagmultupper \geq \constineqLagmultlower > 0$ and, for any $\ineqidx\in\activeineqset\paren*{\ptaccum}$,
                \isextendedversion{\begin{align}\label{ineq:ineqfunLagmultbarrparamboundactive}
                    \frac{1}{\constineqLagmultupper}\paren*{1 - \constcomplupper}\barrparamotriterm \leq \ineqfun[\ineqidx]\paren*{\ptotriter} \leq \frac{1}{\constineqLagmultlower}\paren*{1 + \constcomplupper}\barrparamotriterm \text{ and }
                    \constineqLagmultlower \leq \ineqLagmultotriter \leq \constineqLagmultupper.
                \end{align}}{$\frac{1}{\constineqLagmultupper}\paren*{1 - \constcomplupper}\barrparamotriterm \leq \ineqfun[\ineqidx]\paren*{\ptotriter} \leq \frac{1}{\constineqLagmultlower}\paren*{1 + \constcomplupper}\barrparamotriterm$ and $\constineqLagmultlower \leq \ineqLagmultotriter \leq \constineqLagmultupper$.}
                \label{lemm:activebarrparambounds}
                \item There exist $\constinactivecompllower, \constinactivecomplupper \in \setRpp[]$ such that $\constinactivecomplupper \geq \constinactivecompllower > 0$ and, for any $\ineqidx\notin\activeineqset\paren*{\ptaccum}$,
                \isextendedversion{\begin{align}\label{ineq:ineqfunLagmultbarrparamboundinactive}
                    \frac{1}{\constinactivecomplupper}\paren*{1 - \constcomplupper}\barrparamotriterm \leq \ineqLagmultotriter\leq \frac{1}{\constinactivecompllower}\paren*{1 + \constcomplupper}\barrparamotriterm \text{ and }
                    \constinactivecompllower \leq \ineqfun[\ineqidx]\paren*{\ptotriter} \leq \constinactivecomplupper.
                \end{align}}{$\frac{1}{\constinactivecomplupper}\paren*{1 - \constcomplupper}\barrparamotriterm \leq \ineqLagmultotriter\leq \frac{1}{\constinactivecompllower}\paren*{1 + \constcomplupper}\barrparamotriterm$ and $\constinactivecompllower \leq \ineqfun[\ineqidx]\paren*{\ptotriter} \leq \constinactivecomplupper$.}
                \label{lemm:inactivebarrparambounds}
            \end{enumerate}
        \end{lemma}}
        \isextendedversion{
        \lemmaineqLagmultineqfunbarrparambounds{}
        \begin{proof}
            See \cref{appx:proofineqLagmultineqfunbarrparambounds}.
        \end{proof}}{\lemmaineqLagmultineqfunbarrparambounds{}}
        
        \newcommand{\lemmadirNewtonallvarbarrparambound}{\begin{lemma}\label{lemm:dirNewtonallvarbarrparambound}
            Under \cref{assu:SCLICQOSOSC,assu:forcingfuncomplgradLagbounded,assu:barrparammonononincr}, there exists $\constdirNewton > 0$ such that $\Riemnorm[{\allvarotriter}]{\dirNewtonallvarbarrparamotriter} \leq \constdirNewton\barrparamotriterm$ holds 
            for any $\otriteridx$-th iteration of \cref{algo:RIPMOuter} satisfying that $\allvarotriter$ is sufficiently close to $\allvaraccum$ and $\barrparamotriter$ is sufficiently small.
        \end{lemma}}
        \isextendedversion{
        \lemmadirNewtonallvarbarrparambound{}
        \begin{proof}
            See \cref{appx:proofdirNewtonallvarbarrparambound}.
        \end{proof}}{\lemmadirNewtonallvarbarrparambound{}}
        Using \cref{lemm:ineqLagmultineqfunbarrparambounds,lemm:dirNewtonallvarbarrparambound}, we ensure that the point $\allvarotriterpinit$ defined by \cref{def:extrapolation} satisfies the strict feasibility \cref{eq:stopcondstrictfeasi} in a sufficiently small neighborhood of $\allvaraccum$ as follows: 
        note that by the definition of \cref{def:Newtondirection}, we have that, for all $\ineqidx\in\ineqset$,
        \begin{align}\label{eq:dirNewtoninexidxequation}
             \sbra*{\dirNewtonineqLagmultbarrparamotriter[]}_{\ineqidx}\ineqfun[\ineqidx]\paren*{\ptotriter} = - \ineqLagmultotriter[\ineqidx]\ineqfun[\ineqidx]\paren*{\ptotriter} + \barrparamotriter - \ineqLagmultotriter[\ineqidx]\D\ineqfun[\ineqidx]\paren*{\ptotriter}\sbra*{\dirNewtonptbarrparamotriter}.
        \end{align}
        \newcommand{\lemmaproofdirNewtonineqfunpos}{
        \begin{lemma}\label{lemm:dirNewtonineqfunpos}
            Under \cref{assu:SCLICQOSOSC,assu:forcingfuncomplgradLagbounded,assu:barrparamdecreaseupodr,assu:barrparammonononincr}, $\ineqfun[\ineqidx]\paren*{\retr[\ptotriter]\paren*{\dirNewtonptbarrparamotriter}} > 0$ holds for any $\ineqidx\in\ineqset$ and any $\otriteridx$-th iteration of \cref{algo:RIPMOuter} satisfying that $\allvarotriter$ is sufficiently close to $\allvaraccum$ and $\barrparamotriterm$ is sufficiently small. 
        \end{lemma}
        \begin{proof}
            See \cref{appx:proofdirNewtonineqfunpos}.
        \end{proof}}
        \isextendedversion{\lemmaproofdirNewtonineqfunpos{}}{
        \lemmaproofdirNewtonineqfunpos{}
        We can similarly prove the feasibility of the Lagrange multipliers for the inequality constraints
        by an argument analogous to the Euclidean setting~\cite[Theorem 6.2]{Gouldetal2001SuperlinConvofPDIntptAlgoforNonlinProgram}.}
        
        \newcommand{\lemmadirNewtonineqLagmultpos}{\begin{lemma}\label{lemm:dirNewtonineqLagmultpos}
            Under \cref{assu:SCLICQOSOSC,assu:forcingfuncomplgradLagbounded,assu:barrparamdecreaseupodr,assu:barrparammonononincr}, $\ineqLagmultotriter[] + \dirNewtonineqLagmultbarrparamotriter[] > 0$ holds for any $\otriteridx$-th iteration of \cref{algo:RIPMOuter} satisfying that $\allvarotriter$ is sufficiently close to $\allvaraccum$ and $\barrparamotriterm$ is sufficiently small. 
        \end{lemma}}
        \isextendedversion{
        \lemmadirNewtonineqLagmultpos{}
        \begin{proof}
            See \cref{appx:proofdirNewtonineqLagmultpos}.
        \end{proof}}{\lemmadirNewtonineqLagmultpos{}}
        From \cref{coro:Jacobnonsingularardaccum,lemm:dirNewtonineqfunpos,lemm:dirNewtonineqLagmultpos}, we obtain the following corollary.
        \begin{corollary}\label{coro:extrapolationwelldef}
            Under \cref{assu:SCLICQOSOSC,assu:forcingfuncomplgradLagbounded,assu:barrparamdecreaseupodr,assu:barrparammonononincr}, the point $\allvarotriterpinit$ in~\cref{def:extrapolation} is well-defined and satisfies \cref{eq:stopcondstrictfeasi} for any $\otriteridx$-th iteration of \cref{algo:RIPMOuter} satisfying that $\allvarotriter$ is sufficiently close to $\allvaraccum$ and $\barrparamotriterm$ is sufficiently small.
        \end{corollary}

        We now proceed to show that the initial point $\allvarotriterpinit$ defined by~\cref{def:extrapolation} also satisfies the remaining stopping conditions~\cref{eq:stopcondKKT,eq:stopcondbarrcompl,eq:stopcondeqvio}.
        To this end, we provide one more auxiliary lemma, which is an extension of \cite[Lemma~7.4.9]{Absiletal08OptimAlgoonMatMani}.
        \newcommand{\lemmalincombgradineq}{
        \begin{lemma}\label{lemm:lincombgradineq}
            Let $\retr[]$ be a retraction on $\mani$, and let $\funtwo[1], \ldots, \funtwo[\numfun] \in \funset\paren*{\mani}$ be continuously differentiable functions.
            Given $\ptaccum\in\mani$ and $\coeffsev > 1$, there exist a closed neighborhood $\subsetmaninin\subseteq\mani$ of $\ptaccum$ and $\tholdvalten > 0$ such that, for all $\pt \in \subsetmaninin$, any $\vecsix\in\setR[\numfun]$, and all $\tanvecone[\pt]\in\tanspc[\pt]\mani$ with $\Riemnorm[\pt]{\tanvecone[\pt]} \leq \tholdvalten$, 
            \isextendedversion{\begin{align}\label{ineq:lincombgradpullbackgrad}
                \Riemnorm[\pt]{\sum_{\idxfou=1}^{\numfun}\vecsix[\idxfou]\gradstr\funtwo[\idxfou]\paren*{\retr[\pt]\paren*{\tanvecone[\pt]}}} \leq \coeffsev\Riemnorm[\pt]{\sum_{\idxfou=1}^{\numfun}\vecsix[\idxfou]\gradstr\pullback[\pt]{\funtwo}^{\idxfou}\paren*{\tanvecone[\pt]}}.
            \end{align}}{$\Riemnorm[\pt]{\sum_{\idxfou=1}^{\numfun}\vecsix[\idxfou]\gradstr\funtwo[\idxfou]\paren*{\retr[\pt]\paren*{\tanvecone[\pt]}}} \leq \coeffsev\Riemnorm[\pt]{\sum_{\idxfou=1}^{\numfun}\vecsix[\idxfou]\gradstr\pullback[\pt]{\funtwo}^{\idxfou}\paren*{\tanvecone[\pt]}}$.}
        \end{lemma}}
        \isextendedversion{
        \lemmalincombgradineq{}
        \begin{proof}
            See \cref{appx:prooflincombgradineq}.
            \qed
        \end{proof}}{\lemmalincombgradineq{}}
        Now, we prove the proposition that the initial point $\allvarotriterpinit$ defined by~\cref{def:extrapolation} satisfies \cref{eq:stopcondKKT,eq:stopcondbarrcompl,eq:stopcondeqvio} using \cref{lemm:lincombgradineq}.
        Note that the coefficients of the following bound can be made arbitrarily small by taking $\allvarotriter$ sufficiently close to $\allvaraccum$ and $\barrparamotriterm$ sufficiently small. 
        \begin{proposition}\label{lemm:barrKKTvecflddirNewtonbound}
            Choose $\constgradLagNewton, \constcomplNewton, \consteqvioNewton \in \setRpp[]$ arbitrarily.
            Under \cref{assu:SCLICQOSOSC,assu:forcingfuncomplgradLagbounded,assu:barrparammonononincr,assu:barrparamdecreaseupodr,assu:retrsecondord,assu:pullbackHessobjineqfunLipschitz}, 
            \isextendedversion{
            \begin{align}
                &\Riemnorm[\ptotriter]{\gradstr[\pt]\Lagfun\paren*{\retr[\pt]\paren*{\dirNewtonptbarrparamotriter},  \ineqLagmultotriter[] + \dirNewtonineqLagmultbarrparamotriter[], \eqLagmultotriter[] + \dirNewtoneqLagmultbarrparamotriter[]}} \leq \constgradLagNewton\barrparamotriter,\label{ineq:gradLagfunnormdirptNewtonbound}\\
                &\norm*{\Ineqfunmat[]\paren*{\retr[\ptotriter]\paren*{\dirNewtonptbarrparamotriter}}\paren*{\ineqLagmultotriter[] + \dirNewtonineqLagmultbarrparamotriter[]}- \barrparamotriter\onevec} \leq \constcomplNewton\barrparamotriter,\label{ineq:complnormdirptNewtonbound}\\
                &\norm*{\eqfun[]\paren*{\retr[\ptotriter]\paren*{\dirNewtonptbarrparamotriter}}}\leq\consteqvioNewton\barrparamotriter\label{ineq:eqvionormdirptNewtonbound}
            \end{align}}{$\Riemnorm[\ptotriter]{\gradstr[\pt]\Lagfun\paren*{\retr[\pt]\paren*{\dirNewtonptbarrparamotriter},  \ineqLagmultotriter[] + \dirNewtonineqLagmultbarrparamotriter[], \eqLagmultotriter[] + \dirNewtoneqLagmultbarrparamotriter[]}} \leq \constgradLagNewton\barrparamotriter$, \\
            $\norm*{\Ineqfunmat[]\paren*{\retr[\ptotriter]\paren*{\dirNewtonptbarrparamotriter}}\paren*{\ineqLagmultotriter[] + \dirNewtonineqLagmultbarrparamotriter[]}- \barrparamotriter\onevec} \leq \constcomplNewton\barrparamotriter$, and $\norm*{\eqfun[]\paren*{\retr[\ptotriter]\paren*{\dirNewtonptbarrparamotriter}}}\leq\consteqvioNewton\barrparamotriter$}
            hold for any $\otriteridx$-th iteration of \cref{algo:RIPMOuter} satisfying that $\allvarotriter$ is sufficiently close to $\allvaraccum$ and $\barrparamotriterm$ is sufficiently small.
        \end{proposition}
        \begin{proof}
            See \cref{appx:proofbarrKKTvecflddirNewtonbound}.
            \qed
        \end{proof}
        Using the aforementioned results, we conclude that \cref{algo:RIPMOuter} can proceed without any inner iteration as in the following theorem.
        \newcommand{\subsetmaniele}[1][]{\mathcal{P}_{#1}}  
        \newcommand{\tholdbarrparamfiv}[1][]{\barrparam[]^{\prime}} 
        \begin{theorem}\label{coro:Newtonstopconds}
            Under \cref{assu:SCLICQOSOSC,assu:forcingfuncomplgradLagbounded,assu:barrparammonononincr,assu:barrparamdecreaseupodr,assu:retrsecondord,assu:pullbackHessobjineqfunLipschitz}, the initial point $\allvarotriterpinit$ defined by~\cref{def:extrapolation} satisfies \cref{eq:stopcondKKT,eq:stopcondbarrcompl,eq:stopcondeqvio,eq:stopcondstrictfeasi} with $\barrparamotriter > 0$ for any $\otriteridx$-th iteration satisfying that $\allvarotriter$ is sufficiently close to $\allvaraccum$ and $\barrparamotriterm$ is sufficiently small.
        \end{theorem}
        \begin{proof}
            The theorem immediately follows from \cref{coro:extrapolationwelldef,lemm:barrKKTvecflddirNewtonbound}.
            \qed
        \end{proof}

        \subsection{Local Near-Quadratic Convergence of \cref{algo:RIPMOuter}}\label{subsec:locnearquadconv}

        In this subsection, we prove the local superlinear convergence of \cref{algo:RIPMOuter} using the extrapolation.
        We also present a specific update of the sequence of barrier parameters for the local near-quadratic convergence of \cref{algo:RIPMOuter}.        
        To this end, we additionally assume the following:
        \begin{assumption}\label{assu:usedirexactsolution}
            The inputs $\allvar[0]\in\strictfeasirgn\times\setRpp[\ineqdime]\times\setR[\eqdime]$ and $\barrparam[-1] > 0$  are sufficiently close to $\allvaraccum$ and sufficiently small, respectively, and satisfy \cref{eq:stopcondKKT,eq:stopcondbarrcompl,eq:stopcondeqvio,eq:stopcondstrictfeasi}.
        \end{assumption}
        \begin{assumption}\label{assu:barrparamdecreaselowodr}
            $\barrparamotriterp = \smallO[\barrparamotriter]$ holds.
        \end{assumption}
        In the following lemma, we provide the update rule for the sequence of barrier parameters that satisfies \cref{assu:barrparamdecreaselowodr,assu:barrparamdecreaseupodr,assu:barrparammonononincr} and will be used for \cref{algo:RIPMOuter} to achieve the local near-quadratic convergence.
        The proof is straightforward.
        \newcommand{\defbarrparamupdaterule}{
        \begin{lemma}\label{def:barrparamupdaterule}
            The update rule
            \begin{align}\label{eq:barrparamupdaterule}
                \barrparamotriterp \leftarrow \coefffiv \barrparamotriter^{1+\constnin} \text{ with } 0< \coefffiv < 1, 0 < \constnin < 1, \text{ and } 0 < \barrparam[0] \leq 1
            \end{align}            
            satisfies \cref{assu:barrparamdecreaseupodr,assu:barrparammonononincr,assu:barrparamdecreaselowodr}.
        \end{lemma}}
        \isextendedversion{
        \defbarrparamupdaterule{}
        \begin{proof}
            See \cref{appx:proofbarrparamupdaterule}.
        \end{proof}}{\defbarrparamupdaterule{}}
        Let $\paren*{\opensubset, \chart[\pt]}$ be a chart of $\mani$ with $\pt\in\opensubset$, and let $\paren*{\setR[\ineqdime], \chart[{\ineqLagmult[]}]}$ and $\paren*{\setR[\eqdime], \chart[{\eqLagmult[]}]}$ be charts for $\setR[\ineqdime]$ and $\setR[\eqdime]$, respectively, where $\chart[{\ineqLagmult[]}]$ and $\chart[{\eqLagmult[]}]$ are the identity maps.
        We consider the product manifold $\mani\times\setR[\ineqdime]\times\setR[\eqdime]$ and define the chart $\chart\paren*{\allvar}=\paren*{\chart[\pt]\paren*{\pt}, \chart[{\ineqLagmult[]}]\paren*{\ineqLagmult[]}, \chart[{\eqLagmult[]}]\paren*{\eqLagmult[]}}$.
        \isextendedversion{
        Using the chart, we provide the definitions of the local convergence:
        \begin{definition}[{\hspace{-0.04em}\cite[Definition~4.5.2]{Absiletal08OptimAlgoonMatMani}}]
            Let $\brc*{\allvarotriter}_{\otriteridx}\subseteq\mani\times\setR[\ineqdime]\times\setR[\eqdime]$ be a sequence converging to $\allvaraccum$.
            Let $\paren*{\opensubset, \chart[\pt]}$ be a chart of $\mani$ with $\ptaccum\in\opensubset$ and define $\loccoord{\allvar}=\paren*{\chart[\pt]\paren*{\pt}, \ineqLagmult[], \eqLagmult[]}$ for any $\allvar\in\opensubset\times\setR[\ineqdime]\times\setR[\eqdime]$.
            If 
            \isextendedversion{
            \begin{align}
                \lim_{\otriteridx\to\infty} \frac{\norm*{\loccoord{\allvarotriterp} - \loccoord{\allvaraccum}}}{\norm*{\loccoord{\allvarotriter} - \loccoord{\allvaraccum}}} = 0,
            \end{align}}
            {$\lim_{\otriteridx\to\infty} \frac{\norm*{\loccoord{\allvarotriterp} - \loccoord{\allvaraccum}}}{\norm*{\loccoord{\allvarotriter} - \loccoord{\allvaraccum}}} = 0$,}
            then $\brc*{\allvarotriter}_{\otriteridx}$ is said to converge superlinearly to $\allvaraccum$.
            If there exist constants $p > 0, c \geq 0$, and $K \geq 0$ such that, for all $\otriteridx\geq K$, there holds
            \isextendedversion{
            \begin{align}
                \norm*{\loccoord{\allvarotriterp} - \loccoord{\allvaraccum}} \leq c \norm*{\loccoord{\allvarotriter} - \loccoord{\allvaraccum}}^{p},
            \end{align}}{$\norm*{\loccoord{\allvarotriterp} - \loccoord{\allvaraccum}} \leq \norm*{\loccoord{\allvarotriter} - \loccoord{\allvaraccum}}^{p}$,}
            then $\brc*{\allvarotriter}_{\otriteridx}$ is said to converge to $\allvaraccum$ with order at least $p$.
        \end{definition}
        We proceed to define local near-quadratic convergence, which is our main goal:
        \begin{definition}
            An algorithm is said to have local near-quadratic convergence if, for any given parameter $\constnin\in\setR$ with $0 < \constnin < 1$, a sequence $\brc*{\allvarotriter}_{\otriteridx}$ generated by the algorithm converges to $\allvaraccum$ with order at least  $1+\constnin$.
        \end{definition}}{}
        The coordinate expressions of the barrier \KKT{} vector field \isextendedversion{\cref{eq:barrKKTvecfld}}{} around $\overline{\omega}$
        is
        \begin{equation}
            \begin{aligned}[t]\label{def:loccoordbarrKKTvecfld}
                \loccoordbarrKKTvecfld\paren*{\loccoord{\allvar}, \barrparam[]} &\coloneqq 
                \begin{bmatrix}
                    \D\chart[\pt]\paren*{\inv{\chart[\pt]}\paren*{\loccoord{\pt}}}\sbra*{\gradstr[\pt]\Lagfun\paren*{\inv{\chart}\paren*{\loccoord{\allvar}}}}\\
                    \D\chart[{\ineqLagmult[]}]\paren*{\inv{\chart[{\ineqLagmult[]}]}\paren*{\loccoord{\ineqLagmult[]}}}\sbra*{\Ineqfunmat[]\paren*{\inv{\chart[\pt]}\paren*{\loccoord{\pt}}}\ineqLagmult[] - \barrparam[]\onevec}\\
                    \D\chart[{\eqLagmult[]}]\paren*{\inv{\chart[{\eqLagmult[]}]}\paren*{\loccoord{\eqLagmult[]}}}\sbra*{\eqfun[]\paren*{\inv{\chart[\pt]}\paren*{\loccoord{\pt}}}}
                \end{bmatrix} =
                \begin{bmatrix}
                    \D\chart[\pt]\paren*{\pt}\sbra*{\gradstr[\pt]\Lagfun\paren*{\allvar}}\\
                    \Ineqfunmat[]\paren*{\pt}\ineqLagmult[] - \barrparam[]\onevec\\
                    \eqfun[]\paren*{\pt}
                \end{bmatrix}.
            \end{aligned}
        \end{equation}
        Note that $\loccoordbarrKKTvecfld\paren*{\loccoord{\allvaraccum}, 0} = 0$ holds by definition.
        We also define
        \isextendedversion{
        \begin{align}\label{def:implloccoordbarrKKTvecfld}
            \implloccoordbarrKKTvecfld\paren*{\loccoord{\allvar}, \implvecone, \implvectwo, \implvecthr} \coloneqq \loccoordbarrKKTvecfld\paren*{\loccoord{\allvar}, 0} - 
            \begin{bmatrix}
                \implvecone\\
                \implvectwo\\
                \implvecthr
            \end{bmatrix}
         \end{align}}
         {$\implloccoordbarrKKTvecfld\paren*{\loccoord{\allvar}, \implvecone, \implvectwo, \implvecthr} \coloneqq \loccoordbarrKKTvecfld\paren*{\loccoord{\allvar}, 0} - \sbra*{\implvecone; \implvectwo; \implvecthr}$}
        for $\implvecone\in\setR[\dime]$, $\implvectwo\in\setR[\ineqdime]$, and $\implvecthr\in\setR[\eqdime]$.
        Given $\barrparam[] > 0$, we write $\loccoordbarrKKTvecfld[{\barrparam[]}]\paren*{\loccoord{\allvar}}$ for the restricted function $\loccoord{\allvar}\mapsto\loccoordbarrKKTvecfld\paren*{\loccoord{\allvar}, \barrparam[]}$.
        We first prove the nonsingularity of $\D\loccoordbarrKKTvecfld[{\barrparam[]}]$ at $\loccoord{\allvaraccum}$ using \cref{lemm:Jacobnonsingaccum} as well as \cite[Theorem 3.1]{LaiYoshise2024RiemIntPtMethforCstrOptimonMani}:
        
        \newcommand{\lemmloccoordJacobnonsingaccum}{\begin{lemma}\label{lemm:loccoordJacobnonsingaccum}
            Let $\barrparam[]\in\setR[]$ be an arbitrary scalar.
            Under \cref{assu:SCLICQOSOSC}, $\D\loccoordbarrKKTvecfld[{\barrparam[]}]\paren*{\loccoord{\allvaraccum}}$ is nonsingular and independent of the value of $\barrparam[]$.
        \end{lemma}}
        \isextendedversion{
        \lemmloccoordJacobnonsingaccum{}
        \begin{proof}
            See \cref{appx:proofloccoordJacobnonsingaccum}.
            \qed
        \end{proof}}{\lemmloccoordJacobnonsingaccum{}}
        Then, we prove the existence and the uniqueness of a solution of $\implloccoordbarrKKTvecfld\paren*{\cdot, \implvecone, \implvectwo, \implvecthr} = 0$ using the implicit function theorem on manifolds~\cite[Theorem~C.40]{Lee12IntrotoSmthManibook2ndedn}
        and \cref{lemm:loccoordJacobnonsingaccum}.
        We prove them in a similar manner to the Euclidean setting~\cite[Lemma 3.1]{Gouldetal2001SuperlinConvofPDIntptAlgoforNonlinProgram}:
        
        \newcommand{\lemmasolimplloccoordbarrKKTvecfld}{\begin{lemma}\label{lemm:solimplloccoordbarrKKTvecfld}
            Under \cref{assu:SCLICQOSOSC}, let $\loccoord{\allvar}\paren*{\implvecone,\implvectwo, \implvecthr}$ be a solution to $\implloccoordbarrKKTvecfld\paren*{\cdot, \implvecone, \implvectwo, \implvecthr} = 0$.
            Then, for some $\implconstone > 0$,
            \begin{enumerate}
                \item the solution $\loccoord{\allvar}\paren*{\implvecone,\implvectwo, \implvecthr}$ exists and is unique with respect to $\paren*{\implvecone,\implvectwo, \implvecthr}$ in the neighborhood $\nbhdone\paren*{\implconstone} \coloneqq \brc*{\paren*{\implvecone, \implvectwo, \implvecthr} \in \setR[\dime] \times \setR[\ineqdime] \times\setR[\eqdime] \colon \norm*{\implvecone} + \norm*{\implvectwo} + \norm*{\implvecthr} \leq \implconstone}$.
                Moreover, $\loccoord{\allvar}\paren*{\implvecone,\implvectwo, \implvecthr}$ is a continuously differentiable function of $\paren*{\implvecone, \implvectwo, \implvecthr}$ in the neighborhood $\nbhdone\paren*{\implconstone}$, and \label{lemm:solimplloccoordbarrKKTvecfldexist}
                \item for $\paren*{\implvecone[1], \implvectwo[1], \implvecthr[1]}, \paren*{\implvecone[2], \implvectwo[2], \implvecthr[2]}\in\nbhdone\paren*{\implconstone}$ sufficiently small, we have
                \isextendedversion{
                \begin{align}\label{eq:loccoordbigtheta}
                    \norm*{\loccoord{\allvar}\paren*{\implvecone[1], \implvectwo[1], \implvecthr[1]} - \loccoord{\allvar}\paren*{\implvecone[2], \implvectwo[2], \implvecthr[2]}} = \bigtheta[\norm*{\implvecone[1] - \implvecone[2]} + \norm*{\implvectwo[1] - \implvectwo[2]} + \norm*{\implvecthr[1] - \implvecthr[2]}].
                \end{align}}
                {$\norm*{\loccoord{\allvar}\paren*{\implvecone[1], \implvectwo[1], \implvecthr[1]} - \loccoord{\allvar}\paren*{\implvecone[2], \implvectwo[2], \implvecthr[2]}} = \bigtheta[\norm*{\implvecone[1] - \implvecone[2]} + \norm*{\implvectwo[1] - \implvectwo[2]} + \norm*{\implvecthr[1] - \implvecthr[2]}]$.}
                \label{lemm:solimplloccoordbarrKKTvecfldorder}
            \end{enumerate}
        \end{lemma}}
        \isextendedversion{
        \lemmasolimplloccoordbarrKKTvecfld{}
        \begin{proof}
            See \cref{appx:proofsolimplloccoordbarrKKTvecfld}.
            \qed
        \end{proof}}{\lemmasolimplloccoordbarrKKTvecfld{}}
        Next, we derive an upper bound on the Euclidean distance between $\loccoord{\allvarotriterp}$ and $\loccoord{\allvaraccum}$.
        \begin{lemma}\label{lemm:diffallvarotriterpaccumbarrparambigO}
            Suppose \cref{assu:SCLICQOSOSC,assu:forcingfuncomplgradLagbounded,assu:barrparamdecreaseupodr,assu:barrparammonononincr,assu:retrsecondord,assu:pullbackHessobjineqfunLipschitz,assu:usedirexactsolution}.
            Then,
        \isextendedversion{\begin{align}\label{ineq:diffloccoordallvarotriterpaccumbigO}
                \norm*{\loccoord{\allvarotriterp} - \loccoord{\allvaraccum}} = \bigO[\barrparamotriter].
            \end{align}}{$\norm*{\loccoord{\allvarotriterp} - \loccoord{\allvaraccum}} = \bigO[\barrparamotriter]$.}
        \end{lemma}
        \begin{proof}
            See \cref{appx:proofdiffallvarotriterpaccumbarrparambigO}.
            \qed
        \end{proof}
        Using \cref{lemm:diffallvarotriterpaccumbarrparambigO}, we derive the lower bound on the Euclidean distance between $\loccoord{\allvarotriterp}$ and $\loccoord{\allvaraccum}$, and provide the tight bound.
        \begin{lemma}\label{lemm:diffallvarotriterpaccumbarrparambigtheta}
            Suppose \cref{assu:SCLICQOSOSC,assu:forcingfuncomplgradLagbounded,assu:barrparamdecreaseupodr,assu:barrparammonononincr,assu:retrsecondord,assu:pullbackHessobjineqfunLipschitz,assu:usedirexactsolution}.
            Then, 
            \isextendedversion{
            \begin{align}
                \norm*{\loccoord{\allvarotriterp} - \loccoord{\allvaraccum}} = \bigtheta[\barrparamotriter].
            \end{align}}
            {$\norm*{\loccoord{\allvarotriterp} - \loccoord{\allvaraccum}} = \bigtheta[\barrparamotriter]$.}
        \end{lemma}
        \begin{proof}
            See \cref{appx:proofdiffallvarotriterpaccumbarrparambigtheta}.
            \qed
        \end{proof}
        Now, we establish the local convergence of a class of Riemannian \IPM{}s.
        We first prove its local superlinear convergence and then local near-quadratic convergence when using the update rule \cref{eq:barrparamupdaterule}.
        \isextendedversion{}{We refer readers to \cite[Definition~4.5.2]{Absiletal08OptimAlgoonMatMani} for the definition of local convergence in Riemannian optimization.
        }
        \begin{theorem}\label{theo:locnearquadconv}
            Suppose \cref{assu:SCLICQOSOSC,assu:forcingfuncomplgradLagbounded,assu:barrparamdecreaseupodr,assu:barrparammonononincr,assu:retrsecondord,assu:pullbackHessobjineqfunLipschitz,assu:usedirexactsolution,assu:barrparamdecreaselowodr}.
            Then, the sequence $\brc*{\allvarotriter}_{\otriteridx}$ converges to $\allvaraccum$, and the convergence rate is superlinear.
            Moreover, the sequence $\brc*{\allvarotriter}_{\otriteridx}$ converges near-quadratically
            \isextendedversion{\footnote{Although it is theoretically ideal for $\constnin$ to approach $1$ as closely as possible, it may lead to an overly rapid reduction of the barrier parameters. Since this may impair computational performance, careful tuning of $\constnin$ is important in practice.}}{\footnote{
            We say that the generated sequence $\brc*{\allvarotriter}_{\otriteridx}$ converges near-quadratically if, for any given $\constnin\in\setR$ with $0 < \constnin < 1$, the generated sequence $\brc*{\allvarotriter}_{\otriteridx}$ converges to $\allvaraccum$ with order at least $1+\constnin$.
            We refer \cite[Definition~4.5.2]{Absiletal08OptimAlgoonMatMani} for the definition of the order of convergence.
            Although it is theoretically ideal for $\constnin$ to approach $1$ as closely as possible, it may lead to an overly rapid reduction of the barrier parameters. Since this may impair computational performance, careful tuning of $\constnin$ is important in practice.}}
            if the sequence of the barrier parameters $\brc*{\barrparamotriter}_{\otriteridx}$ is updated according to \cref{eq:barrparamupdaterule}.
        \end{theorem}   
        \begin{proof}
            See \cref{appx:prooflocnearquadconv}.
            \qed
        \end{proof}

    \section{Application to the \RIPTRM{} for \RNLO{}~\cref{prob:RNLO} with $\eqset=\emptyset$}\label{sec:applicationtoRIPTRM}
        In this section, we apply our results to the \RIPTRM{}~\cite{Obaraetal2025APrimalDualIPTRMfor2ndOrdStnryPtofRiemIneqCstrOptimProbs} designed for \RNLO{}~\cref{prob:RNLO} with $\eqset=\emptyset$.
        Although the algorithm can be configured for the global convergence to either a first- or second-order stationary point (\SOSP{}) as described in \cite{Obaraetal2025APrimalDualIPTRMfor2ndOrdStnryPtofRiemIneqCstrOptimProbs}, we focus on the latter case since it theoretically subsumes the former.
        To compute an \SOSP{}, at the $\otriteridx$-th outer iteration, \RIPTRM{} sets $\barrparamotriter > 0$ so that $\lim_{\otriteridx \to \infty} \barrparamotriter = 0$ and finds a point that satisfies \cref{eq:stopcondKKT}, \cref{eq:stopcondbarrcompl}, \cref{eq:stopcondstrictfeasi}, and the condition for the second-order stationarity defined as
        \begin{align}
            \mineigval\sbra*{\trsquad\paren*{\ptotriterp, \ineqLagmultotriterp[]}} \geq - \forcingfunsosp\paren*{\barrparam},\label{eq:stopcondsecondord}
        \end{align}
        where $\forcingfunsosp\colon\setRp[]\to\setRp[]$ is a forcing function and 
        \begin{align}
            &\trsquad\paren*{\pt[], \ineqLagmult[]} \coloneqq \Hess\objfun\paren*{\pt[]} - \sum_{\ineqidx\in\ineqset}\ineqLagmult[\ineqidx]\Hess\ineqfun[\ineqidx]\paren*{\pt[]} + \ineqgradopr[\pt]\IneqLagmultmat[]\inv{\Ineqfunmat[]\paren*{\pt}}\coineqgradopr[\pt].\label{eq:trsquaddef}
        \end{align}
        To find such a point, the algorithm uses a trust region approach in the inner iteration: 
        let $\trradius[] > 0$ be the trust region radius.
        We then introduce the trust region subproblem at $\pt\in\strictfeasirgn$ with $\ineqLagmult[]\in\setRpp[\ineqdime]$,  $\trradius[] > 0$, and $\barrparam[] > 0$ as 
        \isextendedversion{
        \begin{mini}
            {\dir[] \in \tanspc[\pt]\mani}{
            \frac{1}{2}\metr[\pt]{\trsquad\paren*{\pt[], \ineqLagmult[]}\sbra*{\dir[]}}{\dir[]} + \metr[\pt]{\trslin[{\barrparam[]}]\paren*{\pt}}{\dir[]}}
            {\label{prob:TRS}}{}
            \addConstraint{\Riemnorm[\pt]{\dir[]}}{\leq \trradius[],\label{cstr:trradius}}
        \end{mini}}
        {\begin{equation}\label{prob:TRS}
            \minimize[{\dir[] \in \tanspc[\pt]\mani}] ~ \frac{1}{2}\metr[\pt]{\trsquad\paren*{\pt[], \ineqLagmult[]}\sbra*{\dir[]}}{\dir[]} + \metr[\pt]{\trslin[{\barrparam[]}]\paren*{\pt}}{\dir[]} \quad \subjectto \quad \Riemnorm[\pt]{\dir[]} \leq \trradius[]
        \end{equation}}
        where
        \isextendedversion{
        \begin{align}\label{eq:trslindef}
            \trslin[{\barrparam[]}]\paren*{\pt}\coloneqq\gradstr[]\objfun\paren*{\pt} - \barrparam[]\ineqgradopr[\pt]\sbra*{\inv{\Ineqfunmat[]\paren*{\pt}}\onevec}.
        \end{align}}
        {$\trslin[{\barrparam[]}]\paren*{\pt}\coloneqq\gradstr[]\objfun\paren*{\pt} - \barrparam[]\ineqgradopr[\pt]\sbra*{\inv{\Ineqfunmat[]\paren*{\pt}}\onevec}$.}
        Note that, since the feasible region is bounded and closed, and the objective function is continuous, the subproblem~\cref{prob:TRS} admits a global optimum.
        Despite its nonconvexity, we can compute this global optimum~\cite{Adachietal2017SolvingTRSbyGenEigenProb,CarmonDuchi20201stOrdMethforNonconvQuadMin}.
        At the $\inriteridx$-th inner iteration, \RIPTRM{} computes the search direction for $\pt$
        by solving \cref{prob:TRS} at $\ptinriter$ with $\ineqLagmultinriter[]$, $\trradiusinriter$, and $\barrparam[]$. 
        In this article, we use the global optimum of \cref{prob:TRS} as the search direction and call it the \exactsolution{}.
        \isextendedversion{
        The following proposition provides the necessary and sufficient conditions for the \exactsolution{}:
        \begin{proposition}[{\hspace{-0.04em}\cite[Proposition 7.3.1]{Absiletal08OptimAlgoonMatMani}}]\label{prop:trsgloboptimiffcond}
            The vector $\direxact[\pt]\in\tanspc[\pt]\mani$ is a global optimum of \cref{prob:TRS}
            if and only if there exists a scalar $\trsLagmult \geq 0$ such that             \begin{subequations}\label{eq:trsgloboptimiffcond}
                \begin{align}
                    \paren*{\trsquad\paren*{\pt, \ineqLagmult[]}+\trsLagmult\id[{\tanspc[\pt]\mani}]}\direxact[\pt] &= - \trslin[{\barrparam[]}]\paren*{\pt},\label{eq:trsgloboptimiffcondquadlineq}\\
                    \trsLagmult\paren*{\trradius[] - \Riemnorm[\pt]{\direxact[\pt]}} &= 0,\label{eq:trsgloboptimiffcondconmpl}\\
                    \Riemnorm[\pt]{\direxact[\pt]} &\leq \trradius[],\label{eq:trsgloboptimiffconddirtrradius}\\
                    \trsquad\paren*{\pt, \ineqLagmult[]} + \trsLagmult\id[{\tanspc[\pt]\mani}] &\succeq 0.\label{eq:trsgloboptimiffcondquadpsd}
                \end{align}
            \end{subequations}
        \end{proposition}}{}
        Using the \exactsolution{}, we also define the corresponding step for $\ineqLagmult[]$ as
        \begin{align}\label{def:direxactineqLagmult}
            \direxactineqLagmult[] = - \ineqLagmult[] + \barrparam[]\inv{\Ineqfunmat[]\paren*{\pt}}\onevec - \IneqLagmultmat[]\inv{\Ineqfunmat[]\paren*{\pt}}\coineqgradopr[\pt]\sbra*{\direxactpt}.
        \end{align}
        \RIPTRM{} first checks whether the point $\paren*{\retr[\ptinriter]\paren*{\direxactptinriter}, \ineqLagmultinriter[] + \direxactineqLagmultinriter[]}$ satisfies the stopping conditions \cref{eq:stopcondKKT,eq:stopcondbarrcompl,eq:stopcondstrictfeasi,eq:stopcondsecondord} and returns the point if it does.
        Otherwise, it updates the next iterate using these steps.
        We refer readers to \cite{Obaraetal2025APrimalDualIPTRMfor2ndOrdStnryPtofRiemIneqCstrOptimProbs} for details.
        We summarize the method in \cref{algo:RIPTRMOuter}.

        \begin{algorithm2e}[t]
            \SetKwInOut{Require}{Require}
            \SetKwInOut{Input}{Input}
            \SetKwInOut{Output}{Output}
            \SetKwProg{Fn}{Iteration}{:}{end}
            \Require{Riemannian manifold $\mani$, twice continuously differentiable functions $\objfun:\mani\to\setR[]$ and $\brc*{\ineqfun[\ineqidx]}_{\ineqidx\in\ineqset}\colon\mani\to\setR[]$, maximal trust region radius $\maxtrradius > 0$, initial trust region radius $\inittrradius_{0} \in \left(0,\maxtrradius\right]$, 
            minimum initial trust region radius $\mininittrradius\in \left(0,\maxtrradius\right]$, 
            forcing functions $\forcingfungradLag,  \forcingfuncompl, \forcingfunsosp\colon\setRp[]\to\setRp[]$.}
            \Input{Initial point $\paren*{\pt[0], \ineqLagmult[0]} \in \strictfeasirgn \times \setRpp[\ineqdime]$,
            initial barrier parameter $\barrparam[-1] > 0$.}
            \hrule
            \Fn{OUTER}{
                \For{$\otriteridx = 0, 1, \ldots$}{
                    Set $\barrparamotriter > 0$ so that $\lim_{\otriteridx \to \infty} \barrparamotriter = 0$.    
                    
                    Compute $\paren*{\ptotriterp, \ineqLagmultotriterp[]}\in \strictfeasirgn \times \setRpp[\ineqdime]$ satisfying \cref{eq:stopcondKKT,eq:stopcondbarrcompl,eq:stopcondstrictfeasi,eq:stopcondsecondord} and $\finaltrradius_{\otriteridx} > 0$ by calling INNER$\paren*{\ptotriter, \ineqLagmultotriter[], \barrparamotriter, \inittrradius_{\otriteridx}}$.

                    Set $\inittrradius_{\otriteridx+1}\gets\max\paren*{\finaltrradius_{\otriteridx}, \mininittrradius}$.
                }
            }
            \hrule
            \Fn{INNER($\pt^{0}, \ineqLagmult[]^{0}, \barrparam[]^{0}, \trradius[]^{0}$)}{
                \For{$\inriteridx = 0, 1, \ldots$\label{algo:forline}}{
                    Compute $\dirptinriter \in \tanspc[\ptinriter]\mani$ by exactly solving the subproblem~\cref{prob:TRS}.
                    
                    Compute $\dirineqLagmultinriter\in\tanspc[{\ineqLagmultinriter[]}]\setR[\ineqdime]$ according to \cref{def:direxactineqLagmult}.
                    
                    \If{$\paren*{\retr[\ptinriter]\paren*{\dirptinriter}, \ineqLagmultinriter[] + \dirineqLagmultinriter[]}$ satisfies \cref{eq:stopcondKKT,eq:stopcondbarrcompl,eq:stopcondstrictfeasi,eq:stopcondsecondord}}
                    {\Return $\paren*{\retr[\ptinriter]\paren*{\dirptinriter}, \ineqLagmultinriter[] + \dirineqLagmultinriter[]}$ and $\finaltrradius=\trradiusinriter$.}                    
                Update $\paren*{\ptinriterp, \ineqLagmultinriterp[]}\in\strictfeasirgn\times\setRpp[\ineqdime]$, and $\trradiusinriterp > 0$ using $\dirptinriter$ and $\dirineqLagmultinriter[]$. 
                }
            }
            \caption{\RIPTRM{} for finding an \SOSP{} of \RNLO{} with $\eqset=\emptyset$}\label{algo:RIPTRMOuter}
        \end{algorithm2e}

    \subsection{Local Near-Quadratic Convergence of \cref{algo:RIPTRMOuter}}
        In this subsection, we analyze the local superlinear and near-quadratic convergence of \cref{algo:RIPTRMOuter} using the preceding results.
        We additionally assume the following:
        \begin{assumption}\label{assu:HessobjineqfunLipschitz}
            There exist positive scalars $\Lipschitzconstfou[\objfun], \brc*{\Lipschitzconstfou[{\ineqfun[\ineqidx]}]}_{\ineqidx\in\ineqset} \in \setRpp[]$ such that, for all $\pt[]\in\mani$ satisfying $\Riemdist{\pt[]}{\ptaccum} < \injradius\paren*{\ptaccum}$, 
            \isextendedversion{
            \begin{align}
                &\opnorm{\Hess\objfun\paren*{\pt} - \partxp[]{\pt}{\ptaccum}\circ\Hess\objfun\paren*{\ptaccum}\circ\partxp[]{\ptaccum}{\pt}} \leq \Lipschitzconstfou[\objfun] \Riemdist{\pt[]}{\ptaccum},\label{ineq:HessobjfunLipschitzbound} \\
                &\opnorm{\Hess\ineqfun[\ineqidx]\paren*{\pt} - \partxp[]{\pt}{\ptaccum}\circ\Hess\ineqfun[\ineqidx]\paren*{\ptaccum}\circ\partxp[]{\ptaccum}{\pt}} \leq \Lipschitzconstfou[{\ineqfun[\ineqidx]}] \Riemdist{\pt[]}{\ptaccum} \text{ for all } \ineqidx\in\ineqset.\label{ineq:HessineqfunLipschitzbound}
            \end{align}}{
            $\opnorm{\Hess\objfun\paren*{\pt} - \partxp[]{\pt}{\ptaccum}\circ\Hess\objfun\paren*{\ptaccum}\circ\partxp[]{\ptaccum}{\pt}} \leq \Lipschitzconstfou[\objfun] \Riemdist{\pt[]}{\ptaccum}$, $\opnorm{\Hess\ineqfun[\ineqidx]\paren*{\pt} - \partxp[]{\pt}{\ptaccum}\circ\Hess\ineqfun[\ineqidx]\paren*{\ptaccum}\circ\partxp[]{\ptaccum}{\pt}} \leq \Lipschitzconstfou[{\ineqfun[\ineqidx]}] \Riemdist{\pt[]}{\ptaccum}$ for all $\ineqidx\in\ineqset$.}
        \end{assumption}
        \cref{assu:HessobjineqfunLipschitz} holds if $\objfun$ and $\brc*{\ineqfun[\ineqidx]}_{\ineqidx\in\ineqset}$ are of class $C^{3}$; for details, see the arXiv version of \cite{Obaraetal2025APrimalDualIPTRMfor2ndOrdStnryPtofRiemIneqCstrOptimProbs}.

        The following lemma shows that $\trsquad\paren*{\pt, \ineqLagmult[]}$ is positive definite around $\paren*{\ptaccum, \ineqLagmultaccum[]}$\isextendedversion{; see \cref{eq:trsquaddef} for the definition of $\trsquad$.}{.}
        In the proof, we use the parallel transport along the minimizing geodesic since the operator $\trsquad\paren*{\pt, \ineqLagmult[]}$ is defined on the space $\tanspc[\pt]\mani$, which varies depending on $\pt\in\mani$.
        \begin{lemma}\label{lemm:trsquadposdifardallvaraccum}
            Under \cref{assu:SCLICQOSOSC,assu:HessobjineqfunLipschitz}, $\trsquad\paren*{\pt, \ineqLagmult[]}$ is positive definite for all $\paren*{\pt, \ineqLagmult[]}\in\strictfeasirgn\times\setRpp[\ineqdime]$ sufficiently close to $\paren*{\ptaccum, \ineqLagmultaccum[]}\in\feasirgn\times\setRp[\ineqdime]$.
        \end{lemma}
        \begin{proof}
            See \cref{appx:prooftrsquadposdifardallvaraccum}.
            \qed
        \end{proof}

        We next prove that the \exactsolution{} is identical to the Newton step~\cref{def:Newtondirection} under the given assumptions.
        Here, $\Jacobian[\barrKKTvecfld]\paren*{\pt, \ineqLagmult[]}$ and $\paren*{\dirNewtonptbarrparam, \dirNewtonineqLagmultbarrparam[]}\in\tanspc[\pt]\mani\times\setR[\ineqdime]$ denote the Jacobian \cref{def:JacobbarrKKTvecfld} and the Newton step \cref{def:Newtondirection} for the purely inequality-constrained setting, respectively.
        \isextendedversion{}{The proposition 
        is based on the necessary and sufficient conditions for \exactsolution{}~\cite[Proposition 7.3.1]{Absiletal08OptimAlgoonMatMani}; see \cite[Appendix O]{Obaraetal2025LocalConvofRiemIPMs} for the complete proof.}
        
        \newcommand{\propdirNewtonardaccum}{\begin{proposition}\label{prop:dirNewtonardaccum}
            Let $\trradius[], \barrparam[] \in \setRpp[]$ be any positive scalars, and let $\paren*{\pt, \ineqLagmult[]}\in\feasirgn\times\setRp[\ineqdime]$ be any point for which $\Jacobian[\barrKKTvecfld]\paren*{\pt, \ineqLagmult[]}$ is nonsingular. 
            Let $\direxactpt\in\tanspc[\pt]\mani$ denote the global optimum of \cref{prob:TRS} at $\pt\in\mani$ with $\trradius[]$ and $\barrparam[]$, and let $\direxactineqLagmult\in\setR[\ineqdime]$ be the associated vector \cref{def:direxactineqLagmult}.
            Under \cref{assu:SCLICQOSOSC}, if the search direction $\direxactpt\in\tanspc[\pt]\mani$ satisfies $\Riemnorm[\pt]{\direxactpt} < \trradius[]$, 
            then we have $\paren*{\direxactpt, \direxactineqLagmult[]} = \paren*{\dirNewtonptbarrparam, \dirNewtonineqLagmultbarrparam[]}$.
        \end{proposition}}
        \isextendedversion{
        \propdirNewtonardaccum{}
        \begin{proof}
            See \cref{appx:proofdirNewtonardaccum}.
            \qed
        \end{proof}}{\propdirNewtonardaccum{}}
        Using the continuity of the Jacobian and the barrier \KKT{} vector field, we show that the Newton step can be made arbitrarily small in a neighborhood near $\allvaraccum$:
        
        \newcommand{\lemmdirNewtonwelldefinedbound}{\begin{lemma}\label{lemm:dirNewtonwelldefinedbound}
            Choose $\trradius[] > 0$ arbitrarily. 
            Then, $\Riemnorm[{\paren*{\pt, \ineqLagmult[]}}]{\paren*{\dirNewtonptbarrparam, \dirNewtonineqLagmultbarrparam[]}} \leq \trradius[]$ holds for any $\paren*{\pt, \ineqLagmult[]}\in\feasirgn\times\setRp[\ineqdime]$ sufficiently close to $\paren*{\ptaccum, \ineqLagmultaccum[]}$ and any $\barrparam[] > 0$ sufficiently small.
        \end{lemma}}
        \isextendedversion{
        \lemmdirNewtonwelldefinedbound{}
        \begin{proof}
            See \cref{appx:proofdirNewtonwelldefinedbound}.
            \qed
        \end{proof}}{\lemmdirNewtonwelldefinedbound{}}
        For the $\otriteridx$-th OUTER iteration of \cref{algo:RIPTRMOuter}, the algorithm calls the INNER iteration with the iterate $\paren*{\ptotriter, \ineqLagmultotriter[]}\in\mani\times\setR[\ineqdime]$ and the radius $\inittrradius_{\otriteridx+1} \geq \mininittrradius > 0$ set as the initial point and the initial trust region radius, respectively.
        Under \cref{assu:SCLICQOSOSC}, it follows from \cref{prop:dirNewtonardaccum,lemm:dirNewtonwelldefinedbound} that, at the first iteration of INNER($\pt^{\otriteridx}, \ineqLagmult[]^{\otriteridx}, \barrparam[]^{\otriteridx}, \trradius[]^{\otriteridx}$), the search direction is equivalent to the Newton step if $\paren*{\ptotriter, \ineqLagmultotriter[]}$ is sufficiently close to $\paren*{\ptaccum, \ineqLagmultaccum[]}\in\feasirgn\times\setRp[\ineqdime]$ and $\barrparamotriter > 0$ is sufficiently small.
        Combining this with \cref{lemm:dirNewtonineqfunpos,lemm:dirNewtonineqLagmultpos,lemm:barrKKTvecflddirNewtonbound,lemm:trsquadposdifardallvaraccum}
        implies that the first iterate $\paren*{\retr[{\pt[0]}]\paren*{\dir[\pt^{0}]}, \ineqLagmult[]^{0} + \dir[{\ineqLagmult[]^{0}}]}
        $ of the INNER iteration satisfies the stopping conditions \cref{eq:stopcondKKT,eq:stopcondbarrcompl,eq:stopcondstrictfeasi,eq:stopcondsecondord} with $\barrparamotriter > 0$, provided that $\paren*{\ptotriter, \ineqLagmultotriter[]}$ is sufficiently close to $\paren*{\ptaccum, \ineqLagmultaccum[]}\in\feasirgn\times\setRp[\ineqdime]$ and $\barrparamotriter > 0$ is sufficiently small.
        Consequently, the INNER iteration terminates at the first iteration.
        The following corollary formally states the result:
        \begin{corollary}\label{coro:RIPTRMoneNewtonstep}
            Under \cref{assu:SCLICQOSOSC,assu:forcingfuncomplgradLagbounded,assu:barrparammonononincr,assu:barrparamdecreaseupodr,assu:retrsecondord,assu:pullbackHessobjineqfunLipschitz,assu:HessobjineqfunLipschitz}, 
            $\paren*{\ptotriterp, \ineqLagmultotriterp[]} = $ \\ $\paren*{\retr[{\ptotriter}]\paren*{\dir[\ptotriter]},  \ineqLagmultotriter[] + \dir[{\ineqLagmultotriter[]}]}$ holds
            for any $\otriteridx$-th OUTER iteration of \cref{algo:RIPTRMOuter} satisfying that $\paren*{\ptotriter, \ineqLagmultotriter[]}$ is sufficiently close to $\paren*{\ptaccum, \ineqLagmultaccum[]}$ and $\barrparamotriterm$ is sufficiently small.
        \end{corollary}
        Therefore, using the same argument as in \cref{subsec:locnearquadconv}, we obtain the following convergence property:
        \begin{corollary}\label{coro:RIPTRMlocalconv}
            Suppose \cref{assu:SCLICQOSOSC,assu:forcingfuncomplgradLagbounded,assu:barrparamdecreaseupodr,assu:barrparammonononincr,assu:pullbackHessobjineqfunLipschitz,assu:usedirexactsolution,assu:barrparamdecreaselowodr,assu:retrsecondord,assu:HessobjineqfunLipschitz}.
            Then, the sequence $\brc*{\paren*{\ptotriter, \ineqLagmultotriter[]}}_{\otriteridx}$ generated by \cref{algo:RIPTRMOuter} converges to $\paren*{\ptaccum, \ineqLagmultaccum[]}$, and the convergence rate is superlinear.
            Moreover, the sequence $\brc*{\paren*{\ptotriter, \ineqLagmultotriter[]}}_{\otriteridx}$ converges near-quadratically if the sequence of the barrier parameters $\brc*{\barrparamotriter}_{\otriteridx}$ is updated according to \cref{eq:barrparamupdaterule}.
        \end{corollary}

        Let us conclude this section by highlighting the transition from global to local convergence of \cref{algo:RIPTRMOuter}. 
        \Cref{algo:RIPTRMOuter} possesses global convergence, as confirmed in \cite{Obaraetal2025APrimalDualIPTRMfor2ndOrdStnryPtofRiemIneqCstrOptimProbs}. 
        Accordingly, it generates iterates based on the global convergence property when the iterates are far from a solution and switches to updates based on local convergence once they are sufficiently close to the solution. 
        A favorable feature of \Cref{algo:RIPTRMOuter} is that this transition occurs smoothly and eliminates the need for additional computations to shift from the global to the local regime. 
        This property arises from employing the global optimum of \cref{prob:TRS} as the search direction.

    \section{Numerical experiments}\label{sec:experiments}

        In this section, we present numerical experiments.
        We first illustrate the behavior of our methods under different settings for updating the barrier parameter.
        We then solve another instance using our algorithms and existing methods, including two Riemannian interior point methods proposed by Lai and Yoshise~\cite{LaiYoshise2024RiemIntPtMethforCstrOptimonMani}.
        All the experiments were implemented in Python using Pymanopt~2.2.1~\cite{TownsendKoepWeichwald2017PymanoptPyToolboxforOptimonManiusingAutoDiff} and executed on a MacBook Pro with an Apple M1 Max chip and 64~GB of memory.

    \subsection{Problem setting}

        We construct a synthetic optimization problem on the set of positive definite matrices.
        Let $A\in\setR[\dime\times\dime]$ be a target matrix, and define $\symposmani\coloneqq\brc*{\matone\in\setR[\dime\times\dime]\colon\matone\succ 0}$ as the set of $\dime$-dimensional positive definite matrices.
        We approximate $\A$ with a nonnegative positive definite matrix $\matone\in\setR[\dime\times\dime]$ by solving the following optimization problem:
        \begin{align}
            \minimize[\matone \in \symposmani]~\frac{1}{2}\fronorm{\A - \matone}^{2} \quad \subjectto~\matone \geq 0.\label{prob:syntheticPD}
        \end{align}

        \textbf{Input.} 
        We consider the case $\dime = 5$.
        The target matrix $\A$ is generated as follows: first, we generate $S\in\setR[\dime\times\dime]$ whose entries are drawn from the uniform distribution on $[0,1)$.
        We then compute $\trsp{S}S$ and confirm its positive definiteness via eigenvalue decomposition.  
        Gaussian noise with mean~$0$ and standard deviation~$0.01$ is added to each entry to obtain the matrix $\A$.
        We also generate an initial point $\matone[0]$ by sampling $S_{0}\in\setR[\dime\times\dime]$ with entries sampled from the uniform distribution on $[0,1)$, computing $X_{0} = \trsp{S_{0}}S_{0}$, and confirming its positive definiteness.
        The initial point lies in the interior of the feasible set.
        In the first experiment, both $\A$ and $\matone[0]$ are generated randomly.
        In the second, we generate another $\A$ together with $20$ different initial points.

    \subsection{Experimental environment}
        We compare the following algorithms:
        \begin{itemize}
            \item \LERIPM{}: the Riemannian interior point method with extrapolation steps (\cref{algo:RIPMOuter} without inner iterations).
            \item \RIPTRM{}: the Riemannian interior point trust region method (\cref{algo:RIPTRMOuter}).
            \item \GlobalRIPM{}: a globally convergent Riemannian interior point method~\cite[Algorithm~5]{LaiYoshise2024RiemIntPtMethforCstrOptimonMani}.
            \item \LocalRIPM{}: a locally convergent Riemannian interior point method~\cite[Algorithm~2]{LaiYoshise2024RiemIntPtMethforCstrOptimonMani}.
            \item \RALM{}: a Riemannian augmented Lagrangian method~\cite[Algorithm 1]{LiuBoumal2020SimpleAlgoforOptimonRiemManiwithCstr}.
            \item \RSQO{}: a Riemannian sequential quadratic optimization method~\cite[Algorithm 3.1]{ObaraOkunoTakeda2021SQOforNLOonRiemMani}.
        \end{itemize}
        In the implementation of \LERIPM{}, we employ the conjugate residual method, a matrix-free algorithm for solving linear systems, similar to \cite{LaiYoshise2024RiemIntPtMethforCstrOptimonMani}. 
        We use \cref{eq:barrparamupdaterule} for the update of the barrier parameter.
        The implementation of \RIPTRM{} follows that in a companion paper~\cite{Obaraetal2025APrimalDualIPTRMfor2ndOrdStnryPtofRiemIneqCstrOptimProbs} except for $\constnin$ used in \cref{eq:barrparamupdaterule}.
        The implementations of \GlobalRIPM{}, \LocalRIPM{}, \RALM{}, and \RSQO{} follow the MATLAB codes provided by the original authors~\cite{LiuBoumal2019CodeOptimonManiwithExtraCstrs,Obara2021CodeSeqQuadProgramonMani,Lai2024CodeRIPM}.
        We adopt the default parameter settings for them, except that we use gradient descent as the subsolver for \RALM{} in place of the limited-memory BFGS employed in the original MATLAB implementation, since the latter is not available in Pymanopt~2.2.1.

        To measure the deviation of an iterate from the set of \KKT{} points, we introduce the following residual for \cref{prob:syntheticPD}: 
        {\small 
        \begin{align}\label{def:KKTresid}            
            &\sqrt{\Riemnorm[\matone]{\gradstr[\pt]\Lagfun\paren*{\matone, \ineqLagmult[]}}^{2} + \sum_{1\leq\matrowidx, \matcolidx\leq\dime}\paren*{\min\paren*{0, \ineqLagmult[\matrowidx\matcolidx]}^{2} + \min\paren*{0, \matone[\matrowidx\matcolidx]}^{2} + \paren*{\ineqLagmult[\matrowidx\matcolidx]\matone[\matrowidx\matcolidx]}^{2} } + \Manvio\paren*{\matone}^{2}},
        \end{align}}
        where the first two terms correspond to the \KKT{} conditions, and $\Manvio\colon\mani\to\setR[]$ denotes the violation of the manifold constraints; specifically, $\Manvio\paren*{\matone}=+\infty$ if the matrix $\matone$ has a negative eigenvalue, and $0$ otherwise.

    \subsection{Results and discussion}
\begin{figure}[t]
            \centering
            \includegraphics[width=\linewidth]{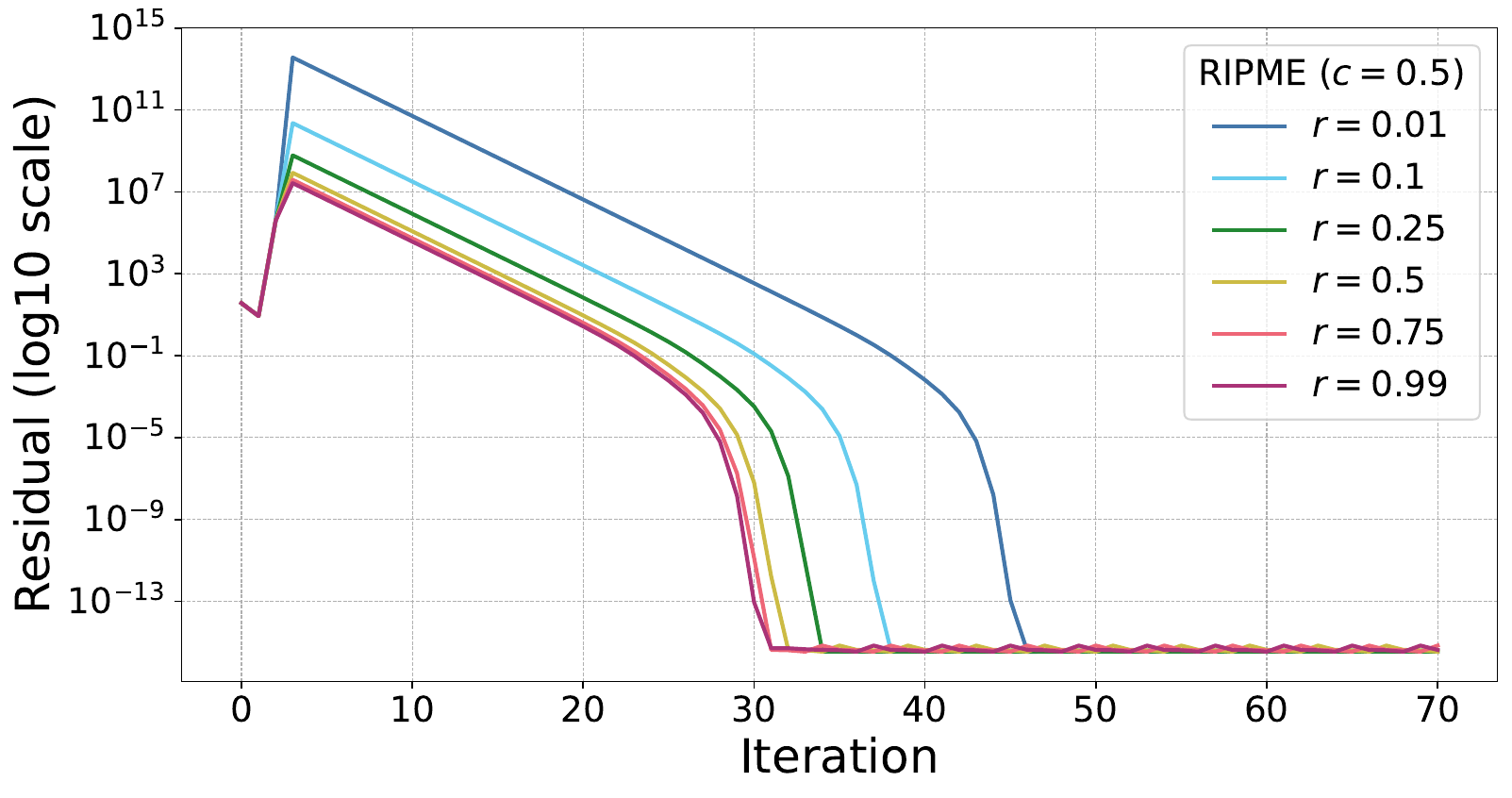}
            \caption{Residuals over iterations of \LERIPM{} for different values of $\constnin$.}
            \label{fig:residiterripm}
        \end{figure}

\begin{figure}[t]
            \centering
            \includegraphics[width=\linewidth]{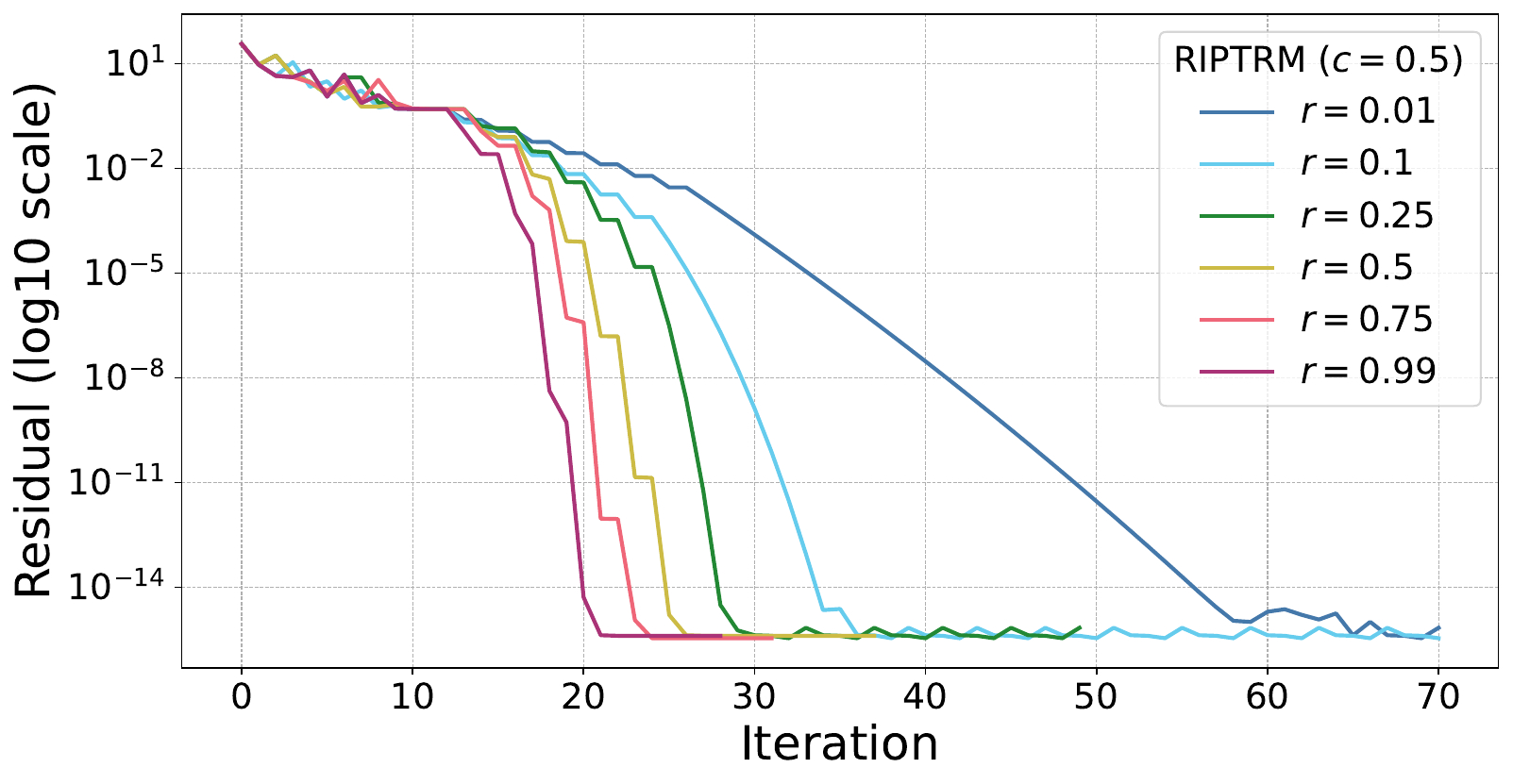}
            \caption{Residuals over iterations of \RIPTRM{} for different values of $\constnin$.}
            \label{fig:residiterriptrm}
        \end{figure}

        We first applied our proposed methods, \LERIPM{} and \RIPTRM{}, to a randomly generated instance of \cref{prob:syntheticPD}.
        For \LERIPM{} and \RIPTRM{},
        we varied the parameter $\constnin \in \{0.01, 0.1, 0.25, 0.5, 0.75, 0.99\}$.
        \Cref{fig:residiterripm} shows the residuals over iterations for \LERIPM{} under different values of $\constnin$.
        For any $\constnin$, once the barrier parameter decreased to the level of numerical error ($10^{-16}$), \LERIPM{} exhibited rapid local convergence to a solution.
        \Cref{fig:residiterriptrm} presents the corresponding results for \RIPTRM{}. 
        In all cases, we observed that \RIPTRM{} seamlessly transitioned from the global to the local regime. 
        \RIPTRM{} demonstrated stable local convergence behavior largely consistent with the theory, reducing the residual at the rate of $1+\constnin$; 
        for $\constnin \leq 0.25$, the algorithm repeated with one inner iteration followed by an update of $\barrparam$, as confirmed in \cref{coro:RIPTRMoneNewtonstep}. 
        For $r \geq 0.5$, two inner iterations and an update of $\barrparam$ were alternated, producing a stepwise decrease of the residuals. 
        This was not fully consistent with the theoretical prediction but remained broadly in line with it.
    
\begin{figure}[t]
            \centering
            \includegraphics[width=\linewidth]{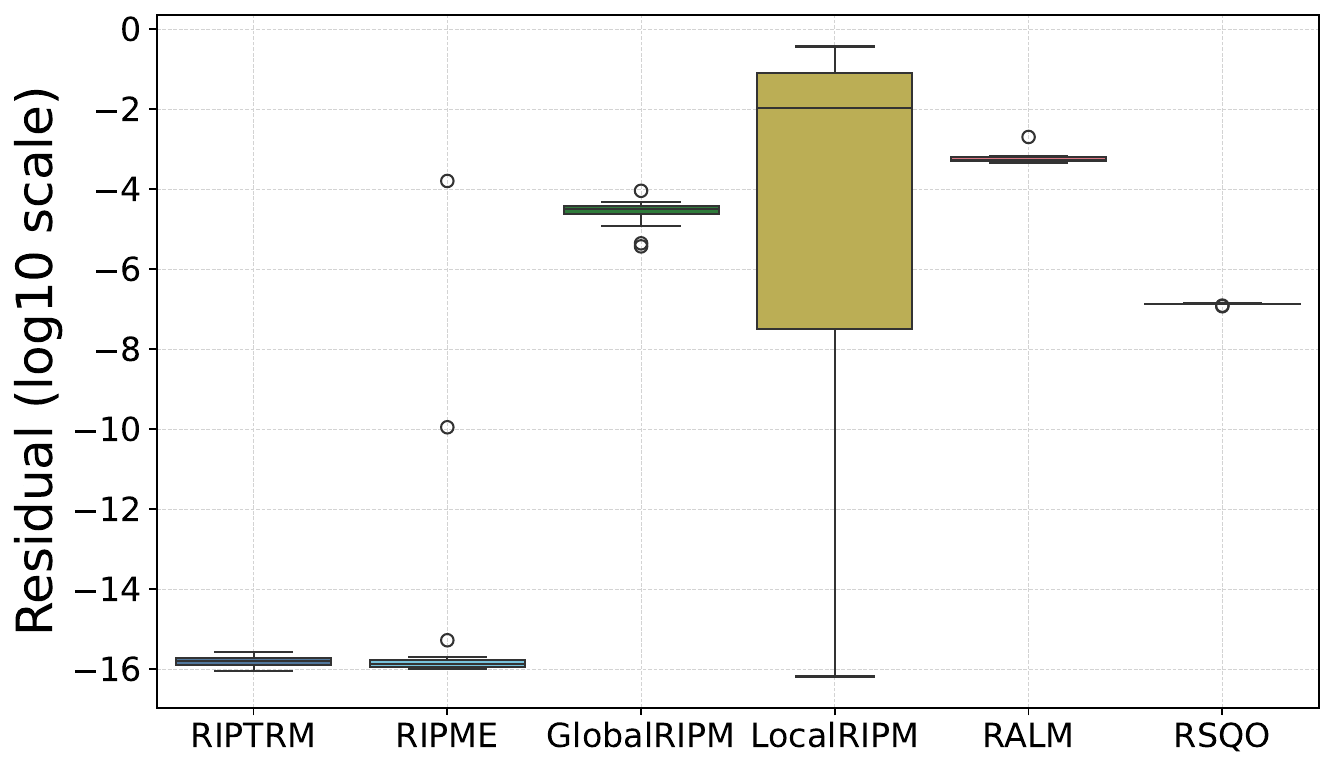}
            \caption{Box plots of the residuals across $20$ runs from different initial points.}
            \label{fig:boxplotresid}
        \end{figure}

        Next, we generated another $\A$ along with $20$ initial points and compared our methods with existing ones.  
        We set $\constnin = 0.25$ for both \LERIPM{} and \RIPTRM{}.  
        \Cref{fig:boxplotresid} presents, for each solver, box plots of the residuals from $20$ runs with different initial points, where the residual for a run is defined as the minimum value of \cref{def:KKTresid} attained during that run.  
        We observed that \LERIPM{} and \RIPTRM{} successfully solved the instances with high accuracy.

    \section{Conclusion}\label{sec:conclusion}
        \isextendedversion{
        In this paper, we considered \RNLO{}~\cref{prob:RNLO} and analyzed a class of \RIPM{}s— namely, \cref{algo:RIPMOuter}—for solving it.
        We introduced the extrapolation for computing the initial point~\cref{def:extrapolation} in each outer iteration and proved that this point satisfies the stopping conditions as stated in \cref{coro:Newtonstopconds}.
        This implies that \cref{algo:RIPMOuter} can proceed without performing inner iterations.
        In \cref{theo:locnearquadconv}, we established that the local superlinear convergence of \cref{algo:RIPMOuter}, as well as the local near-quadratic convergence under the update rule~\cref{eq:barrparamupdaterule} for the barrier parameters.
        We also applied these results to a \RIPTRM{} proposed in \cite{Obaraetal2025APrimalDualIPTRMfor2ndOrdStnryPtofRiemIneqCstrOptimProbs}, \cref{algo:RIPTRMOuter}.
        We proved that \cref{algo:RIPTRMOuter} also satisfies the stopping conditions when using the \exactsolution{}, which is the global optimum of \cref{prob:TRS}, and the second-order stationarity around a solution.
        This implies the local superlinear and near-quadratic convergence of \cref{algo:RIPTRMOuter}.
        Numerical results supported the theoretical analyses of the proposed methods.
        }{In this paper, we considered \RNLO{}~\cref{prob:RNLO} and analyzed a class of \RIPM{}s employing the extrapolation step, \cref{algo:RIPMOuter}, for solving it.
        We proved the local superlinear convergence of the algorithm as well as the local near-quadratic convergence under certain assumptions on the update rule of barrier parameters. 
        We also applied these results to the \RIPTRM{} proposed in \cite{Obaraetal2025APrimalDualIPTRMfor2ndOrdStnryPtofRiemIneqCstrOptimProbs}, \cref{algo:RIPTRMOuter}, to establish its local superlinear and near-quadratic convergence.
        Numerical results supported the theoretical analyses of the proposed methods.
        }
        In closing, we discuss future directions for the further development of \RIPM{}s. 
        Since this work mainly focuses on the theoretical aspects of \RIPM{}s, further practical improvements would be valuable:
        \begin{enumerate}
            \item \textbf{Dynamic update of barrier parameters:} 
            In the numerical experiments, we adopt monotone updates of the barrier parameter \cref{eq:barrparamupdaterule} together with stopping conditions for the inner iterations. 
            In Euclidean optimization, Armand et al.~\cite{Armandetal2012FromGlobtoLocConvofIntMethsforNonlinOptim} proposed sophisticated update procedures, accompanied by an analysis of the smooth transition between global and local convergence. 
            Incorporating such techniques into our methods could improve practical performance.

            \item \textbf{Development of an inner iteration with global convergence:}
            Currently, 
            \RIPTRM{} in \cite{Obaraetal2025APrimalDualIPTRMfor2ndOrdStnryPtofRiemIneqCstrOptimProbs} is the only available method whose inner iteration exhibits a global convergence in accordance with the stopping conditions.
            However, this algorithm applies only to \RNLO{}~\cref{prob:RNLO} with $\eqset=\emptyset$, and no algorithms are known to achieve such global convergence for the full \RNLO{}~\cref{prob:RNLO} setting (\RNLO{} with $\eqset=\emptyset$).
            Thus, an important avenue for future work is to design methods that establish global convergence in the complete \RNLO{}~\cref{prob:RNLO} setting.
            \isextendedversion{One possibility would be to extend \RIPTRM{} in \cite{Obaraetal2025APrimalDualIPTRMfor2ndOrdStnryPtofRiemIneqCstrOptimProbs} to the full \RNLO{}~\cref{prob:RNLO} setting based on Euclidean \IPTRM{}s that manage both inequality and equality constraints~\cite{ByrdGilbertNocedal2000TRMethBsdIntPtTechniqueforNLP,ByrdLiuNocedal1997LocalBehavIntPtMethforNLP,YamashitaYabeTanabe2005GlobalSuperlinConvPrimalDualIntPtTRNethforLgeScaleCstrOptim}.
            Other types of algorithms could also serve as inner iterations, provided they possess the global convergence.}{}

            \item \textbf{Improvement in the computation of the search direction:} 
            Numerical errors are unavoidable when solving linear systems and may compromise the theoretical guarantees in practice. 
            Moreover, inexact computation of search directions would be important in large-scale optimization, since the exact computation can be inefficient in such settings.
            Therefore, addressing such inexactness constitutes an important direction for future research in our setting, as well as in the Euclidean case~\cite{Wright2001EfctofFinitePrecArithonIntPtMethforNonlinProgram,CurtisSchenkWumlautchter2010IntPtAlgoforLgeScaleNLOwithInexactStepComput}.
            Another promising direction is to exploit problem-specific structures in applications; for instance, leveraging the sparsity of the Hessian can significantly improve computational efficiency.
            \isextendedversion{
            \item \textbf{Handling of degenerate problems:} 
            In this paper, we assume the \LICQ{} in \cref{assu:SCLICQOSOSC}.
            In the Euclidean setting, several works~\cite{WrightOrban2002PropofLogBarrFunonDegenNLP,VicenteWright2002LocalConvPDMethforDegenNonlinProg,YamashitaYabe2005QuadConvofPDIPMforDegenNonlinOptimProb} analyze the local behavior of \IPM{}s around a degenerate solution and design \IPM{}s to achieve the local convergence under the Mangasarian-Fromovitz constraint qualifications, a weaker assumption than the \LICQ{}.
            Extending these techniques to the Riemannian setting is an important future avenue.
            \item \textbf{Local convergence to non-strict local minima:}
            In our local convergence analysis, we assume the \SOSC{} condition in \cref{assu:SCLICQOSOSC}, which implies the strict local optimality.
            Considering local convergence to non-strict local minima is a challenging direction for future research.
            \isextendedversion{For example, in the unconstrained setting, the local convergence properties of trust region methods to such minima are shown to significantly depend on the choice of the search direction.
            Rebjock and Boumal~\cite{RebjockBoumal2024FastConvtoNonIsolMinFourEquivCondsforCtwoFuns} demonstrated that trust region methods can fail to achieve local convergence when using the \exactsolution{}, while using a direction by a truncated conjugate gradient method preserves local convergence under certain assumptions~\cite{RebjockBoumal2024FastConvofTRforNonIsolMinviaAnalofCGonIndefMat}.}{
            For example, in the unconstrained setting, the local convergence properties of trust region methods to such minima are shown to significantly depend on the choice of the search direction~\cite{RebjockBoumal2024FastConvtoNonIsolMinFourEquivCondsforCtwoFuns,RebjockBoumal2024FastConvofTRforNonIsolMinviaAnalofCGonIndefMat}.}
            Building on these insights in the constrained setting would be an interesting avenue.
            }{}
        \end{enumerate}

    \begin{acknowledgements}
        This work was supported by JSPS KAKENHI Grant Numbers 20K19748, 22KJ0563, 23H03351, and 25K15008 and is part of the results of Value Exchange Engineering, a joint research project between R4D, Mercari, Inc. and the RIISE.
        The authors are grateful to the two anonymous referees for their insightful comments and constructive suggestions, and to the editor for the careful handling of the manuscript.
    \end{acknowledgements}

    {\small \noindent
    \textbf{Data Availability}
    The datasets generated and/or analyzed during the current study are available from the corresponding author upon reasonable request.
    }

    
    \appendix  
    \setcounter{section}{0}
    \renewcommand{\thesection}{Appendix~\Alph{section}}
    \isextendedversion{

\section{Proof of \cref{lemm:Jacobnonsingaccum}}\label{appx:proofJacobnonsingaccum}

    \begin{proof}
        Let $\dirallvaraccum \in \tanspc[\allvaraccum]\mani$ satisfy $\Jacobian[\barrKKTvecfld]\paren*{\allvaraccum}\sbra*{\dirallvaraccum} = \zerovec[\allvaraccum]$;
        that is, under $\tanspc[{\ineqLagmult[]}]\setR[\ineqdime]\cid\setR[\ineqdime]$ and $\tanspc[{\eqLagmult[]}]\setR[\eqdime]\cid\setR[\eqdime]$,
        \begin{align}
            &\Hess[\pt]\Lagfun\paren*{\allvaraccum}\sbra*{\dirptaccum} - \sum_{\ineqidx\in\ineqset}\dirineqLagmultaccum[\ineqidx] \gradstr\ineqfun[\ineqidx]\paren*{\ptaccum} + \sum_{\eqidx\in\eqset}\direqLagmultaccum[\eqidx] \gradstr\eqfun[\eqidx]\paren*{\ptaccum} = \zerovec[\ptaccum],\label{eq:Jacobzeroeqdx}\\
            &\ineqLagmultaccum[\ineqidx]\metr[\ptaccum]{\gradstr\ineqfun[\ineqidx]\paren*{\ptaccum}}{\dirptaccum} + \dirineqLagmultaccum[\ineqidx]\ineqfun[\ineqidx]\paren*{\ptaccum} = 0 \text{ for all } \ineqidx\in\ineqset,\label{eq:Jacobzeroeqdy}\\
            &\metr[\ptaccum]{\gradstr\eqfun[\eqidx]\paren*{\ptaccum}}{\dirptaccum}=0 \text{ for all } \eqidx\in\eqset.\label{eq:Jacobzeroeqgradeqfun}
            
        \end{align}
        We will show that such $\dirallvaraccum$ is actually the zero vector.
        For each $\ineqidx\notin\activeineqset\paren*{\ptaccum}$, $\ineqLagmultaccum[\ineqidx] = 0$ and $\ineqfun[\ineqidx]\paren*{\ptaccum} > 0$ hold by the complementarity condition in \cref{eq:KKTconditions}. 
        This, together with \cref{eq:Jacobzeroeqdy}, implies 
        \begin{align}\label{eq:dirineqzeroinactive}
            \dirineqLagmultaccum[\ineqidx]=0 \text{ for any } \ineqidx\notin\activeineqset\paren*{\ptaccum}.
        \end{align}
        On the other hand, for each $\ineqidx\in\activeineqset\paren*{\ptaccum}$, $\ineqLagmultaccum[\ineqidx] > 0$ and $\ineqfun[\ineqidx]\paren*{\ptaccum} = 0$ hold by the \SC{}.
        Together with \cref{eq:Jacobzeroeqdy} again, this implies
        \begin{align}\label{eq:metrgradineqdxactive}
            \metr[\ptaccum]{\gradstr\ineqfun[\ineqidx]\paren*{\ptaccum}}{\dirptaccum} = 0 \text{ for any } \ineqidx \in \activeineqset\paren*{\ptaccum}.
        \end{align}
        From \cref{eq:Jacobzeroeqgradeqfun}, \cref{eq:metrgradineqdxactive}, and the SC, again,  $\dirptaccum\in\criticalcone\paren*{\allvaraccum}$ holds, which, together with \cref{eq:Jacobzeroeqdx,eq:dirineqzeroinactive}, implies $\metr[\ptaccum]{\Hess[\pt]\Lagfun\paren*{\allvaraccum}\sbra*{\dirptaccum}}{\dirptaccum} = \zerovec[\ptaccum]$.
        Therefore, $\dirptaccum=\zerovec[\ptaccum]$ holds by the \SOSC{}.
        Substituting $\dirptaccum=\zerovec[\ptaccum]$ and \cref{eq:dirineqzeroinactive}
        into \cref{eq:Jacobzeroeqdx}, we obtain $\sum_{\ineqidx\in\activeineqset\paren*{\ptaccum}} \dirineqLagmultaccum[\ineqidx] \gradstr\ineqfun[\ineqidx]\paren*{\ptaccum} + \sum_{\eqidx\in\eqset}\direqLagmultaccum[\eqidx] \gradstr\eqfun[\eqidx]\paren*{\ptaccum} = \zerovec[\ptaccum]$.
        Thus, by the \LICQ{}, it holds that $\direqLagmultaccum=0$ for $\eqidx\in\eqset$ and  $\dirineqLagmultaccum[\ineqidx]=0$ for $\ineqidx\in\activeineqset\paren*{\ptaccum}$.
        Hence we have $\dirallvaraccum=\zerovec[\allvaraccum]$, which completes the proof.
        \qed
    \end{proof}

\section{Proof of \cref{lemm:ineqLagmultineqfunbarrparambounds}}\label{appx:proofineqLagmultineqfunbarrparambounds}
    \begin{proof}
        We first prove \cref{lemm:activebarrparambounds}.
        Let $\constineqLagmultlower\coloneqq\min_{\ineqidx\in\activeineqset\paren*{\ptaccum}}\frac{1}{2}\ineqLagmultaccum[\ineqidx]$ and  $\constineqLagmultupper\coloneqq\max_{\ineqidx\in\activeineqset\paren*{\ptaccum}} 2\ineqLagmultaccum[\ineqidx]$.
        Under the \SC{}, for each $\ineqidx\in\activeineqset\paren*{\ptaccum}$ and any $\otriteridx$-th iteration of \cref{algo:RIPMOuter} satisfying that $\allvarotriter$ is sufficiently close to $\allvaraccum$, it holds that $0 < \constineqLagmultlower \leq \ineqLagmultotriter[\ineqidx] \leq \constineqLagmultupper$ by definition.
        From \cref{ineq:forcingfuncomplbounded}, we have $\abs*{\ineqLagmultotriter\ineqfun[\ineqidx]\paren*{\ptotriter} - \barrparamotriterm} \leq \constcomplupper\barrparamotriterm$, which implies
        \isextendedversion{
        \begin{align}
            \frac{1}{\constineqLagmultupper}\paren*{1 - \constcomplupper}\barrparamotriterm \leq \frac{1}{\ineqLagmultotriter}\paren*{1 - \constcomplupper}\barrparamotriterm \leq \ineqfun[\ineqidx]\paren*{\ptotriter} \leq \frac{1}{\ineqLagmultotriter}\paren*{1 + \constcomplupper}\barrparamotriterm \leq \frac{1}{\constineqLagmultlower}\paren*{1 + \constcomplupper}\barrparamotriterm.
        \end{align}}
        {$\frac{1}{\constineqLagmultupper}\paren*{1 - \constcomplupper}\barrparamotriterm \leq \frac{1}{\ineqLagmultotriter}\paren*{1 - \constcomplupper}\barrparamotriterm \leq \ineqfun[\ineqidx]\paren*{\ptotriter} \leq \frac{1}{\ineqLagmultotriter}\paren*{1 + \constcomplupper}\barrparamotriterm \leq \frac{1}{\constineqLagmultlower}\paren*{1 + \constcomplupper}\barrparamotriterm$.}
        The proof of \cref{lemm:activebarrparambounds} is complete.

        Similarly, we prove \cref{lemm:inactivebarrparambounds}.
        Define $\constinactivecompllower\coloneqq\min_{\ineqidx\notin\activeineqset\paren*{\ptaccum}}\frac{1}{2}\ineqfun[\ineqidx]\paren*{\ptotriter}$ and let $\constinactivecomplupper\coloneqq\max_{\ineqidx\notin\activeineqset\paren*{\ptaccum}} 2\ineqfun[\ineqidx]\paren*{\ptotriter}$.
        Under the SC{}, for each $\ineqidx\notin\activeineqset\paren*{\ptaccum}$ and any $\otriteridx$-th iteration of \cref{algo:RIPMOuter} satisfying that $\allvarotriter$ is sufficiently close to $\allvaraccum$, it holds that $0 < \constinactivecompllower \leq \ineqfun[\ineqidx]\paren*{\ptotriter} \leq \constinactivecomplupper$.
        By \cref{ineq:forcingfuncomplbounded}, we have $\abs*{\ineqLagmultotriter\ineqfun[\ineqidx]\paren*{\ptotriter} - \barrparamotriterm} \leq \constcomplupper\barrparamotriterm$, which implies
        \isextendedversion{
        \begin{align}
            &\frac{1}{\constinactivecomplupper}\paren*{1 - \constcomplupper}\barrparamotriterm \leq \frac{1}{\ineqfun[\ineqidx]\paren*{\ptotriter}}\paren*{1 - \constcomplupper}\barrparamotriterm \leq \ineqLagmultotriter \text{ and } \ineqLagmultotriter \leq \frac{1}{\ineqfun[\ineqidx]\paren*{\ptotriter}}\paren*{1 + \constcomplupper}\barrparamotriterm \leq \frac{1}{\constinactivecompllower}\paren*{1 + \constcomplupper}\barrparamotriterm.
        \end{align}}
        {$\frac{1}{\constinactivecomplupper}\paren*{1 - \constcomplupper}\barrparamotriterm \leq \frac{1}{\ineqfun[\ineqidx]\paren*{\ptotriter}}\paren*{1 - \constcomplupper}\barrparamotriterm \leq \ineqLagmultotriter$ and $\ineqLagmultotriter \leq \frac{1}{\ineqfun[\ineqidx]\paren*{\ptotriter}}\paren*{1 + \constcomplupper}\barrparamotriterm \leq \frac{1}{\constinactivecompllower}\paren*{1 + \constcomplupper}\barrparamotriterm$.}
        The proof of \cref{lemm:inactivebarrparambounds} is complete.
        \qed
    \end{proof}

\section{Proof of \cref{lemm:dirNewtonallvarbarrparambound}}\label{appx:proofdirNewtonallvarbarrparambound}
    \begin{proof}
        From \cref{coro:Jacobnonsingularardaccum}, we see that the Newton step \cref{def:Newtondirection} is well-defined in a sufficiently small neighborhood of $\allvaraccum$.
        In addition, 
        from \cref{lemm:Jacobnonsingbouundardaccum}, 
        there exists $\consteig > 0$ such that, for any $\otriteridx$-th iteration of \cref{algo:RIPMOuter} satisfying that $\allvarotriter$ is sufficiently close to $\allvaraccum$ and $\barrparamotriter$ is sufficiently small, $\opnorm{\inv{\Jacobian[\barrKKTvecfld]\paren*{\allvarotriter}}}\leq \consteig$  holds and hence
        \isextendedversion{
        \begin{align}
            &\Riemnorm[{\allvarotriter}]{\dirNewtonallvarbarrparamotriter} \leq \opnorm{\inv{\Jacobian[\barrKKTvecfld]\paren*{\allvarotriter}}} \Riemnorm[{\allvarotriter}]{\barrKKTvecfld\paren*{\allvarotriter;\barrparamotriter}}\\
            &\leq \consteig \paren*{\Riemnorm[{\ptotriter}]{\gradstr[\pt]\Lagfun\paren*{\allvarotriter}} + \norm*{\IneqLagmultmatotriter\ineqfun[]\paren*{\ptotriter} - \barrparamotriterm\onevec[\ineqdime]} + \paren*{\barrparamotriterm - \barrparamotriter}\Riemnorm[\allvarotriter]{\begin{bmatrix} \zerovec[\ptotriter] \\ \onevec[\ineqdime] \end{bmatrix}} + \norm*{\eqfun[]\paren*{\ptotriter}}} \\
            &\leq \consteig \paren*{\forcingfungradLag\paren*{\barrparamotriterm} + \forcingfuncompl\paren*{\barrparamotriterm} + \sqrt{\ineqdime}\paren*{\barrparamotriterm + \barrparamotriter} + \forcingfuneqvio\paren*{\barrparamotriterm}}\\
            &\leq \consteig \paren*{\constgradLagupper + \constcomplupper + 2 \sqrt{\ineqdime} + \consteqvioupper}\barrparamotriterm,
        \end{align}}
        {$\Riemnorm[{\allvarotriter}]{\dirNewtonallvarbarrparamotriter} \leq \opnorm{\inv{\Jacobian[\barrKKTvecfld]\paren*{\allvarotriter}}} \Riemnorm[{\allvarotriter}]{\barrKKTvecfld\paren*{\allvarotriter;\barrparamotriter}} \leq \consteig \paren*{\Riemnorm[{\ptotriter}]{\gradstr[\pt]\Lagfun\paren*{\allvarotriter}} + \norm*{\IneqLagmultmatotriter\ineqfun[]\paren*{\ptotriter} - \barrparamotriterm\onevec[\ineqdime]} + \paren*{\barrparamotriterm - \barrparamotriter}\Riemnorm[\allvarotriter]{\begin{bmatrix} \zerovec[\ptotriter] \\ \onevec[\ineqdime] \end{bmatrix}} + \norm*{\eqfun[]\paren*{\ptotriter}}} \leq \consteig \paren*{\forcingfungradLag\paren*{\barrparamotriterm} + \forcingfuncompl\paren*{\barrparamotriterm} + \sqrt{\ineqdime}\paren*{\barrparamotriterm + \barrparamotriter} + \forcingfuneqvio\paren*{\barrparamotriterm}}\leq \consteig \paren*{\constgradLagupper + \constcomplupper + 2 \sqrt{\ineqdime} + \consteqvioupper}\barrparamotriterm$,}
        where the third inequality follows from \cref{eq:stopcondKKT,eq:stopcondbarrcompl} and the fourth one from \cref{ineq:forcingfungradLagbounded,ineq:forcingfuncomplbounded,ineq:forcingfuneqviobounded,assu:barrparammonononincr}.
        Letting $\constdirNewton\coloneqq\constgradLagupper + \constcomplupper + 2 \sqrt{\ineqdime} + \consteqvioupper > 0$, we complete the proof.
        \qed
    \end{proof}

    }{}

\section{Proof of \cref{lemm:dirNewtonineqfunpos}}\label{appx:proofdirNewtonineqfunpos}
    \begin{proof}
        For each $\ineqidx\notin\activeineqset\paren*{\ptaccum}$, it follows from the continuities of $\ineqfun[\ineqidx]$ and $\retr[]$ that 
        \isextendedversion{
        \begin{align}
            \ineqfun[\ineqidx]\paren*{\ptotriter} \geq \frac{1}{2}\ineqfun[\ineqidx]\paren*{\ptaccum} > 0 \text{ and } \abs*{\ineqfun[\ineqidx]\circ\retr[]\paren*{\ptotriter, \dirNewtonptbarrparamotriter} -  \ineqfun[\ineqidx]\circ\retr[]\paren*{\ptotriter, \zerovec[\ptotriter]}} \leq \frac{1}{3}\ineqfun[\ineqidx]\paren*{\ptaccum}
        \end{align}}
        {$\ineqfun[\ineqidx]\paren*{\ptotriter} \geq \frac{1}{2}\ineqfun[\ineqidx]\paren*{\ptaccum} > 0 \text{ and } \abs*{\ineqfun[\ineqidx]\circ\retr[]\paren*{\ptotriter, \dirNewtonptbarrparamotriter} -  \ineqfun[\ineqidx]\circ\retr[]\paren*{\ptotriter, \zerovec[\ptotriter]}} \leq \frac{1}{3}\ineqfun[\ineqidx]\paren*{\ptaccum}$}
        for any $\otriteridx$-th iteration of \cref{algo:RIPMOuter} satisfying that $\allvarotriter$ is sufficiently close to $\allvaraccum$ and $\barrparamotriterm$ is sufficiently small.
        This, together with \isextendedversion{\cref{eq:retrzero}}{the definition of the retraction}, implies that
        \isextendedversion{
        \begin{align}
            \ineqfun[\ineqidx]\paren*{\retr[\ptotriter]\paren*{\dirNewtonptbarrparamotriter}} \geq \ineqfun[\ineqidx]\paren*{\ptotriter} - \frac{1}{3}\ineqfun[\ineqidx]\paren*{\ptaccum} \geq \frac{1}{6}\ineqfun[\ineqidx]\paren*{\ptaccum} > 0.
        \end{align}}
        {$\ineqfun[\ineqidx]\paren*{\retr[\ptotriter]\paren*{\dirNewtonptbarrparamotriter}} \geq \ineqfun[\ineqidx]\paren*{\ptotriter} - \frac{1}{3}\ineqfun[\ineqidx]\paren*{\ptaccum} \geq \frac{1}{6}\ineqfun[\ineqidx]\paren*{\ptaccum} > 0$.}

        Next, we consider the case $\ineqidx\in\activeineqset\paren*{\ptaccum}$.
        We have
        \isextendedversion{
        \begin{align}
            \begin{split}\label{eq:pullbackineqfunTaylor}
                &\pullback[\ptotriter]{\ineqfun[\ineqidx]}\paren*{\dirNewtonptbarrparamotriter}\\
                &= \pullback[\ptotriter]{\ineqfun[\ineqidx]}\paren*{\zerovec[\ptotriter]} + \D\pullback[\ptotriter]{\ineqfun[\ineqidx]}\paren*{\zerovec[\ptotriter]}\sbra*{\dirNewtonptbarrparamotriter} + \int_{0}^{1}\D\paren*{\pullback[\ptotriter]{\ineqfun[\ineqidx]}\paren*{\tmesix\dirNewtonptbarrparamotriter} - \pullback[\ptotriter]{\ineqfun[\ineqidx]}\paren*{\zerovec[\ptotriter]}}\sbra*{\dirNewtonptbarrparamotriter}\dd{\tmesix}.
            \end{split}
        \end{align}}
        {\begin{equation}
        \begin{aligned}[t]\label{eq:pullbackineqfunTaylor}
                &\pullback[\ptotriter]{\ineqfun[\ineqidx]}\paren*{\dirNewtonptbarrparamotriter}\\
                &= \pullback[\ptotriter]{\ineqfun[\ineqidx]}\paren*{\zerovec[\ptotriter]} + \D\pullback[\ptotriter]{\ineqfun[\ineqidx]}\paren*{\zerovec[\ptotriter]}\sbra*{\dirNewtonptbarrparamotriter} + \int_{0}^{1}\D\paren*{\pullback[\ptotriter]{\ineqfun[\ineqidx]}\paren*{\tmesix\dirNewtonptbarrparamotriter} - \pullback[\ptotriter]{\ineqfun[\ineqidx]}\paren*{\zerovec[\ptotriter]}}\sbra*{\dirNewtonptbarrparamotriter}\dd{\tmesix}.
        \end{aligned}
        \end{equation}}
        In what follows, we provide the bounds on the first two terms and on the last one.
        As for the first two terms, for any $\otriteridx$-th iteration of \cref{algo:RIPMOuter} satisfying that $\allvarotriter$ is sufficiently close to $\allvaraccum$ and $\barrparamotriterm$ is sufficiently small, we have 
        \isextendedversion{
        \begin{align}
            \begin{split}\label{ineq:ineqfunbarrparamboundone}
                &\pullback[\ptotriter]{\ineqfun[\ineqidx]}\paren*{\zerovec[\ptotriter]} + \D\pullback[\ptotriter]{\ineqfun[\ineqidx]}\paren*{\zerovec[\ptotriter]}\sbra*{\dirNewtonptbarrparamotriter} = \ineqfun[\ineqidx]\paren*{\ptotriter} + \D\ineqfun[\ineqidx]\paren*{\ptotriter}\sbra*{\dirNewtonptbarrparamotriter}\\
                &= - \frac{\sbra*{\dirNewtonineqLagmultbarrparamotriter[]}_{\ineqidx}\ineqfun[\ineqidx]\paren*{\ptotriter}}{\ineqLagmultotriter[\ineqidx]} + \frac{\barrparamotriter}{\ineqLagmultotriter[\ineqidx]} \geq - \frac{\constdirNewton\paren*{1 + \constcomplupper}}{\constineqLagmultlower}\barrparamotriterm^{2} + \frac{\barrparamotriter}{\ineqLagmultotriter[\ineqidx]},
            \end{split}
        \end{align}
        }{
        \begin{align}
            \begin{split}\label{ineq:ineqfunbarrparamboundone}
                &\pullback[\ptotriter]{\ineqfun[\ineqidx]}\paren*{\zerovec[\ptotriter]} + \D\pullback[\ptotriter]{\ineqfun[\ineqidx]}\paren*{\zerovec[\ptotriter]}\sbra*{\dirNewtonptbarrparamotriter} = \ineqfun[\ineqidx]\paren*{\ptotriter} + \D\ineqfun[\ineqidx]\paren*{\ptotriter}\sbra*{\dirNewtonptbarrparamotriter}\\
                &= - \frac{\sbra*{\dirNewtonineqLagmultbarrparamotriter[]}_{\ineqidx}\ineqfun[\ineqidx]\paren*{\ptotriter}}{\ineqLagmultotriter[\ineqidx]} + \frac{\barrparamotriter}{\ineqLagmultotriter[\ineqidx]} \geq - \frac{\constdirNewton\paren*{1 + \constcomplupper}}{\constineqLagmultlower}\barrparamotriterm^{2} + \frac{\barrparamotriter}{\ineqLagmultotriter[\ineqidx]},
            \end{split}
        \end{align}}
        where the first equality follows from \isextendedversion{\cref{eq:retrzero,eq:retrdiffzero}}{the definition of the retraction}, the second one from \cref{eq:dirNewtoninexidxequation}, and the inequality from \isextendedversion{\cref{ineq:ineqfunLagmultbarrparamboundactive}}{\enumicref{lemm:ineqLagmultineqfunbarrparambounds}{lemm:activebarrparambounds}}
        and \cref{lemm:dirNewtonallvarbarrparambound}.
        We now consider the third term.
        By the twice continuous differentiability of $\brc*{\ineqfun[\ineqidx]}_{\ineqidx\in\ineqset}$ and \cite[Lemma~10.57]{Boumal23IntroOptimSmthMani}, there exist positive scalars $\brc*{\Lipschitzconstsix[{\ineqfun[\ineqidx]}]}_{\ineqidx\in\ineqset}$ such that, for all $\pt\in\mani$ sufficiently close to $\ptaccum$ and any $\tanvecone[\pt]\in\tanspc[\pt]\mani$ sufficiently small,  
        \begin{equation}
            \Riemnorm[]{\gradstr\pullback[\pt]{\ineqfun[\ineqidx]}\paren*{\tanvecone[\pt]} - \gradstr\pullback[\pt]{\ineqfun[\ineqidx]}\paren*{\zerovec[\pt]}} \leq \Lipschitzconstsix[{\ineqfun[\ineqidx]}]\Riemnorm[\pt]{\tanvecone[\pt]} \text{ for all } \ineqidx \in \ineqset.\label{ineq:pullbackgradineqfunLipschitzcontone}
        \end{equation}            
        Then, for any $\otriteridx$-th iteration of \cref{algo:RIPMOuter} satisfying that $\allvarotriter$ is sufficiently close to $\allvaraccum$ and $\barrparamotriterm$ is sufficiently small, we have
        \isextendedversion{
        \begin{align}
            \begin{split}\label{ineq:ineqfunbarrparamboundtwo}
                &\abs*{\int_{0}^{1}\D\paren*{\pullback[\ptotriter]{\ineqfun[\ineqidx]}\paren*{\tmesix\dirNewtonptbarrparamotriter} - \pullback[\ptotriter]{\ineqfun[\ineqidx]}\paren*{\zerovec[\ptotriter]}}\sbra*{\dirNewtonptbarrparamotriter}\dd{\tmesix}}\\
                &\leq \int_{0}^{1}\abs*{\metr[\ptotriter]{\gradstr\pullback[\ptotriter]{\ineqfun[\ineqidx]}\paren*{\tmesix\dirNewtonptbarrparamotriter} - \gradstr\pullback[\ptotriter]{\ineqfun[\ineqidx]}\paren*{\zerovec[\ptotriter]}}{\dirNewtonptbarrparamotriter}}\dd{\tmesix}\\
                &\leq \int_{0}^{1}\Riemnorm[\ptotriter]{\gradstr\pullback[\ptotriter]{\ineqfun[\ineqidx]}\paren*{\tmesix\dirNewtonptbarrparamotriter} - \gradstr\pullback[\ptotriter]{\ineqfun[\ineqidx]}\paren*{\zerovec[\ptotriter]}}\Riemnorm[\ptotriter]{\dirNewtonptbarrparamotriter}\dd{\tmesix}\\
                &\leq \int_{0}^{1}\Lipschitzconstsix[{\ineqfun[\ineqidx]}]\tmesix\Riemnorm[\ptotriter]{\dirNewtonptbarrparamotriter}^{2} \dd{\tmesix} = \frac{\Lipschitzconstsix[{\ineqfun[\ineqidx]}]}{2}\Riemnorm[\ptotriter]{\dirNewtonptbarrparamotriter}^{2} \leq \frac{\constdirNewton^{2}\Lipschitzconstsix[{\ineqfun[\ineqidx]}]}{2}\barrparamotriterm^{2},
            \end{split}
        \end{align}}{
        \begin{align}
            \begin{split}\label{ineq:ineqfunbarrparamboundtwo}
                &\abs*{\int_{0}^{1}\D\paren*{\pullback[\ptotriter]{\ineqfun[\ineqidx]}\paren*{\tmesix\dirNewtonptbarrparamotriter} - \pullback[\ptotriter]{\ineqfun[\ineqidx]}\paren*{\zerovec[\ptotriter]}}\sbra*{\dirNewtonptbarrparamotriter}\dd{\tmesix}} \leq \int_{0}^{1}\abs*{\metr[\ptotriter]{\gradstr\pullback[\ptotriter]{\ineqfun[\ineqidx]}\paren*{\tmesix\dirNewtonptbarrparamotriter}\\
                &- \gradstr\pullback[\ptotriter]{\ineqfun[\ineqidx]}\paren*{\zerovec[\ptotriter]}}{\dirNewtonptbarrparamotriter}}\dd{\tmesix} \leq \int_{0}^{1}\Riemnorm[\ptotriter]{\gradstr\pullback[\ptotriter]{\ineqfun[\ineqidx]}\paren*{\tmesix\dirNewtonptbarrparamotriter} - \gradstr\pullback[\ptotriter]{\ineqfun[\ineqidx]}\paren*{\zerovec[\ptotriter]}}\\
                &\Riemnorm[\ptotriter]{\dirNewtonptbarrparamotriter}\dd{\tmesix} \leq \int_{0}^{1}\Lipschitzconstsix[{\ineqfun[\ineqidx]}]\tmesix\Riemnorm[\ptotriter]{\dirNewtonptbarrparamotriter}^{2} \dd{\tmesix} = \frac{\Lipschitzconstsix[{\ineqfun[\ineqidx]}]}{2}\Riemnorm[\ptotriter]{\dirNewtonptbarrparamotriter}^{2} \leq \frac{\constdirNewton^{2}\Lipschitzconstsix[{\ineqfun[\ineqidx]}]}{2}\barrparamotriterm^{2},
            \end{split}
        \end{align}}
        where the third inequality follows from \cref{ineq:pullbackgradineqfunLipschitzcontone} and \cref{lemm:dirNewtonallvarbarrparambound} and the last one from \cref{lemm:dirNewtonallvarbarrparambound}, again.
        Combining \cref{eq:pullbackineqfunTaylor} with \cref{ineq:ineqfunbarrparamboundone,ineq:ineqfunbarrparamboundtwo} yields 
        \isextendedversion{
        \begin{align}
            \pullback[\ptotriter]{\ineqfun[\ineqidx]}\paren*{\dirNewtonptbarrparamotriter} \geq 
            \barrparamotriter \paren*{\frac{1}{\ineqLagmultotriter[\ineqidx]} - \paren*{\frac{\constdirNewton\paren*{1 + \constcomplupper}}{\constineqLagmultlower} + \frac{\constdirNewton^{2}\Lipschitzconstsix[{\ineqfun[\ineqidx]}]}{2}}\frac{\barrparamotriterm^{2}}{\barrparamotriter}}.
        \end{align}}
        {$\pullback[\ptotriter]{\ineqfun[\ineqidx]}\paren*{\dirNewtonptbarrparamotriter} \geq 
            \barrparamotriter \paren*{\frac{1}{\ineqLagmultotriter[\ineqidx]} - \paren*{\frac{\constdirNewton\paren*{1 + \constcomplupper}}{\constineqLagmultlower} + \frac{\constdirNewton^{2}\Lipschitzconstsix[{\ineqfun[\ineqidx]}]}{2}}\frac{\barrparamotriterm^{2}}{\barrparamotriter}}$.}
        From \cref{assu:barrparamdecreaseupodr}, the right-hand side is positive for $\barrparamotriter > 0$ sufficiently small, 
        which implies $\pullback[\ptotriter]{\ineqfun[\ineqidx]}\paren*{\dirNewtonptbarrparamotriter} > 0$ for any $\barrparamotriter > 0$ sufficiently small.
        We complete the proof.
        \qed
    \end{proof}

    \isextendedversion{

\section{Proof of \cref{lemm:dirNewtonineqLagmultpos}}\label{appx:proofdirNewtonineqLagmultpos}
    \begin{proof}
        Recall the definition \cref{def:activeineqset} of $\activeineqset\paren*{\pt}$ with $\pt\in\mani$.
        For each $\ineqidx\in\activeineqset\paren*{\ptaccum}$, since the point $\ineqLagmultotriter[\ineqidx] + \sbra*{\dirNewtonineqLagmultbarrparamotriter[]}_{\ineqidx}$ can be made arbitrarily close to $\ineqLagmultaccum[\ineqidx]$ by considering the $\otriteridx$-th iterations of \cref{algo:RIPMOuter} satisfying that $\allvarotriter$ is sufficiently close to $\allvaraccum$ and $\barrparamotriterm$ is sufficiently small, it holds that $\ineqLagmultotriter[\ineqidx] + \sbra*{\dirNewtonineqLagmultbarrparamotriter[]}_{\ineqidx} > 0$ for such $\allvarotriter$ and $\barrparamotriter$ under the \SC{}.

        Next, we consider the case $\ineqidx\notin\activeineqset\paren*{\ptaccum}$.
        Note that the function $\pt\mapsto\Riemnorm[\pt]{\gradstr\ineqfun[\ineqidx]\paren*{\pt}}$ is bounded around $\ptaccum\in\mani$ for each $\ineqidx\in\ineqset$ by the continuity of the gradient. 
        Letting $\tholdvalnin > 0$ be a sufficiently large scalar, we obtain
        \isextendedversion{
        \begin{align}
                &\ineqLagmultotriter[\ineqidx] + \sbra*{\dirNewtonineqLagmultbarrparamotriter[]}_{\ineqidx} = \inv{\ineqfun[\ineqidx]\paren*{\ptotriter}}\paren*{\barrparamotriter - \ineqLagmultotriter[\ineqidx]\metr[\ptotriter]{\gradstr\ineqfun[\ineqidx]\paren*{\ptotriter}}{\dirNewtonptbarrparamotriter}}\\
                &\geq \inv{\ineqfun[\ineqidx]\paren*{\ptotriter}}\paren*{\barrparamotriter - \ineqLagmultotriter[\ineqidx]\Riemnorm[\ptotriter]{\gradstr\ineqfun[\ineqidx]\paren*{\ptotriter}}\Riemnorm[\ptotriter]{\dirNewtonptbarrparamotriter}}\\
                &\geq \inv{\ineqfun[\ineqidx]\paren*{\ptotriter}}\paren*{\barrparamotriter - \frac{\constdirNewton\paren*{1 + \constcomplupper}}{\constinactivecompllower}\Riemnorm[\ptotriter]{\gradstr\ineqfun[\ineqidx]\paren*{\ptotriter}}\barrparamotriterm^{2}}\\
                &\geq \frac{\barrparamotriter}{\ineqfun[\ineqidx]\paren*{\ptotriter}}\paren*{1 - \tholdvalnin\frac{\barrparamotriterm^{2}}{\barrparamotriter}} \geq \frac{\barrparamotriter}{\constinactivecompllower}\paren*{1 - \tholdvalnin\frac{\barrparamotriterm^{2}}{\barrparamotriter}} > 0,
        \end{align}}
        {$\ineqLagmultotriter[\ineqidx] + \sbra*{\dirNewtonineqLagmultbarrparamotriter[]}_{\ineqidx} = \inv{\ineqfun[\ineqidx]\paren*{\ptotriter}}\paren*{\barrparamotriter - \ineqLagmultotriter[\ineqidx]\metr[\ptotriter]{\gradstr\ineqfun[\ineqidx]\paren*{\ptotriter}}{\dirNewtonptbarrparamotriter}} \geq \inv{\ineqfun[\ineqidx]\paren*{\ptotriter}}\paren*{\barrparamotriter - \ineqLagmultotriter[\ineqidx]\Riemnorm[\ptotriter]{\gradstr\ineqfun[\ineqidx]\paren*{\ptotriter}} \Riemnorm[\ptotriter]{\dirNewtonptbarrparamotriter}} \geq \inv{\ineqfun[\ineqidx]\paren*{\ptotriter}}\paren*{\barrparamotriter - \frac{\constdirNewton\paren*{1 + \constcomplupper}}{\constinactivecompllower}\Riemnorm[\ptotriter]{\gradstr\ineqfun[\ineqidx]\paren*{\ptotriter}}\barrparamotriterm^{2}} \geq \frac{\barrparamotriter}{\ineqfun[\ineqidx]\paren*{\ptotriter}}\paren*{1 - \tholdvalnin\frac{\barrparamotriterm^{2}}{\barrparamotriter}} \geq \frac{\barrparamotriter}{\constinactivecompllower}\paren*{1 - \tholdvalnin\frac{\barrparamotriterm^{2}}{\barrparamotriter}} > 0$,}
        where the first equality follows from \cref{eq:dirNewtoninexidxequation}, the second inequality follows from \cref{ineq:ineqfunLagmultbarrparamboundinactive} and \cref{lemm:dirNewtonallvarbarrparambound}, the third one from the boundedness of $\Riemnorm[\paren*{\cdot}]{\gradstr\ineqfun[\ineqidx]\paren*{\cdot}}$ around $\allvaraccum$, and the fourth one from \cref{ineq:ineqfunLagmultbarrparamboundinactive} again and the positivity of $\paren*{1 - \tholdvalnin\frac{\barrparamotriterm^{2}}{\barrparamotriter}}$ by \cref{assu:barrparamdecreaseupodr}.
        The proof is complete.
        \qed
    \end{proof}

\section{Proof of \cref{lemm:lincombgradineq}}\label{appx:prooflincombgradineq}

    \begin{proof}
        Let $\paren*{\chart, \opensubset}$ be a chart with $\ptaccum\in\opensubset\subseteq\setR[\dime]$, and let $\stdbasis[\idxfiv]$ be the $\idxfiv$-th standard basis of $\setR[\dime]$.
        For any $\funeig[]\in\funset\paren*{\mani}$, any $\pt\in\opensubset$ and all $\tanvecone[\pt]\in\tanspc[\pt]\mani$, we have
        \begin{align}
            \partderiv[\coordidx]\pullback[\pt]{\funeig}\paren*{\tanvecone[\pt]} = \sum_{\idxfiv}\partderiv[\idxfiv]\funeig\paren*{\retr[\pt]\paren*{\tanvecone[\pt]}}\diffretrmat_{\coordidx}^{\idxfiv}\paren*{\tanvecone[\pt]},
        \end{align}
        where $\partderiv[\coordidx]\pullback[\pt]{\funeig}\paren*{\tanvecone[\pt]}\colon\tanspc[\pt]\mani\to\setR[]$ is the partial derivative of $\pullback[\pt]{\funeig}$ at $\tanvecone[\pt]$ with respect to the $\coordidx$-th variable, $\partderiv[\idxfiv]\funeig\paren*{\cdot}\coloneqq\lim_{\tmeele\downarrow 0}\frac{\funeig[]\circ\inv{\chart}\paren*{\chart\paren*{\cdot} + \tmeele\stdbasis[\idxfiv]} - \funeig[]\circ\inv{\chart}\paren*{\chart\paren*{\cdot}}}{\tmeele}$, $\diffretrmat\paren*{\tanvecone[\pt]}$ denotes the differential of $\retr[\pt]$ at $\tanvecone[\pt]\in\tanspc[\pt]\mani$, and $\diffretrmat_{\coordidx}^{\idxfiv}\paren*{\tanvecone[\pt]}$ is its $\paren*{\coordidx, \idxfiv}$-th element.
        Let $\coordmetr[]\paren*{\pt}\in\setR[\dime\times\dime]$ be the coordinate expression of the Riemannian metric at $\pt\in\opensubset$, and let $\coordmetr[\coordidx\idxfiv]\paren*{\pt}$ be its $\paren*{\coordidx, \idxfiv}$-th element.
        Note that $\coordmetr[]\paren*{\pt}$ is positive-definite.
        Then,
        \begin{align}
            \begin{split}\label{eq:lincombpullbaclgradbound}
                &\Riemnorm[\pt]{\sum_{\idxfou=1}^{\numfun}\vecsix[\idxfou]\gradstr\pullback[\pt]{\funtwo[]}^{\idxfou}\paren*{\tanvecone[\pt]}}^{2} = \sum_{\idxfou[1], \idxfou[2], \coordidx, \idxfiv}\vecsix[{\idxfou[1]}]\vecsix[{\idxfou[2]}]\partderiv[\coordidx]\pullback[\pt]{\funtwo}^{\idxfou[1]}\paren*{\tanvecone[\pt]}\coordmetr[\coordidx\idxfiv]\paren*{\pt}\partderiv[\idxfiv]\pullback[\pt]{\funtwo}^{\idxfou[2]}\paren*{\tanvecone[\pt]}\\
                &= \sum_{\idxfou[1], \idxfou[2], \coordidx, \idxfiv, \idxsix, \idxsev}\vecsix[{\idxfou[1]}]\vecsix[{\idxfou[2]}]\partderiv[\idxsix]\funtwo[{\idxfou[1]}]\paren*{\retr[\pt]\paren*{\tanvecone[\pt]}}\diffretrmat_{\coordidx}^{\idxsix}\paren*{\tanvecone[\pt]}\coordmetr[\coordidx\idxfiv]\paren*{\pt}\diffretrmat_{\idxfiv}^{\idxsev}\paren*{\tanvecone[\pt]}\partderiv[\idxsev]\funtwo[{\idxfou[2]}]\paren*{\retr[\pt]\paren*{\tanvecone[\pt]}}\\
                &= \trsp{\vecsix[]}\trsp{\partderiv[]\funtwo\paren*{\retr[\pt]\paren*{\tanvecone[\pt]}}}\diffretrmat\paren*{\tanvecone[\pt]}\coordmetr[]\paren*{\pt}\diffretrmat\paren*{\tanvecone[\pt]}\partderiv[]\funtwo\paren*{\retr[\pt]\paren*{\tanvecone[\pt]}}\vecsix[],
            \end{split}
        \end{align}            
        where $\partderiv[]\funtwo\paren*{\cdot}\in\setR[\dime\times\numfun]$ is a matrix whose $\paren*{i, j}$-th element is $\partderiv[i]\funtwo[j]\paren*{\cdot}$.
        We also have
        \begin{align}
            \begin{split}\label{eq:lincombgradbound}
                &\Riemnorm[\pt]{\sum_{\idxfou=1}^{\numfun}\vecsix[\idxfou]\gradstr{\funtwo}^{\idxfou}\paren*{\retr[\pt]\paren*{\tanvecone[\pt]}}}^{2}\\
                &= \sum_{\idxfou[1], \idxfou[2], \coordidx, \idxfiv}\vecsix[{\idxfou[1]}]\vecsix[{\idxfou[2]}]\partderiv[\coordidx]\funtwo[{\idxfou[1]}]\paren*{\retr\paren*{\tanvecone[\pt]}}\coordmetr[\coordidx\idxfiv]\paren*{\retr[\pt]\paren{\tanvecone[\pt]}}\partderiv[\idxfiv]\funtwo[{\idxfou[2]}]\paren*{\retr\paren*{\tanvecone[\pt]}}\\
                &=\trsp{\vecsix[]}\trsp{\partderiv[]\funtwo\paren*{\retr[\pt]\paren*{\tanvecone[\pt]}}}\coordmetr[]\paren*{\pt}\partderiv[]\funtwo\paren*{\retr[\pt]\paren*{\tanvecone[\pt]}}\vecsix[].
            \end{split}
        \end{align}
        It follows from \cref{eq:lincombpullbaclgradbound,eq:lincombgradbound} that 
        \begin{align}
            \begin{split}\label{eq:diffnormpullbackgradlincombgradretrlincomb}
                 &\coeffsev^{2}\Riemnorm[\pt]{\sum_{\idxfou=1}^{\numfun}\vecsix[\idxfou]\gradstr\pullback[\pt]{\funtwo}^{\idxfou}\paren*{\tanvecone[\pt]}}^{2} - \Riemnorm[\pt]{\sum_{\idxfou=1}^{\numfun}\vecsix[\idxfou]\gradstr\funtwo[\idxfou]\paren*{\retr[\pt]\paren*{\tanvecone[\pt]}}}^{2}\\
                 &= \trsp{\vecsix[]}\trsp{\partderiv[]\funtwo\paren*{\retr[\pt]\paren*{\tanvecone[\pt]}}}\paren*{\coeffsev^{2}\diffretrmat\paren*{\tanvecone[\pt]}\coordmetr[]\paren*{\pt}\diffretrmat\paren*{\tanvecone[\pt]} - \coordmetr[]\paren*{\pt}  }\partderiv[]\funtwo\paren*{\retr[\pt]\paren*{\tanvecone[\pt]}}\vecsix[].
            \end{split}
        \end{align}
        Since $\diffretrmat\paren*{\zerovec[\ptaccum]}=\id[{\tanspc[\ptaccum]\mani}]$ holds by \cref{eq:retrdiffzero}, we have 
        \begin{align}
            \coeffsev^{2}\diffretrmat\paren*{\zerovec[\ptaccum]}\coordmetr[]\paren*{\ptaccum}\diffretrmat\paren*{\zerovec[\ptaccum]} - \coordmetr[]\paren*{\ptaccum} = \paren*{\coeffsev^{2} - 1}\coordmetr[]\paren*{\ptaccum}\succ 0,
        \end{align}
        implying that $\coeffsev\Riemnorm[\ptaccum]{\sum_{\idxfou=1}^{\numfun}\vecsix[\idxfou]\gradstr\pullback[\ptaccum]{\funtwo}^{\idxfou}\paren*{\zerovec[\ptaccum]}} \geq \Riemnorm[\ptaccum]{\sum_{\idxfou=1}^{\numfun}\vecsix[\idxfou]\gradstr\funtwo[\idxfou]\paren*{\retr[\ptaccum]\paren*{\zerovec[\ptaccum]}}}$.
        Let $\subsetmaninin\subseteq\mani$ denote a closed neighborhood of $\ptaccum$.
        Since the functions $\retr[]$ and $\brc*{\funtwo[\idxfou]}_{\idxfou}$ are continuously differentiable and the set $\brc*{\paren*{\pt, \tanvecone[\pt]} \colon \pt \in\subsetmaninin, \tanvecone[\pt]\in\tanspc[\pt]\mani, \Riemnorm[\pt]{\tanvecone[\pt]}\leq \tholdvalten}$ is compact, we obtain the positive definiteness of $\coeffsev^{2}\diffretrmat\paren*{\tanvecone[\pt]}\coordmetr[]\paren*{\pt}\diffretrmat\paren*{\tanvecone[\pt]} - \coordmetr[]\paren*{\pt}$ for all $\pt\in\subsetmaninin$ by taking $\subsetmaninin$ sufficiently small and $\tholdvalten$ sufficiently small if necessary.
        This implies \cref{ineq:lincombgradpullbackgrad} for all $\pt \in \subsetmaninin$, any $\vecsix\in\setR[\numfun]$, and all $\tanvecone[\pt]\in\tanspc[\pt]\mani$ with $\Riemnorm[\pt]{\tanvecone[\pt]} \leq \tholdvalten$.
        The proof is complete.
        \qed
    \end{proof}

    }{}

\section{Proof of \cref{lemm:barrKKTvecflddirNewtonbound}}\label{appx:proofbarrKKTvecflddirNewtonbound}
    \begin{proof}
        Given $\allvar\in\mani\times\setR[\ineqdime]\times\setR[\eqdime]$ and $\barrparam[] > 0$, we define the operator $\linbarrKKTvecfld[\allvar]{\barrparam[]}\colon\tanspc[\allvar]\mani\to\tanspc[\pt]\mani\times\setR[\ineqdime]\times\setR[\eqdime]$ as 
        \begin{align}\label{def:linbarrKKTvecfld}
            &\linbarrKKTvecfld[\allvar]{\barrparam[]}\paren*{\tanvecone[\allvar]} \coloneqq
            \begin{bmatrix}
                \gradstr[]\pullback[\pt]{\objfun}\paren*{\tanvecone[\pt]} - \sum_{\ineqidx\in\ineqset} \paren*{\ineqLagmult[\ineqidx] + \tanvecone[{\ineqLagmult[\ineqidx]}]}\gradstr[]\pullback[\pt]{\ineqfun[\ineqidx]}\paren*{\tanvecone[\pt]} + \sum_{\eqidx\in\eqset} \paren*{\eqLagmult[\eqidx] + \tanvecone[{\eqLagmult[\eqidx]}]}\gradstr[]\pullback[\pt]{\eqfun[\eqidx]}\paren*{\tanvecone[\pt]}\\
                \Ineqfunmat[]\paren*{\retr[\pt]\paren*{\tanvecone[\pt]}}\paren*{\ineqLagmult[] + \tanvecone[{\ineqLagmult[]}]}- \barrparam[]\onevec\\
                \eqfun[]\paren*{\retr[\pt]\paren*{\tanvecone[\pt]}}
            \end{bmatrix}
        \end{align}
        for $\tanvecone[\allvar]=\paren*{\tanvecone[\pt], \tanvecone[{\ineqLagmult[]}], \tanvecone[{\eqLagmult[]}]}\in\tanspc[\allvar]\mani$.
        Note that, for any $\tanvecone[\allvar], \tanvectwo[\allvar]\in\tanspc[\allvar]\mani$, its directional derivative is
        \isextendedversion{\begin{align}\label{def:derivlinbarrKKTvecfld}
            &\D\linbarrKKTvecfld[\allvar]{\barrparam[]}\paren*{\tanvecone[\allvar]}\sbra*{\tanvectwo[\allvar]}= 
            \begin{bmatrix}
                A_{1}\\
                B_{1}\\
                C_{1}
            \end{bmatrix},
        \end{align}}{
        \begin{align}\label{def:derivlinbarrKKTvecfld}
            \D\linbarrKKTvecfld[\allvar]{\barrparam[]}\paren*{\tanvecone[\allvar]}\sbra*{\tanvectwo[\allvar]}= \sbra*{{A_{1}; B_{1}; C_{1}}},
        \end{align}}
        where
        \begin{align}
            A_{1} &\coloneqq \Hess\pullback[\pt]{\objfun}\paren*{\tanvecone[\pt]}\sbra*{\tanvectwo[\pt]} - \sum_{\ineqidx\in\ineqset}\paren*{\ineqLagmult[\ineqidx] + \tanvecone[{\ineqLagmult[\ineqidx]}]}\Hess\pullback[\pt]{\ineqfun[\ineqidx]}\paren*{\tanvecone[\pt]}\sbra*{\tanvectwo[\pt]} + \sum_{\eqidx\in\eqset}\paren*{\eqLagmult[\eqidx] + \tanvecone[{\eqLagmult[\eqidx]}]}\Hess\pullback[\pt]{\eqfun[\eqidx]}\paren*{\tanvecone[\pt]}\sbra*{\tanvectwo[\pt]}\\
            &\quad - \sum_{\ineqidx\in\ineqset} \tanvectwo[{\ineqLagmult[\ineqidx]}]\gradstr[]\pullback[\pt]{\ineqfun[\ineqidx]}\paren*{\tanvecone[\pt]} + \sum_{\eqidx\in\eqset} \tanvectwo[{\eqLagmult[\eqidx]}]\gradstr[]\pullback[\pt]{\eqfun[\eqidx]}\paren*{\tanvecone[\pt]},\\
            B_{1} &\coloneqq \sbra*{\paren*{\ineqLagmult[\ineqidx] + \tanvecone[{\ineqLagmult[\ineqidx]}]}\D\pullback[\pt]{\ineqfun[\ineqidx]}\paren*{\tanvecone[\pt]}\sbra*{\tanvectwo[\pt]} + \tanvectwo[{\ineqLagmult[\ineqidx]}]\pullback[\pt]{\ineqfun[\ineqidx]}\paren*{\tanvecone[\pt]}}_{\ineqidx=1,\ldots,\ineqdime}, \text{ and } C_{1} \coloneqq \sbra*{\D\pullback[\pt]{\eqfun[\eqidx]}\paren*{\tanvecone[\pt]}\sbra*{\tanvectwo[\pt]}}_{\eqidx=1,\ldots,\eqdime}.
        \end{align}
        From the continuous differentiability of $\linbarrKKTvecfld[\allvarotriter]{\barrparamotriter}$, we have
        \isextendedversion{
        \begin{align}
            \begin{split}\label{eq:linbarrKKTvecTaylor}
                \linbarrKKTvecfld[\allvarotriter]{\barrparamotriter}\paren*{\dirNewtonallvarbarrparamotriter} &= \linbarrKKTvecfld[\allvarotriter]{\barrparamotriter}\paren*{\zerovec[\allvarotriter]} + \D\linbarrKKTvecfld[\allvarotriter]{\barrparamotriter}\paren*{\zerovec[\allvarotriter]}\sbra*{\dirNewtonallvarbarrparamotriter}\\
                &\quad + \int_{0}^{1} \paren*{\D\linbarrKKTvecfld[\allvarotriter]{\barrparamotriter}\paren*{\tmeeig\dirNewtonallvarbarrparamotriter} - \D\linbarrKKTvecfld[\allvarotriter]{\barrparamotriter}\paren*{\zerovec[\allvarotriter]}}\sbra*{\dirNewtonallvarbarrparamotriter}\dd{\tmeeig}.
            \end{split}
        \end{align}}{
        \begin{equation}
            \begin{aligned}[t]\label{eq:linbarrKKTvecTaylor}
                &\linbarrKKTvecfld[\allvarotriter]{\barrparamotriter}\paren*{\dirNewtonallvarbarrparamotriter} \\
                &= \linbarrKKTvecfld[\allvarotriter]{\barrparamotriter}\paren*{\zerovec[\allvarotriter]} + \D\linbarrKKTvecfld[\allvarotriter]{\barrparamotriter}\paren*{\zerovec[\allvarotriter]}\sbra*{\dirNewtonallvarbarrparamotriter} + \int_{0}^{1} \paren*{\D\linbarrKKTvecfld[\allvarotriter]{\barrparamotriter}\paren*{\tmeeig\dirNewtonallvarbarrparamotriter} - \D\linbarrKKTvecfld[\allvarotriter]{\barrparamotriter}\paren*{\zerovec[\allvarotriter]}}\sbra*{\dirNewtonallvarbarrparamotriter}\dd{\tmeeig}.
            \end{aligned}
        \end{equation}}
        In the following, we analyze each term of \cref{eq:linbarrKKTvecTaylor}.
        For the first term, it follows from \isextendedversion{\cref{def:linbarrKKTvecfld,eq:barrKKTvecfld,eq:pullbackgradzero}}{\cref{def:linbarrKKTvecfld,eq:pullbackgradzero}, and the definition of $\barrKKTvecfld$} that
        \begin{align}\label{eq:linbarrKKTvecfldzerovec}
            \linbarrKKTvecfld[\allvarotriter]{\barrparamotriter}\paren*{\zerovec[\allvarotriter]} = \barrKKTvecfld\paren*{\allvarotriter;\barrparamotriter}.
        \end{align}
        As for the second term, substituting $\paren*{\allvar, \barrparam[], \tanvecone[\allvar], \tanvectwo[\allvar]}=\paren*{\allvarotriter, \barrparamotriter, \zerovec[\allvarotriter], \dirNewtonallvarbarrparamotriter}$ into \cref{def:derivlinbarrKKTvecfld} yields
        \begin{align}
              \begin{split}\label{eq:difflinbarrKKTvecfldzerovec}
                &\D\linbarrKKTvecfld[\allvarotriter]{\barrparamotriter}\paren*{\zerovec[\allvarotriter]}\sbra*{\dirNewtonallvarbarrparamotriter}\\
                &= 
                \left[\begin{array}{l}
                    \Hess\pullback[\ptotriter]{\objfun}\paren*{\zerovec[\ptotriter]}\sbra*{\dirNewtonptbarrparamotriter} - \sum_{\ineqidx\in\ineqset}\ineqLagmultotriter[\ineqidx]\Hess\pullback[\ptotriter]{\ineqfun[\ineqidx]}\paren*{\zerovec[\ptotriter]}\sbra*{\dirNewtonptbarrparamotriter} \\
                    \quad + \sum_{\eqidx\in\eqset}\eqLagmultotriter[\eqidx]\Hess\pullback[\ptotriter]{\eqfun[\eqidx]}\paren*{\zerovec[\ptotriter]}\sbra*{\dirNewtonptbarrparamotriter} - \sum_{\ineqidx\in\ineqset}\sbra*{\dirNewtonineqLagmultbarrparamotriter[]}_{\ineqidx}\gradstr\pullback[\ptotriter]{\ineqfun[\ineqidx]}\paren*{\zerovec[\ptotriter]}\\
                    \quad + \sum_{\eqidx\in\eqset}\sbra*{\dirNewtoneqLagmultbarrparamotriter[]}_{\eqidx}\gradstr\pullback[\ptotriter]{\eqfun[\eqidx]}\paren*{\zerovec[\ptotriter]}\\
                    \sbra*{\ineqLagmultotriter[\ineqidx] \D\pullback[\ptotriter]{\ineqfun[\ineqidx]}\paren*{\zerovec[\ptotriter]}\sbra*{\dirNewtonptbarrparamotriter} + \sbra*{\dirNewtonineqLagmultbarrparamotriter[]}_{\ineqidx}\pullback[\ptotriter]{\ineqfun[\ineqidx]}\paren*{\zerovec[\ptotriter]}}_{\ineqidx=1,\ldots,\ineqdime}\\
                    \sbra*{\D\pullback[\ptotriter]{\eqfun[\eqidx]}\paren*{\zerovec[\ptotriter]}\sbra*{\dirNewtonptbarrparamotriter}}_{\eqidx=1,\ldots,\eqdime}
                \end{array}\right]\\
                &=
                \left[\begin{array}{l}
                    \Hess[\pt]\Lagfun\paren*{\allvarotriter}\sbra*{\dirNewtonptbarrparamotriter} - \sum_{\ineqidx\in\ineqset}\sbra*{\dirNewtonineqLagmultbarrparamotriter[]}_{\ineqidx}\gradstr\ineqfun[\ineqidx]\paren*{\ptotriter} + \sum_{\eqidx\in\eqset}\sbra*{\dirNewtoneqLagmultbarrparamotriter[]}_{\eqidx}\gradstr\eqfun[\eqidx]\paren*{\ptotriter} \\
                    \sbra*{\ineqLagmultotriter[\ineqidx] \metr[\ptotriter]{\gradstr\ineqfun[\ineqidx]\paren*{\ptotriter}}{\dirNewtonptbarrparamotriter} + \sbra*{\dirNewtonineqLagmultbarrparamotriter[]}_{\ineqidx}\ineqfun[\ineqidx]\paren*{\ptotriter}}_{\ineqidx=1,\ldots,\ineqdime}\\
                    \sbra*{\metr[\ptotriter]{\gradstr\eqfun[\eqidx]\paren*{\ptotriter}}{\dirNewtonptbarrparamotriter}}_{\eqidx = 1,\ldots,\eqdime}
                \end{array}\right]\\
                &= \Jacobian[\barrKKTvecfld]\paren*{\allvarotriter}\sbra*{\dirNewtonallvarbarrparamotriter[]},  
            \end{split}
        \end{align}
        where the second equality follows from \isextendedversion{\cref{eq:secondretrHesszero,eq:riemgrad,eq:pullbackgradzero,eq:retrzero,eq:retrdiffzero}}{\cref{eq:riemgrad}, \cref{eq:pullbackgradzero}, the definition of the retraction, and \cref{prop:secodrretrHess}} and the third one from \cref{def:JacobbarrKKTvecfld}.
        Next, we consider the third term.
        Note that, by the twice continuous differentiability of $\brc*{\ineqfun[\ineqidx]}_{\ineqidx\in\ineqset}$ and $\brc*{\eqfun[\eqidx]}_{\eqidx\in\eqset}$ and \cite[Lemma~10.57]{Boumal23IntroOptimSmthMani}, there exist positive scalars $\brc*{\Lipschitzconstone[{\ineqfun[\ineqidx]}]}_{\ineqidx\in\ineqset}$, 
        $\brc*{\Lipschitzconstsix[{\ineqfun[\ineqidx]}]}_{\ineqidx\in\ineqset}$, and $\brc*{\Lipschitzconstsix[{\eqfun[\eqidx]}]}_{\eqidx\in\eqset}$ such that, for all $\pt\in\mani$ sufficiently close to $\ptaccum$ and any $\tanvecone[\pt]\in\tanspc[\pt]\mani$ sufficiently small,
        \isextendedversion{
        \begin{align}
            &\abs*{\pullback[\pt]{\ineqfun[\ineqidx]}\paren*{\tanvecone[\pt]} - \pullback[\pt]{\ineqfun[\ineqidx]}\paren*{\zerovec[\pt]}} \leq \Lipschitzconstone[{\ineqfun[\ineqidx]}]\Riemdist{\retr[\pt]\paren*{\tanvecone[\pt]}}{\pt} \text{ for all } \ineqidx \in \ineqset,\label{ineq:ineqfunNewtonLipschitzcont}\\
            &\Riemnorm[]{\gradstr\pullback[\pt]{\ineqfun[\ineqidx]}\paren*{\tanvecone[\pt]} - \gradstr\pullback[\pt]{\ineqfun[\ineqidx]}\paren*{\zerovec[\pt]}} \leq \Lipschitzconstsix[{\ineqfun[\ineqidx]}]\Riemnorm[\pt]{\tanvecone[\pt]} \text{ for all } \ineqidx \in \ineqset.\label{ineq:pullbackgradineqfunLipschitzcont}\\
            &\Riemnorm[]{\gradstr\pullback[\pt]{\eqfun[\eqidx]}\paren*{\tanvecone[\pt]} - \gradstr\pullback[\pt]{\eqfun[\eqidx]}\paren*{\zerovec[\pt]}} \leq \Lipschitzconstsix[{\eqfun[\eqidx]}]\Riemnorm[\pt]{\tanvecone[\pt]} \text{ for all } \eqidx \in \eqset.\label{ineq:pullbackgradeqfunLipschitzcont}
        \end{align}}{
        \begin{align}
            \begin{split}\label{ineq:pullbackgradeqineqfunLipschitzcont}
                &\abs*{\pullback[\pt]{\ineqfun[\ineqidx]}\paren*{\tanvecone[\pt]} - \pullback[\pt]{\ineqfun[\ineqidx]}\paren*{\zerovec[\pt]}} \leq \Lipschitzconstone[{\ineqfun[\ineqidx]}]\Riemdist{\retr[\pt]\paren*{\tanvecone[\pt]}}{\pt} \text{ for all } \ineqidx \in \ineqset,\\
                &\Riemnorm[]{\gradstr\pullback[\pt]{\ineqfun[\ineqidx]}\paren*{\tanvecone[\pt]} - \gradstr\pullback[\pt]{\ineqfun[\ineqidx]}\paren*{\zerovec[\pt]}} \leq \Lipschitzconstsix[{\ineqfun[\ineqidx]}]\Riemnorm[\pt]{\tanvecone[\pt]} \text{ for all } \ineqidx \in \ineqset,\label{ineq:pullbackgradineqfunLipschitzcont}\\
                &\Riemnorm[]{\gradstr\pullback[\pt]{\eqfun[\eqidx]}\paren*{\tanvecone[\pt]} - \gradstr\pullback[\pt]{\eqfun[\eqidx]}\paren*{\zerovec[\pt]}} \leq \Lipschitzconstsix[{\eqfun[\eqidx]}]\Riemnorm[\pt]{\tanvecone[\pt]} \text{ for all } \eqidx \in \eqset.\label{ineq:pullbackgradeqfunLipschitzcont}
            \end{split}
        \end{align}}
        It also follows from \cref{def:derivlinbarrKKTvecfld} with $\paren*{\allvar, \barrparam[], \tanvecone[\allvar], \tanvectwo[\allvar]}$ replaced by $\paren*{\allvarotriter, \barrparamotriter, \tmeeig\dirNewtonallvarbarrparamotriter, \dirNewtonallvarbarrparamotriter}$ and \cref{eq:difflinbarrKKTvecfldzerovec} that
        \isextendedversion{
        \begin{align}\label{eq:diffderivlinbarrKKTvecfld}
                &\paren*{\D\linbarrKKTvecfld[\allvarotriter]{\barrparamotriter}\paren*{\tmeeig\dirNewtonallvarbarrparamotriter} - \D\linbarrKKTvecfld[\allvarotriter]{\barrparamotriter}\paren*{\zerovec[\allvarotriter]}}\sbra*{\dirNewtonallvarbarrparamotriter} = \begin{bmatrix}
                    A_{2}\\
                    B_{2}\\
                    C_{2}
                \end{bmatrix},
        \end{align}}
        {$\paren*{\D\linbarrKKTvecfld[\allvarotriter]{\barrparamotriter}\paren*{\tmeeig\dirNewtonallvarbarrparamotriter} - \D\linbarrKKTvecfld[\allvarotriter]{\barrparamotriter}\paren*{\zerovec[\allvarotriter]}}\sbra*{\dirNewtonallvarbarrparamotriter} = \sbra*{A_{2}; B_{2}; C_{2}}$,}
        where
        \isextendedversion{
        \begin{align}
            &A_{2} \coloneqq \paren*{\Hess\pullback[\ptotriter]{\objfun}\paren*{\tmeeig\dirNewtonptbarrparamotriter} - \Hess\pullback[\ptotriter]{\objfun}\paren*{\zerovec[\ptotriter]}}\sbra*{\dirNewtonptbarrparamotriter}\\
            &\quad - \sum_{\ineqidx\in\ineqset}\paren*{\ineqLagmultotriter[\ineqidx] + \tmeeig\sbra*{\dirNewtonineqLagmultbarrparamotriter[]}_{\ineqidx}} \paren*{\Hess\pullback[\ptotriter]{\ineqfun[\ineqidx]}\paren*{\tmeeig\dirNewtonptbarrparamotriter} - \Hess\pullback[\ptotriter]{\ineqfun[\ineqidx]}\paren*{\zerovec[\ptotriter]}}\sbra*{\dirNewtonptbarrparamotriter}\\
            &\quad + \sum_{\eqidx\in\eqset}\paren*{\eqLagmultotriter[\eqidx] + \tmeeig\sbra*{\dirNewtoneqLagmultbarrparamotriter[]}_{\eqidx}} \paren*{\Hess\pullback[\ptotriter]{\eqfun[\eqidx]}\paren*{\tmeeig\dirNewtonptbarrparamotriter} - \Hess\pullback[\ptotriter]{\eqfun[\eqidx]}\paren*{\zerovec[\ptotriter]}}\sbra*{\dirNewtonptbarrparamotriter}\\
            &\quad- \sum_{\ineqidx\in\ineqset} \tmeeig\sbra*{\dirNewtonineqLagmultbarrparamotriter[]}_{\ineqidx}\Hess\pullback[\ptotriter]{\ineqfun[\ineqidx]}\paren*{\zerovec[\ptotriter]}\sbra*{\dirNewtonptbarrparamotriter}\\
            &\quad - \sum_{\ineqidx\in\ineqset}\sbra*{\dirNewtonineqLagmultbarrparamotriter[]}_{\ineqidx}\paren*{\gradstr\pullback[\ptotriter]{\ineqfun[\ineqidx]}\paren*{\tmeeig\dirNewtonptbarrparamotriter} - \gradstr\pullback[\ptotriter]{\ineqfun[\ineqidx]}\paren*{\zerovec[\ptotriter]}}\\
            &\quad + \sum_{\eqidx\in\eqset} \tmeeig\sbra*{\dirNewtoneqLagmultbarrparamotriter[]}_{\eqidx}\Hess\pullback[\ptotriter]{\eqfun[\eqidx]}\paren*{\zerovec[\ptotriter]}\sbra*{\dirNewtonptbarrparamotriter}\\
            &\quad + \sum_{\eqidx\in\eqset}\sbra*{\dirNewtoneqLagmultbarrparamotriter[]}_{\eqidx}\paren*{\gradstr\pullback[\ptotriter]{\eqfun[\eqidx]}\paren*{\tmeeig\dirNewtonptbarrparamotriter} - \gradstr\pullback[\ptotriter]{\eqfun[\eqidx]}\paren*{\zerovec[\ptotriter]}},\\
            &B_{2} \coloneqq \sbra*{\paren*{\ineqLagmultotriter[\ineqidx] + \tmeeig\sbra*{\dirNewtonineqLagmultbarrparamotriter}_{\ineqidx}}\metr[\ptotriter]{\gradstr\pullback[\ptotriter]{\ineqfun[\ineqidx]}\paren*{\tmeeig\dirNewtonptbarrparamotriter} - \gradstr\pullback[\ptotriter]{\ineqfun[\ineqidx]}\paren*{\zerovec[\ptotriter]}}{\dirNewtonptbarrparamotriter}\\
            &\quad  + \tmeeig\sbra*{\dirNewtonineqLagmultbarrparamotriter}_{\ineqidx}\metr[\ptotriter]{\gradstr\pullback[\ptotriter]{\ineqfun[\ineqidx]}\paren*{\zerovec[\ptotriter]}}{\dirNewtonptbarrparamotriter} \\
            &\quad + \sbra*{\dirNewtonineqLagmultbarrparamotriter[]}_{\ineqidx}\paren*{\pullback[\ptotriter]{\ineqfun[\ineqidx]}\paren*{\tmeeig\dirNewtonptbarrparamotriter} - \pullback[\ptotriter]{\ineqfun[\ineqidx]}\paren*{\zerovec[\ptotriter]}}}_{\ineqidx=1,\ldots,\ineqdime},\\
            &C_{2} \coloneqq \sbra*{\metr[\ptotriter]{\gradstr\pullback[\ptotriter]{\eqfun[\eqidx]}\paren*{\tmeeig\dirNewtonptbarrparamotriter} - \gradstr\pullback[\ptotriter]{\eqfun[\eqidx]}\paren*{\zerovec[\ptotriter]}}{\dirNewtonptbarrparamotriter}}_{\eqidx=1,\ldots,\eqdime}.
        \end{align}}
        {
        \begin{align}
            &A_{2} \coloneqq \paren*{\Hess\pullback[\ptotriter]{\objfun}\paren*{\tmeeig\dirNewtonptbarrparamotriter} - \Hess\pullback[\ptotriter]{\objfun}\paren*{\zerovec[\ptotriter]}}\sbra*{\dirNewtonptbarrparamotriter}\\
            &\quad - \sum_{\ineqidx\in\ineqset}\paren*{\ineqLagmultotriter[\ineqidx] + \tmeeig\sbra*{\dirNewtonineqLagmultbarrparamotriter[]}_{\ineqidx}}\paren*{\Hess\pullback[\ptotriter]{\ineqfun[\ineqidx]}\paren*{\tmeeig\dirNewtonptbarrparamotriter} - \Hess\pullback[\ptotriter]{\ineqfun[\ineqidx]}\paren*{\zerovec[\ptotriter]}}\sbra*{\dirNewtonptbarrparamotriter}\\
            &\quad + \sum_{\eqidx\in\eqset}\paren*{\eqLagmultotriter[\eqidx] + \tmeeig\sbra*{\dirNewtoneqLagmultbarrparamotriter[]}_{\eqidx}} \paren*{\Hess\pullback[\ptotriter]{\eqfun[\eqidx]}\paren*{\tmeeig\dirNewtonptbarrparamotriter} - \Hess\pullback[\ptotriter]{\eqfun[\eqidx]}\paren*{\zerovec[\ptotriter]}}\sbra*{\dirNewtonptbarrparamotriter}\\
            &\quad - \sum_{\ineqidx\in\ineqset} \tmeeig\sbra*{\dirNewtonineqLagmultbarrparamotriter[]}_{\ineqidx}\Hess\pullback[\ptotriter]{\ineqfun[\ineqidx]}\paren*{\zerovec[\ptotriter]}\sbra*{\dirNewtonptbarrparamotriter} \\
            &\quad - \sum_{\ineqidx\in\ineqset}\sbra*{\dirNewtonineqLagmultbarrparamotriter[]}_{\ineqidx}\paren*{\gradstr\pullback[\ptotriter]{\ineqfun[\ineqidx]}\paren*{\tmeeig\dirNewtonptbarrparamotriter} - \gradstr\pullback[\ptotriter]{\ineqfun[\ineqidx]}\paren*{\zerovec[\ptotriter]}}\\
            &\quad + \sum_{\eqidx\in\eqset} \tmeeig\sbra*{\dirNewtoneqLagmultbarrparamotriter[]}_{\eqidx}\Hess\pullback[\ptotriter]{\eqfun[\eqidx]}\paren*{\zerovec[\ptotriter]}\sbra*{\dirNewtonptbarrparamotriter}\\
            &\quad + \sum_{\eqidx\in\eqset}\sbra*{\dirNewtoneqLagmultbarrparamotriter[]}_{\eqidx}\paren*{\gradstr\pullback[\ptotriter]{\eqfun[\eqidx]}\paren*{\tmeeig\dirNewtonptbarrparamotriter} - \gradstr\pullback[\ptotriter]{\eqfun[\eqidx]}\paren*{\zerovec[\ptotriter]}},\\
            &B_{2} \coloneqq \sbra*{\paren*{\ineqLagmultotriter[\ineqidx] + \tmeeig\sbra*{\dirNewtonineqLagmultbarrparamotriter}_{\ineqidx}}\metr[\ptotriter]{\gradstr\pullback[\ptotriter]{\ineqfun[\ineqidx]}\paren*{\tmeeig\dirNewtonptbarrparamotriter} - \gradstr\pullback[\ptotriter]{\ineqfun[\ineqidx]}\paren*{\zerovec[\ptotriter]}}{\dirNewtonptbarrparamotriter}\\
            &\quad + \tmeeig\sbra*{\dirNewtonineqLagmultbarrparamotriter}_{\ineqidx}\metr[\ptotriter]{\gradstr\pullback[\ptotriter]{\ineqfun[\ineqidx]}\paren*{\zerovec[\ptotriter]}}{\dirNewtonptbarrparamotriter}\\
            &\quad + \sbra*{\dirNewtonineqLagmultbarrparamotriter[]}_{\ineqidx}\paren*{\pullback[\ptotriter]{\ineqfun[\ineqidx]}\paren*{\tmeeig\dirNewtonptbarrparamotriter} - \pullback[\ptotriter]{\ineqfun[\ineqidx]}\paren*{\zerovec[\ptotriter]}}}_{\ineqidx=1,\ldots,\ineqdime},\\
            &C_{2} \coloneqq \sbra*{\metr[\ptotriter]{\gradstr\pullback[\ptotriter]{\eqfun[\eqidx]}\paren*{\tmeeig\dirNewtonptbarrparamotriter} - \gradstr\pullback[\ptotriter]{\eqfun[\eqidx]}\paren*{\zerovec[\ptotriter]}}{\dirNewtonptbarrparamotriter}}_{\eqidx=1,\ldots,\eqdime}.
        \end{align}}
        Let $\coeffeig[a], \coeffeig[b], \coeffeig[c] > 0$ be sufficiently large.
        Note that $\brc*{\abs*{\ineqLagmultotriter[\ineqidx]}}_{\ineqidx}, \brc*{\abs*{\eqLagmultotriter[\eqidx]}}_{\eqidx}, \brc*{\abs*{\sbra*{\dirNewtonineqLagmultbarrparamotriter[]}_{\ineqidx}}}_{\ineqidx}$, and $\brc*{\abs*{\sbra*{\dirNewtoneqLagmultbarrparamotriter[]}_{\eqidx}}}_{\eqidx}$ are bounded around $\allvaraccum$ and $\barrparam[]=0$.
        Thus, for any $\otriteridx$-th iteration of \cref{algo:RIPMOuter} satisfying that $\allvarotriter$ is sufficiently close to $\allvaraccum$ and $\barrparamotriterm$ is sufficiently small, we have
            \begin{align}
                &\Riemnorm[\ptotriter]{A_{2}} \leq \opnorm{\Hess\pullback[\ptotriter]{\objfun}\paren*{\tmeeig\dirNewtonptbarrparamotriter} - \Hess\pullback[\ptotriter]{\objfun}\paren*{\zerovec[\ptotriter]}}\Riemnorm[\ptotriter]{\dirNewtonptbarrparamotriter}\notag\\
                &\quad + \sum_{\ineqidx\in\ineqset}\abs*{\ineqLagmultotriter[\ineqidx] + \tmeeig\sbra*{\dirNewtonineqLagmultbarrparamotriter[]}_{\ineqidx}}\opnorm{\Hess\pullback[\ptotriter]{\ineqfun[\ineqidx]}\paren*{\tmeeig\dirNewtonptbarrparamotriter} - \Hess\pullback[\ptotriter]{\ineqfun[\ineqidx]}\paren*{\zerovec[\ptotriter]}}\Riemnorm[\ptotriter]{\dirNewtonptbarrparamotriter} \notag\\
                &\quad+ \sum_{\eqidx\in\eqset}\abs*{\eqLagmultotriter[\eqidx] + \tmeeig\sbra*{\dirNewtoneqLagmultbarrparamotriter[]}_{\eqidx}}\opnorm{\Hess\pullback[\ptotriter]{\eqfun[\eqidx]}\paren*{\tmeeig\dirNewtonptbarrparamotriter} - \Hess\pullback[\ptotriter]{\eqfun[\eqidx]}\paren*{\zerovec[\ptotriter]}}\Riemnorm[\ptotriter]{\dirNewtonptbarrparamotriter}\notag\\
                &\quad+ \sum_{\ineqidx\in\ineqset}\tmeeig\abs*{\sbra*{\dirNewtonineqLagmultbarrparamotriter[]}_{\ineqidx}}\opnorm{\Hess\pullback[\ptotriter]{\ineqfun[\ineqidx]}\paren*{\zerovec[\ptotriter]}}\Riemnorm[\ptotriter]{\dirNewtonptbarrparamotriter} \notag\\
                &\quad + \sum_{\ineqidx\in\ineqset}\abs*{\sbra*{\dirNewtonineqLagmultbarrparamotriter[]}_{\ineqidx}}\Riemnorm[]{ \gradstr\pullback[\ptotriter]{\ineqfun[\ineqidx]}\paren*{\tmeeig\dirNewtonptbarrparamotriter} - \gradstr\pullback[\ptotriter]{\ineqfun[\ineqidx]}\paren*{\zerovec[\ptotriter]}}\notag\\
                &\quad + \sum_{\eqidx\in\eqset}\tmeeig\abs*{\sbra*{\dirNewtoneqLagmultbarrparamotriter[]}_{\eqidx}}\opnorm{\Hess\pullback[\ptotriter]{\eqfun[\eqidx]}\paren*{\zerovec[\ptotriter]}}\Riemnorm[\ptotriter]{\dirNewtonptbarrparamotriter}\notag\\
                &\quad + \sum_{\eqidx\in\eqset}\abs*{\sbra*{\dirNewtoneqLagmultbarrparamotriter[]}_{\eqidx}}\Riemnorm[]{ \gradstr\pullback[\ptotriter]{\eqfun[\eqidx]}\paren*{\tmeeig\dirNewtonptbarrparamotriter} - \gradstr\pullback[\ptotriter]{\eqfun[\eqidx]}\paren*{\zerovec[\ptotriter]}}\notag\\
                &\leq\Lipschitzconstfiv[\objfun] \Riemnorm[\ptotriter]{\dirNewtonptbarrparamotriter}^{2} + \sum_{\ineqidx\in\ineqset}\Lipschitzconstfiv[{\ineqfun[\ineqidx]}] \paren*{\abs*{\ineqLagmultotriter[\ineqidx]} + \abs*{\sbra*{\dirNewtonineqLagmultbarrparamotriter[]}_{\ineqidx}}}\Riemnorm[\ptotriter]{\dirNewtonptbarrparamotriter}^{2}\notag\\
                &\quad + \sum_{\eqidx\in\eqset}\Lipschitzconstfiv[{\eqfun[\eqidx]}] \paren*{\abs*{\eqLagmultotriter[\eqidx]} + \abs*{\sbra*{\dirNewtoneqLagmultbarrparamotriter[]}_{\eqidx}}}\Riemnorm[\ptotriter]{\dirNewtonptbarrparamotriter}^{2} + \sum_{\ineqidx\in\ineqset}\abs*{\sbra*{\dirNewtonineqLagmultbarrparamotriter[]}_{\ineqidx}}\opnorm{\Hess\ineqfun[\ineqidx]\paren*{\ptotriter}}\Riemnorm[\ptotriter]{\dirNewtonptbarrparamotriter}\notag\\
                &\quad + \sum_{\ineqidx\in\ineqset}\Lipschitzconstsix[{\ineqfun[\ineqidx]}]\abs*{\sbra*{\dirNewtonineqLagmultbarrparamotriter[]}_{\ineqidx}} \Riemnorm[\ptotriter]{\dirNewtonptbarrparamotriter} + \sum_{\eqidx\in\eqset}\abs*{\sbra*{\dirNewtoneqLagmultbarrparamotriter[]}_{\eqidx}}\opnorm{\Hess\eqfun[\eqidx]\paren*{\ptotriter}}\Riemnorm[\ptotriter]{\dirNewtonptbarrparamotriter}\notag\\
                &\quad+ \sum_{\eqidx\in\eqset}\Lipschitzconstsix[{\eqfun[\eqidx]}] 
                \abs*{\sbra*{\dirNewtoneqLagmultbarrparamotriter[]}_{\eqidx}} \Riemnorm[\ptotriter]{\dirNewtonptbarrparamotriter}\leq \coeffeig[a]\Riemnorm[\ptotriter]{\dirNewtonallvarbarrparamotriter}^{2},\label{ineq:derivlinbarKKTvecflddirNewtonA}
            \end{align}
            where the first inequality follows from the definition of $A_{2}$ and the second one from \isextendedversion{\cref{ineq:pullbackHessobjfunLipschitz,ineq:pullbackHessineqfunLipschitz,ineq:pullbackHesseqfunLipschitz}, 
            \cref{eq:secondretrHesszero}, \cref{ineq:pullbackgradineqfunLipschitzcont,ineq:pullbackgradeqfunLipschitzcont}
            with $\tanvecone[\pt]$ replaced by $\tmeeig\dirNewtonptbarrparamotriter$, and $\tmeeig \leq 1$,}{\cref{assu:pullbackHessobjineqfunLipschitz}, 
            \cref{prop:secodrretrHess}, \cref{ineq:pullbackgradeqineqfunLipschitzcont}
            with $\tanvecone[\pt]$ replaced by $\tmeeig\dirNewtonptbarrparamotriter$, and $\tmeeig \leq 1$,} the third one from $\abs*{\sbra*{\dirNewtonineqLagmultbarrparamotriter[]}_{\ineqidx}} \leq \Riemnorm[\allvarotriter]{\dirNewtonallvarbarrparamotriter}$, $\abs*{\sbra*{\dirNewtoneqLagmultbarrparamotriter[]}_{\eqidx}}\leq \Riemnorm[\allvarotriter]{\dirNewtonallvarbarrparamotriter}$, $\Riemnorm[\ptotriter]{\dirNewtonptbarrparamotriter} \leq \Riemnorm[\allvarotriter]{\dirNewtonallvarbarrparamotriter}$, 
            the boundedness of $\brc*{\abs*{\ineqLagmultotriter[\ineqidx]}}_{\ineqidx\in\ineqset}, \brc*{\abs*{\eqLagmultotriter[\eqidx]}}_{\eqidx\in\eqset}, \brc*{\abs*{\sbra*{\dirNewtonineqLagmultbarrparamotriter[]}_{\ineqidx}}}_{\ineqidx\in\ineqset}$, and $\brc*{\abs*{\sbra*{\dirNewtoneqLagmultbarrparamotriter[]}_{\eqidx}}}_{\eqidx\in\eqset}$, and the continuities of $\brc*{\opnorm{\Hess\ineqfun[\ineqidx]\paren*{\cdot}}}_{\ineqidx\in\ineqset}$ and $\brc*{\opnorm{\Hess\eqfun[\eqidx]\paren*{\cdot}}}_{\eqidx\in\eqset}$.
            Similarly, for such $\otriteridx$-th iteration of \cref{algo:RIPMOuter}, we have
            \begin{align}
            \begin{split}\label{ineq:derivlinbarKKTvecflddirNewtonB}
                 &\norm*{B_{2}} \leq \sum_{\ineqidx\in\ineqset}\paren*{\ineqLagmultotriter[\ineqidx] + \tmeeig\sbra*{\dirNewtonineqLagmultbarrparamotriter[]}_{\ineqidx}}\Riemnorm[\ptotriter]{\gradstr\pullback[\ptotriter]{\ineqfun[\ineqidx]}\paren*{\tmeeig\dirNewtonptbarrparamotriter} - \gradstr\pullback[\ptotriter]{\ineqfun[\ineqidx]}\paren*{\zerovec[\ptotriter]}}\Riemnorm[\ptotriter]{\dirNewtonptbarrparamotriter}\\
                 &\quad + \sum_{\ineqidx\in\ineqset}\tmeeig\abs*{\sbra*{\dirNewtonineqLagmultbarrparamotriter[]}_{\ineqidx}\metr[\ptotriter]{\gradstr\pullback[\ptotriter]{\ineqfun[\ineqidx]}\paren*{\zerovec[\ptotriter]}}{\dirNewtonptbarrparamotriter}}\\
                 &\quad+ \sum_{\ineqidx\in\ineqset}\abs*{\sbra*{\dirNewtonineqLagmultbarrparamotriter[]}_{\ineqidx}\pullback[\ptotriter]{\ineqfun[\ineqidx]}\paren*{\tmeeig\dirNewtonptbarrparamotriter} - \pullback[\ptotriter]{\ineqfun[\ineqidx]}\paren*{\zerovec[\ptotriter]}}\\
                &\leq \sum_{\ineqidx\in\ineqset}\Lipschitzconstsix[{\ineqfun[\ineqidx]}]
                \paren*{\abs*{\ineqLagmultotriter[\ineqidx]} + \abs*{\sbra*{\dirNewtonineqLagmultbarrparamotriter[]}_{\ineqidx}}}\Riemnorm[\ptotriter]{\dirNewtonptbarrparamotriter}^{2} + \sum_{\ineqidx\in\ineqset}\abs*{\sbra*{\dirNewtonineqLagmultbarrparamotriter[]}_{\ineqidx}} \Riemnorm[\ptotriter]{\gradstr\ineqfun[\ineqidx]\paren*{\ptotriter}}\Riemnorm[\ptotriter]{\dirNewtonptbarrparamotriter} \\
                &\quad+ \sum_{\ineqidx\in\ineqset}\Lipschitzconstone[{\ineqfun[\ineqidx]}]\abs*{\sbra*{\dirNewtonineqLagmultbarrparamotriter[]}_{\ineqidx}} \Riemdist{\retr[\ptotriter]\paren*{\tmeeig\dirNewtonptbarrparamotriter}}{\ptotriter} \leq \coeffeig[b]\Riemnorm[\ptotriter]{\dirNewtonallvarbarrparamotriter}^{2},
            \end{split}
            \end{align}
            where the first inequality follows from the definition of $B_{2}$, the second one from \isextendedversion{
            \cref{ineq:ineqfunNewtonLipschitzcont,ineq:pullbackgradineqfunLipschitzcont}
            with $\tanvecone[\pt]$ replaced by $\tmeeig\dirNewtonptbarrparamotriter$, $\tmeeig \leq 1$ and \cref{eq:pullbackgradzero},}{
            \cref{ineq:pullbackgradeqineqfunLipschitzcont}
            with $\tanvecone[\pt]$ replaced by $\tmeeig\dirNewtonptbarrparamotriter$, $\tmeeig \leq 1$, and \cref{eq:pullbackgradzero},} and the third one from $\abs*{\sbra*{\dirNewtonineqLagmultbarrparamotriter[]}_{\ineqidx}} \leq \Riemnorm[\allvarotriter]{\dirNewtonallvarbarrparamotriter}$, $\Riemnorm[\ptotriter]{\dirNewtonptbarrparamotriter} \leq \Riemnorm[\allvarotriter]{\dirNewtonallvarbarrparamotriter}$, the boundedness of 
            $\brc*{\abs*{\ineqLagmultotriter[\ineqidx]}}_{\ineqidx\in\ineqset}$ and $\brc*{\abs*{\sbra*{\dirNewtonineqLagmultbarrparamotriter[]}_{\ineqidx}}}_{\ineqidx\in\ineqset}$, and the continuity of $\brc*{\Riemnorm[\cdot]{\gradstr\ineqfun[\ineqidx]\paren*{\cdot}}}_{\ineqidx\in\ineqset}$, and \isextendedversion{\cref{lemm:retrtanvecequiv}}{\cite[Lemma~2]{HuangAbsilGallivan2015RiemSymRankOneTRM}}.
            We also have, for such $\otriteridx$-th iteration of \cref{algo:RIPMOuter},
            \begin{align}
            \begin{split}\label{ineq:derivlinbarKKTvecflddirNewtonC}
                \norm*{C_{2}} &\leq \sum_{\eqidx\in\eqset}\Riemnorm[\ptotriter]{\gradstr\pullback[\ptotriter]{\eqfun[\eqidx]}\paren*{\tmeeig\dirNewtonptbarrparamotriter} - \gradstr\pullback[\ptotriter]{\eqfun[\eqidx]}\paren*{\zerovec[\ptotriter]}}\Riemnorm[\ptotriter]{\dirNewtonptbarrparamotriter} \\
                &\leq \sum_{\eqidx\in\eqset}\Lipschitzconstsix[{\eqfun[\eqidx]}]\Riemnorm[\ptotriter]{\dirNewtonptbarrparamotriter}^{2} \leq \coeffeig[c]\Riemnorm[\ptotriter]{\dirNewtonallvarbarrparamotriter}^{2}, 
            \end{split}
        \end{align}
        where the first inequality follows from the definition of $C_{2}$, the second one from \isextendedversion{\cref{ineq:pullbackgradeqfunLipschitzcont} with $\tanvecone[\pt]$ replaced by $\tmeeig\dirNewtonptbarrparamotriter$}{\cref{ineq:pullbackgradeqineqfunLipschitzcont} with $\tanvecone[\pt]$ replaced by $\tmeeig\dirNewtonptbarrparamotriter$}, and the third one from $\Riemnorm[\ptotriter]{\dirNewtonptbarrparamotriter} \leq \Riemnorm[\allvarotriter]{\dirNewtonallvarbarrparamotriter}$.
        From \cref{ineq:derivlinbarKKTvecflddirNewtonA,ineq:derivlinbarKKTvecflddirNewtonB,ineq:derivlinbarKKTvecflddirNewtonC}, we have
        \begin{equation}\label{ineq:derivlinbarKKTvecflddirNewton}
            \Riemnorm[\allvarotriter]{\paren*{\D\linbarrKKTvecfld[\allvarotriter]{\barrparamotriter}\paren*{\tmeeig\dirNewtonallvarbarrparamotriter} - \D\linbarrKKTvecfld[\allvarotriter]{\barrparamotriter}\paren*{\zerovec[\allvarotriter]}}\sbra*{\dirNewtonallvarbarrparamotriter}}\leq \paren*{\coeffeig[a] + \coeffeig[b] + \coeffeig[c]}\Riemnorm[\allvarotriter]{\dirNewtonallvarbarrparamotriter}^{2}\leq \coeffeig[]\constdirNewton^{2}\barrparamotriterm^{2},
        \end{equation}
        where $\coeffeig[]\coloneqq\paren*{\coeffeig[a] + \coeffeig[b] + \coeffeig[c]} > 0$ and the second inequality follows from \cref{lemm:dirNewtonallvarbarrparambound}.
        Hence, from \cref{eq:linbarrKKTvecTaylor}, we obtain
        \isextendedversion{
        \begin{equation}
            \begin{aligned}[t]\label{ineq:linbarrKKTvecflddirNewtonboundbarrparam}
                &\Riemnorm[\allvarotriter]{\linbarrKKTvecfld[\allvarotriter]{\barrparamotriter}\paren*{\dirNewtonallvarbarrparamotriter}}\\
                &= \Riemnorm[\allvarotriter]{\barrKKTvecfld\paren*{\allvarotriter;\barrparamotriter} + \D\linbarrKKTvecfld[\allvarotriter]{\barrparamotriter}\paren*{\zerovec[\allvarotriter]}\sbra*{\dirNewtonallvarbarrparamotriter}
                + \int_{0}^{1} \paren*{\D\linbarrKKTvecfld[\allvarotriter]{\barrparamotriter}\paren*{\tmeeig\dirNewtonallvarbarrparamotriter} - \D\linbarrKKTvecfld[\allvarotriter]{\barrparamotriter}\paren*{\zerovec[\allvarotriter]}}\sbra*{\dirNewtonallvarbarrparamotriter}\dd{\tmeeig}}\\
                &\leq \int_{0}^{1} \Riemnorm[\allvarotriter]{\paren*{\D\linbarrKKTvecfld[\allvarotriter]{\barrparamotriter}\paren*{\tmeeig\dirNewtonallvarbarrparamotriter} - \D\linbarrKKTvecfld[\allvarotriter]{\barrparamotriter}\paren*{\zerovec[\allvarotriter]}}\sbra*{\dirNewtonallvarbarrparamotriter}}\dd{\tmeeig}\leq \coeffeig[]\constdirNewton^{2}\barrparamotriterm^{2},
            \end{aligned}
        \end{equation}}
        {\begin{align}
            &\Riemnorm[\allvarotriter]{\linbarrKKTvecfld[\allvarotriter]{\barrparamotriter}\paren*{\dirNewtonallvarbarrparamotriter}}\notag\\
            &= \Riemnorm[\allvarotriter]{\barrKKTvecfld\paren*{\allvarotriter;\barrparamotriter} + \D\linbarrKKTvecfld[\allvarotriter]{\barrparamotriter}\paren*{\zerovec[\allvarotriter]}\sbra*{\dirNewtonallvarbarrparamotriter}
            + \int_{0}^{1} \paren*{\D\linbarrKKTvecfld[\allvarotriter]{\barrparamotriter}\paren*{\tmeeig\dirNewtonallvarbarrparamotriter} -\D\linbarrKKTvecfld[\allvarotriter]{\barrparamotriter}\paren*{\zerovec[\allvarotriter]}}\sbra*{\dirNewtonallvarbarrparamotriter}\dd{\tmeeig}}\notag\\
            &\leq \int_{0}^{1} \Riemnorm[\allvarotriter]{\paren*{\D\linbarrKKTvecfld[\allvarotriter]{\barrparamotriter}\paren*{\tmeeig\dirNewtonallvarbarrparamotriter} - \D\linbarrKKTvecfld[\allvarotriter]{\barrparamotriter}\paren*{\zerovec[\allvarotriter]}}\sbra*{\dirNewtonallvarbarrparamotriter}}\dd{\tmeeig}\leq \coeffeig[]\constdirNewton^{2}\barrparamotriterm^{2},\label{ineq:linbarrKKTvecflddirNewtonboundbarrparam}
        \end{align}}
        where the first inequality follows from \cref{eq:RiemNewtoneq,eq:linbarrKKTvecfldzerovec,eq:difflinbarrKKTvecfldzerovec}, and the second one from \cref{ineq:derivlinbarKKTvecflddirNewton}.
        Hence, we have that, for any $\otriteridx$-th iteration of \cref{algo:RIPMOuter} satisfying that $\allvarotriter$ is sufficiently close to $\allvaraccum$ and $\barrparamotriterm$ is sufficiently small,
        \begin{align}
            &\frac{1}{\barrparamotriter}\Riemnorm[\ptotriter]{\gradstr[\pt]\Lagfun\paren*{\retr[\pt]\paren*{\dirNewtonptbarrparamotriter},  \ineqLagmultotriter[] + \dirNewtonineqLagmultbarrparamotriter[], \eqLagmultotriter[] + \dirNewtoneqLagmultbarrparamotriter[]}}\\ &\leq\frac{\coeffsev}{\barrparamotriter}\Riemnorm[\ptotriter]{\gradstr[]\paren*{\objfun\circ\retr[\ptotriter]}\paren*{\dirNewtonptbarrparamotriter}- \sum_{\ineqidx\in\ineqset} \paren*{\ineqLagmultotriter[\ineqidx] + \sbra*{\dirNewtonineqLagmultbarrparamotriter[]}_{\ineqidx}}\gradstr[]\paren*{\ineqfun[\ineqidx]\circ\retr[\ptotriter]}\paren*{\dirNewtonptbarrparamotriter}\\
            &\quad+ \sum_{\eqidx\in\eqset} \paren*{\eqLagmultotriter[\eqidx] + \sbra*{\dirNewtoneqLagmultbarrparamotriter[]}_{\eqidx}}\gradstr[]\paren*{\eqfun[\eqidx]\circ\retr[\ptotriter]}\paren*{\dirNewtonptbarrparamotriter}}\leq \frac{\coeffsev}{\barrparamotriter}\Riemnorm[\allvarotriter]{\linbarrKKTvecfld[\allvarotriter]{\barrparamotriter}\paren*{\dirNewtonallvarbarrparamotriter}}\\
            &\leq \frac{\coeffsev\coeffeig[]\constdirNewton^{2}\barrparamotriterm^{2}}{\barrparamotriter},
        \end{align}
        where the first inequality follows from \isextendedversion{\cref{ineq:lincombgradpullbackgrad}}{\cref{lemm:lincombgradineq}} with some $\coeffsev > 1$ and the third one from \cref{ineq:linbarrKKTvecflddirNewtonboundbarrparam}.
        Similarly, 
        \isextendedversion{
        \begin{align}
            &\frac{1}{\barrparamotriter}\norm*{\Ineqfunmat[]\paren*{\retr[\ptotriter]\paren*{\dirNewtonptbarrparamotriter}}\paren*{\ineqLagmultotriter[] + \dirNewtonineqLagmultbarrparamotriter[]}- \barrparamotriter\onevec} \leq \frac{1}{\barrparamotriter}\Riemnorm[\allvarotriter]{\linbarrKKTvecfld[\allvarotriter]{\barrparamotriter}\paren*{\dirNewtonallvarbarrparamotriter}} \leq \frac{\coeffeig[]\constdirNewton^{2}\barrparamotriterm^{2}}{\barrparamotriter}
        \end{align}}
        {$\frac{1}{\barrparamotriter}\norm*{\Ineqfunmat[]\paren*{\retr[\ptotriter]\paren*{\dirNewtonptbarrparamotriter}}\paren*{\ineqLagmultotriter[] + \dirNewtonineqLagmultbarrparamotriter[]}- \barrparamotriter\onevec} \leq \frac{1}{\barrparamotriter}\Riemnorm[\allvarotriter]{\linbarrKKTvecfld[\allvarotriter]{\barrparamotriter}\paren*{\dirNewtonallvarbarrparamotriter}} \leq \frac{\coeffeig[]\constdirNewton^{2}\barrparamotriterm^{2}}{\barrparamotriter}$}
        and
        \isextendedversion{
        \begin{equation}
            \frac{1}{\barrparamotriter}\norm*{\eqfun[]\paren*{\retr[\ptotriter]\paren*{\dirNewtonptbarrparamotriter}}} \leq \frac{1}{\barrparamotriter} \Riemnorm[\allvarotriter]{\linbarrKKTvecfld[\allvarotriter]{\barrparamotriter}\paren*{\dirNewtonallvarbarrparamotriter}} \leq \frac{\coeffeig[]\constdirNewton^{2}\barrparamotriterm^{2}}{\barrparamotriter}
        \end{equation}}
        {$\frac{1}{\barrparamotriter}\norm*{\eqfun[]\paren*{\retr[\ptotriter]\paren*{\dirNewtonptbarrparamotriter}}} \leq \frac{1}{\barrparamotriter} \Riemnorm[\allvarotriter]{\linbarrKKTvecfld[\allvarotriter]{\barrparamotriter}\paren*{\dirNewtonallvarbarrparamotriter}} \leq \frac{\coeffeig[]\constdirNewton^{2}\barrparamotriterm^{2}}{\barrparamotriter}$}
        hold for any $\otriteridx$-th iteration of \cref{algo:RIPMOuter} satisfying that $\allvarotriter$ is sufficiently close to $\allvaraccum$ and $\barrparamotriterm$ is sufficiently small.
        All the right-hand sides of the three inequalities converge to zero as $\barrparamotriterm$ tends to zero under \cref{assu:barrparamdecreaseupodr}, \isextendedversion{which imply that \cref{ineq:gradLagfunnormdirptNewtonbound,ineq:complnormdirptNewtonbound,ineq:eqvionormdirptNewtonbound} hold for any $\otriteridx$-th iteration of \cref{algo:RIPMOuter} satisfying that $\allvarotriter$ is sufficiently close to $\allvaraccum$ and $\barrparamotriterm$ is sufficiently small.
        The proof is complete.}{
        which completes the proof.}
        \qed
    \end{proof}

    \isextendedversion{

\section{Proof of \cref{def:barrparamupdaterule}}\label{appx:proofbarrparamupdaterule}

    \begin{proof}
        Note that $\brc*{\barrparamotriter}_{\otriteridx}$ is monotonically decreasing by definition.
        It follows that 
        \isextendedversion{
        \begin{align}
            \lim_{\otriteridx\to\infty}\frac{\barrparamotriter^{2}}{\barrparamotriterp} = \lim_{\otriteridx\to\infty}\frac{\barrparamotriter^{1-\constnin}}{\coefffiv} = 0 \text{ and } \lim_{\otriteridx\to\infty}\frac{\barrparamotriterp}{\barrparamotriter} = \lim_{\otriteridx\to\infty} \coefffiv \barrparamotriter^{\constnin} = 0,
        \end{align}}
        {$\lim_{\otriteridx\to\infty}\frac{\barrparamotriter^{2}}{\barrparamotriterp} = \lim_{\otriteridx\to\infty}\frac{\barrparamotriter^{1-\constnin}}{\coefffiv} = 0 \text{ and } \lim_{\otriteridx\to\infty}\frac{\barrparamotriterp}{\barrparamotriter} = \lim_{\otriteridx\to\infty} \coefffiv \barrparamotriter^{\constnin} = 0$,}
        which imply \cref{assu:barrparamdecreaseupodr,assu:barrparamdecreaselowodr}, respectively.
        \qed
    \end{proof}

\section{Proof of \cref{lemm:loccoordJacobnonsingaccum}}\label{appx:proofloccoordJacobnonsingaccum}
    \begin{proof}
        For any $\vecfiv=\paren*{\vecfiv[\pt], \vecfiv[{\ineqLagmult[]}], \vecfiv[{\eqLagmult[]}]}\in\tanspc[\pt]\mani\times\setR[\ineqdime]\times\setR[\eqdime]$, we have
        \begin{align}
            \begin{split}\label{eq:derivloccoordbarrKKTvecfld}
                \D\loccoordbarrKKTvecfld[{\barrparam[]}]\paren*{\loccoord{\allvar}}\sbra*{\loccoord{\vecfiv}} = 
                \begin{bmatrix}
                    \D\paren*{\D\chart[\pt]\paren*{\inv{\chart[\pt]}\paren*{\cdot_{\pt}}}\sbra*{\gradstr[\pt]\Lagfun\paren*{\inv{\chart}\paren*{\cdot}}}}\paren*{\loccoord{\allvar}}\sbra*{\loccoord{\vecfiv}}\\
                    \D\paren*{\Ineqfunmat[]\paren*{\inv{\chart[\pt]}\paren*{\cdot_{\pt}}}\inv{\chart[{\ineqLagmult[]}]}\paren*{\cdot_{\ineqLagmult[]}} - \barrparam[]\onevec}\paren*{\loccoord{\allvar}}\sbra*{\loccoord{\vecfiv}}\\
                    \D\paren*{\eqfun[]\paren*{\inv{\chart[\pt]}\paren*{\cdot_{\pt}}}}\paren*{\loccoord{\allvar}}\sbra*{\loccoord{\vecfiv}}
                \end{bmatrix},
            \end{split}
        \end{align}
        where $\inv{\chart[\pt]}\paren*{\cdot_{\pt}}$ and $\inv{\chart[{\ineqLagmult[]}]}\paren*{\cdot_{\ineqLagmult[]}}$ denote the maps $\allvar\mapsto\inv{\chart[\pt]}\paren*{\pt}$ and $\allvar\mapsto\inv{\chart[{\ineqLagmult[]}]}\paren*{\ineqLagmult[]}=\ineqLagmult[]$, respectively. 
        Note that \eqcref{eq:derivloccoordbarrKKTvecfld} is independent of the value of $\barrparam[]$.
        In the following, we analyze each component of \cref{eq:derivloccoordbarrKKTvecfld}.
        For the first component, it follows from the product rule that
        \begin{align}
            \begin{split}\label{eq:derivloccoordbarrKKTvecfldfirstcomp}
                &\D\paren*{\D\chart[\pt]\paren*{\inv{\chart[\pt]}\paren*{\cdot_{\pt}}}\sbra*{\gradstr[\pt]\Lagfun\paren*{\inv{\chart}\paren*{\cdot }}}}\paren*{\loccoord{\allvar}}\sbra*{\loccoord{\vecfiv}}\\
                &= \D\paren*{\D\chart[\pt]\paren*{\inv{\chart[\pt]}\paren*{\cdot_{\pt}}}}\paren*{\loccoord{\allvar}}\sbra*{\loccoord{\vecfiv}}\sbra*{\gradstr[\pt]\Lagfun\paren*{\inv{\chart}\paren*{\loccoord{\allvar}}}}\\
                &\quad + \D\chart[\pt]\paren*{\inv{\chart[\pt]}\paren*{\loccoord{\pt}}}\sbra*{ \D\paren*{\gradstr[\pt]\Lagfun\paren*{\inv{\chart}\paren*{\cdot}}}\paren*{\loccoord{\allvar}}\sbra*{\loccoord{\vecfiv}}}\\
                &= \D\paren*{\D\chart[\pt]\paren*{\inv{\chart[\pt]}\paren*{\cdot_{\pt}}}}\paren*{\loccoord{\allvar}}\sbra*{\loccoord{\vecfiv}}\sbra*{\gradstr[\pt]\Lagfun\paren*{\inv{\chart}\paren*{\loccoord{\allvar}}}} \\
                &\quad + \D\chart[\pt]\paren*{\inv{\chart[\pt]}\paren*{\loccoord{\pt}}}\sbra*{\D\paren*{\gradstr[\pt]\Lagfun\paren*{\inv{\chart}\paren*{\cdot_{\pt}, \loccoord{\ineqLagmult[]}, \loccoord{\eqLagmult[]}}}}\paren*{\loccoord{\pt}}\sbra*{\loccoord{\vecfiv[\pt]}}}\\
                &\quad + \D\chart[\pt]\paren*{\inv{\chart[\pt]}\paren*{\loccoord{\pt}}}\sbra*{\D\paren*{\gradstr[\pt]\Lagfun\paren*{\inv{\chart}\paren*{\loccoord{\pt}, \cdot_{\ineqLagmult[]}, \loccoord{\eqLagmult[]}}}}\paren*{\loccoord{\ineqLagmult[]}}\sbra*{\loccoord{\vecfiv[{\ineqLagmult[]}]}}}\\
                &\quad + \D\chart[\pt]\paren*{\inv{\chart[\pt]}\paren*{\loccoord{\pt}}}\sbra*{\D\paren*{\gradstr[\pt]\Lagfun\paren*{\inv{\chart}\paren*{\loccoord{\pt},  \loccoord{\ineqLagmult[]}, \cdot_{\eqLagmult[]}}}}\paren*{\loccoord{\eqLagmult[]}}\sbra*{\loccoord{\vecfiv[{\eqLagmult[]}]}}},
            \end{split}
        \end{align}
        where $\inv{\chart}\paren*{\cdot_{\pt}, \loccoord{\ineqLagmult[]}, \loccoord{\eqLagmult[]}}$, $\inv{\chart}\paren*{\loccoord{\pt}, \cdot_{\ineqLagmult[]}, \loccoord{\eqLagmult[]}}$, and $\inv{\chart}\paren*{\loccoord{\pt}, \loccoord{\ineqLagmult[]}, \cdot_{\eqLagmult[]}}$ denote the maps $\loccoord{\pt}\mapsto\inv{\chart}\paren*{\paren*{\loccoord{\pt}, \loccoord{\ineqLagmult[]}, \loccoord{\eqLagmult[]}}} = \paren*{\pt, \ineqLagmult[], \eqLagmult[]}$, $\loccoord{\ineqLagmult[]}\mapsto\inv{\chart}\paren*{\paren*{\loccoord{\pt}, \loccoord{\ineqLagmult[]}, \loccoord{\eqLagmult[]}}} = \paren*{\pt, \ineqLagmult[], \eqLagmult[]}$, and $\loccoord{\eqLagmult[]}\mapsto\inv{\chart}\paren*{\paren*{\loccoord{\pt}, \loccoord{\ineqLagmult[]}, \loccoord{\eqLagmult[]}}} = \paren*{\pt, \ineqLagmult[], \eqLagmult[]}$, respectively.
        We have 
        \begin{align}
            \begin{split}\label{eq:loccoordderivgradLagfunineqLagmult}
                &\D\paren*{\gradstr[\pt]\Lagfun\paren*{\inv{\chart}\paren*{\loccoord{\pt}, \cdot_{\ineqLagmult[]}}}}\paren*{\loccoord{\ineqLagmult[]}}\sbra*{\loccoord{\vecfiv[{\ineqLagmult[]}]}}\\
                &= - \sum_{\ineqidx\in\ineqset}\sbra*{\D\inv{\chart[{\ineqLagmult[]}]}\paren*{\loccoord{\ineqLagmult[]}}\sbra*{\loccoord{\vecfiv[{\ineqLagmult[]}]}}}_{\ineqidx}\gradstr\ineqfun[\ineqidx]\paren*{\inv{\chart[\pt]}\paren*{\loccoord{\pt}}} = - \ineqgradopr[{\inv{\chart[\pt]}\paren*{\loccoord{\pt}}}]\sbra*{\loccoord{\vecfiv[{\ineqLagmult[]}]}},
            \end{split}
        \end{align}
        where the second equality follows since $\D\inv{\chart[{\ineqLagmult[]}]}\paren*{\loccoord{\ineqLagmult[]}}$ is the identity map.
        Similarly, it follows that 
        \begin{align}
            \begin{split}\label{eq:loccoordderivgradLagfuneqLagmult}
                &\D\paren*{\gradstr[\pt]\Lagfun\paren*{\inv{\chart}\paren*{\loccoord{\pt},  \loccoord{\ineqLagmult[]}, \cdot_{\eqLagmult[]}}}}\paren*{\loccoord{\eqLagmult[]}}\sbra*{\loccoord{\vecfiv[{\eqLagmult[]}]}}\\
                &= \sum_{\eqidx\in\eqset}\sbra*{\D\inv{\chart[{\eqLagmult[]}]}\paren*{\loccoord{\eqLagmult[]}}\sbra*{\loccoord{\vecfiv[{\eqLagmult[]}]}}}_{\eqidx}\gradstr\eqfun[\eqidx]\paren*{\inv{\chart[\pt]}\paren*{\loccoord{\pt}}} = \eqgradopr[{\inv{\chart[\pt]}\paren*{\loccoord{\pt}}}]\sbra*{\loccoord{\vecfiv[{\eqLagmult[]}]}},
            \end{split}
        \end{align}
        where the second equality follows since $\D\inv{\chart[{\eqLagmult[]}]}\paren*{\loccoord{\eqLagmult[]}}$ is the identity map.
        In addition, from $\gradstr[\pt]\Lagfun\paren*{\allvaraccum}=\zerovec[\ptaccum]$ and \cite[Equation~(5.7)]{Absiletal08OptimAlgoonMatMani}, we have
        \begin{align}
            \begin{split}\label{eq:derivgradHessptaccum}
                &\D\paren*{\gradstr[\pt]\Lagfun\paren*{\inv{\chart}\paren*{\cdot_{\pt}, \loccoord{\ineqLagmultaccum[]}}}}\paren*{\loccoord{\ptaccum}}\sbra*{\loccoord{\vecfiv[\ptaccum]}}\\
                &= \Riemcxt{\D\inv{\chart[\pt]}\paren*{\ptaccum}\sbra*{\loccoord{\vecfiv[\ptaccum]}}}{\gradstr[\pt]\Lagfun\paren*{\inv{\chart}\paren*{\loccoord{\allvaraccum}}}}\\
                &= \Hess[\pt]\Lagfun\paren*{\inv{\chart}\paren*{\loccoord{\allvaraccum}}} \circ \D \inv{\chart[\pt]}\paren*{\loccoord{\ptaccum}}\sbra*{\loccoord{\vecfiv[\pt]}},
            \end{split}
        \end{align}
        where the second equality follows from the definition  of the Riemannian Hessian.
        Combining \cref{eq:derivgradHessptaccum,eq:loccoordderivgradLagfunineqLagmult,eq:loccoordderivgradLagfuneqLagmult} with \cref{eq:derivloccoordbarrKKTvecfldfirstcomp} and $\gradstr[\pt]\Lagfun\paren*{\allvaraccum}=\zerovec[\ptaccum]$, again, yields 
        \begin{align}
            \begin{split}\label{eq:derivloccoordbarrKKTvecfldfirstcompallvaraccum}
                &\D\paren*{\D\chart[\pt]\paren*{\inv{\chart[\pt]}\paren*{\cdot_{\pt}}}\sbra*{\gradstr[\pt]\Lagfun\paren*{\inv{\chart}\paren*{\cdot_{\allvar}}}}}\paren*{\loccoord{\allvaraccum}}\sbra*{\loccoord{\vecfiv}}\\
                &= \D\chart[\pt]\paren*{\inv{\chart[\pt]}\paren*{\loccoord{\ptaccum}}}\circ\Hess[\pt]\Lagfun\paren*{\inv{\chart}\paren*{\loccoord{\allvaraccum}}} \circ \D \inv{\chart[\pt]}\paren*{\loccoord{\ptaccum}}\sbra*{\loccoord{\vecfiv[\pt]}}\\
                &\quad - \D\chart[\pt]\paren*{\inv{\chart[\pt]}\paren*{\loccoord{\ptaccum}}}\sbra*{\ineqgradopr[{\inv{\chart[\pt]}\paren*{\loccoord{\ptaccum}}}]\sbra*{\loccoord{\vecfiv[{\ineqLagmult[]}]}}} + \D\chart[\pt]\paren*{\inv{\chart[\pt]}\paren*{\loccoord{\ptaccum}}}\sbra*{\eqgradopr[{\inv{\chart[\pt]}\paren*{\loccoord{\ptaccum}}}]\sbra*{\loccoord{\vecfiv[{\eqLagmult[]}]}}}.
            \end{split}
        \end{align}
        Next, we consider the second component of \cref{eq:derivloccoordbarrKKTvecfld}: for each $\ineqidx\in\ineqset$, 
        \begin{align}
            \begin{split}
                &\sbra*{\D\paren*{\Ineqfunmat[]\paren*{\inv{\chart[\pt]}\paren*{\cdot_{\pt}}}\inv{\chart[{\ineqLagmult[]}]}\paren*{\cdot_{\ineqLagmult[]}} - \barrparam[]\onevec}\paren*{\loccoord{\allvar}}\sbra*{\loccoord{\vecfiv}}}_{\ineqidx}\\
                &=
                \D\paren*{\ineqfun[\ineqidx]\paren*{\inv{\chart[\pt]}\paren*{\cdot_{\pt}}}\sbra*{\inv{\chart[{\ineqLagmult[]}]}\paren*{\cdot_{\ineqLagmult[]}}}_{\ineqidx} - \barrparam[]\onevec}\paren*{\loccoord{\allvar}}\sbra*{\loccoord{\vecfiv}}\\
                &= \D\ineqfun[\ineqidx]\paren*{\inv{\chart[\pt]}\paren*{\loccoord{\pt}}}\sbra*{\D\inv{\chart[\pt]}\paren*{\loccoord{\pt}}\sbra*{\loccoord{\vecfiv[\pt]}}}\cdot\ineqLagmult[\ineqidx] + \ineqfun[\ineqidx]\paren*{\inv{\chart[\pt]}\paren*{\loccoord{\pt}}} \cdot\sbra*{\D\inv{\chart[{\ineqLagmult[]}]}\paren*{\loccoord{\ineqLagmult[]}}\sbra*{\loccoord{\vecfiv[{\ineqLagmult[]}]}}}_{\ineqidx}\\
                &= \ineqLagmult[\ineqidx]\metr[\inv{\chart[\pt]}\paren*{\loccoord{\pt}}]{\gradstr\ineqfun[\ineqidx]\paren*{\inv{\chart[\pt]}\paren*{\loccoord{\pt}}}}{\D\inv{\chart[\pt]}\paren*{\loccoord{\pt}}\sbra*{\loccoord{\vecfiv[\pt]}}}
                + \ineqfun[\ineqidx]\paren*{\inv{\chart[\pt]}\paren*{\loccoord{\pt}}} \cdot\sbra*{\loccoord{\vecfiv[{\ineqLagmult[]}]}}_{\ineqidx},
            \end{split}
        \end{align}
        which implies
        \begin{align}
            \begin{split}\label{eq:derivloccoordbarrKKTvecfldsecondcomp}
                &\D\paren*{\Ineqfunmat[]\paren*{\inv{\chart[\pt]}\paren*{\cdot_{\pt}}}\inv{\chart[{\ineqLagmult[]}]}\paren*{\cdot_{\ineqLagmult[]}} - \barrparam[]\onevec}\paren*{\loccoord{\allvar}}\sbra*{\loccoord{\vecfiv}}\\
                &= \IneqLagmultmat[]\coineqgradopr[{\inv{\chart[\pt]}\paren*{\loccoord{\pt}}}]\sbra*{\D\inv{\chart[\pt]}\paren*{\loccoord{\pt}}\sbra*{\loccoord{\vecfiv[\pt]}}} + \Ineqfunmat[]\paren*{\inv{\chart[\pt]}\paren*{\loccoord{\pt}}}\loccoord{\vecfiv[{\ineqLagmult[]}]}.
            \end{split}
        \end{align}
        As for the third component of \cref{eq:derivloccoordbarrKKTvecfld}, we have that, for each $\eqidx\in\eqset$,
        \begin{align}
            &\sbra*{\D\paren*{\eqfun[]\paren*{\inv{\chart[\pt]}\paren*{\cdot_{\pt}}}}\paren*{\loccoord{\allvar}}\sbra*{\loccoord{\vecfiv}}}_{\eqidx}\\
            &= \D\paren*{\eqfun[\eqidx]\paren*{\inv{\chart[\pt]}\paren*{\loccoord{\pt}}}}\sbra*{\D\inv{\chart[\pt]}\paren*{\loccoord{\pt}}\sbra*{\loccoord{\vecfiv[\pt]}}} = \metr[\pt]{\gradstr\eqfun[\eqidx]\paren*{\inv{\chart[\pt]}\paren*{\loccoord{\pt}}}}{\loccoord{\vecfiv[\pt]}},
        \end{align}
        where the second equality follows from \cref{eq:riemgrad} and $\chart[\pt]$ being the identity map.
        This implies
        \begin{align}\label{eq:derivloccoordbarrKKTvecfldthirdcomp}
            \D\paren*{\eqfun[]\paren*{\inv{\chart[\pt]}\paren*{\cdot_{\pt}}}}\paren*{\loccoord{\allvar}}\sbra*{\loccoord{\vecfiv}} = \coeqgradopr[{\inv{\chart[\pt]}\paren*{\loccoord{\pt}}}]\sbra*{\D\inv{\chart[\pt]}\paren*{\loccoord{\pt}}\sbra*{\loccoord{\vecfiv[\pt]}}}.
        \end{align}
        Substituting \cref{eq:derivloccoordbarrKKTvecfldfirstcompallvaraccum,eq:derivloccoordbarrKKTvecfldsecondcomp,eq:derivloccoordbarrKKTvecfldthirdcomp} into \cref{eq:derivloccoordbarrKKTvecfld} yields  
        \begin{align}
                \begin{split}\label{eq:derivloccoordbarrKKTvecfldaccum}
                    &\D\loccoordbarrKKTvecfld[{\barrparam[]}]\paren*{\loccoord{\allvaraccum}}\\ 
                    &=\D\chart\paren*{\inv{\chart}\paren*{\loccoord{\allvaraccum}}} \circ 
                    \begin{bmatrix}
                        \Hess[\pt]\Lagfun\paren*{\inv{\chart}\paren*{\loccoord{\allvaraccum}}}
                        &- \ineqgradopr[{\inv{\chart[\pt]}\paren*{\loccoord{\ptaccum}}}] &\eqgradopr[{\inv{\chart[\pt]}\paren*{\loccoord{\ptaccum}}}] \\ 
                        {\IneqLagmultmat[]}^{\accumsymbol} \coineqgradopr[{\inv{\chart[\pt]}\paren*{\loccoord{\ptaccum}}}] & \Ineqfunmat[]\paren*{\inv{\chart[\pt]}\paren*{\loccoord{\ptaccum}}} & 0\\
                        \coeqgradopr[{\inv{\chart[\pt]}\paren*{\loccoord{\ptaccum}}}] & 0 & 0
                    \end{bmatrix} \circ \D\inv{\chart}\paren*{\loccoord{\allvaraccum}}
                \end{split}\\
                &= \D\chart\paren*{\inv{\chart}\paren*{\loccoord{\allvaraccum}}} \circ \Jacobian[\barrKKTvecfld]\paren*{\inv{\chart}\paren*{\loccoord{\allvaraccum}}} \circ \D\inv{\chart}\paren*{\loccoord{\allvaraccum}},
        \end{align}
        where we write ${\IneqLagmultmat[]}^{\accumsymbol}\in\setR[\ineqdime\times\ineqdime]$ for $\diag\paren*{\ineqLagmultaccum[]}$ and
        \begin{align}
            \begin{aligned}
                \D\chart\paren*{\inv{\chart}\paren*{\loccoord{\allvaraccum}}}\colon\tanspc[\ptaccum]\mani\times\setR[\ineqdime]\times\setR[\eqdime]&\to\setR[\dime]\times\setR[\ineqdime]\times\setR[\eqdime]\\
                \paren*{\vecfiv[\ptaccum], \vecfiv[{\ineqLagmultaccum[]}], \vecfiv[{\eqLagmultaccum[]}]}&\mapsto
                \begin{bmatrix}
                    \D\chart[\pt]\paren*{\inv{\chart[\pt]}\paren*{\loccoord{\ptaccum}}}\sbra*{\vecfiv[\ptaccum]}\\
                    \vecfiv[{\ineqLagmultaccum[]}]\\
                    \vecfiv[{\eqLagmultaccum[]}]
                \end{bmatrix},
            \end{aligned}
            \\
            \begin{aligned}
                \D\inv{\chart}\paren*{\loccoord{\allvaraccum}}\colon\setR[\dime]\times\setR[\ineqdime]\times\setR[\eqdime]&\to\tanspc[\ptaccum]\mani\times\setR[\ineqdime]\times\setR[\eqdime]\\
                \paren*{\loccoord{\vecfiv[\ptaccum]}, \loccoord{\vecfiv[{\ineqLagmultaccum[]}]}, \loccoord{\vecfiv[{\eqLagmultaccum[]}]}}&\mapsto
                \begin{bmatrix}
                    \D\chart[\pt]\paren*{\inv{\chart[\pt]}\paren*{\loccoord{\ptaccum}}}\sbra*{\loccoord{\vecfiv[\ptaccum]}}\\ \loccoord{\vecfiv[{\ineqLagmultaccum[]}]}\\
                    \loccoord{\vecfiv[{\eqLagmultaccum[]}]}
                \end{bmatrix}.
            \end{aligned}
        \end{align}
        Note that, since the maps $\D\chart\paren*{\inv{\chart}\paren*{\loccoord{\allvaraccum}}}$ and $\D\inv{\chart}\paren*{\loccoord{\allvaraccum}}$ are bijective and the operator $\Jacobian[\barrKKTvecfld]\paren*{\inv{\chart}\paren*{\loccoord{\allvaraccum}}}$ is nonsingular by \cref{lemm:Jacobnonsingaccum},  $\D\loccoordbarrKKTvecfld[{\barrparam[]}]\paren*{\loccoord{\allvaraccum}}$ is also nonsingular.
        \qed
    \end{proof}

\section{Proof of \cref{lemm:solimplloccoordbarrKKTvecfld}}\label{appx:proofsolimplloccoordbarrKKTvecfld}
    \begin{proof}
        We first consider \cref{lemm:solimplloccoordbarrKKTvecfldexist}.
        Since $\implloccoordbarrKKTvecfld\paren*{\loccoord{\allvar}, \zerovec[], \zerovec[], \zerovec[]} = \loccoordbarrKKTvecfld\paren*{\loccoord{\allvar}, 0}$ holds by definition, it follows that $\implloccoordbarrKKTvecfld\paren*{\loccoord{\allvaraccum}, \zerovec[], \zerovec[], \zerovec[]} = \loccoordbarrKKTvecfld\paren*{\loccoord{\allvaraccum}, \zerovec[]} = \zerovec[]$ and $
        \D\implloccoordbarrKKTvecfld\paren*{\cdot, \zerovec[], \zerovec[], \zerovec[]}\paren*{\loccoord{\allvaraccum}}=\D\loccoordbarrKKTvecfld[{\zerovec[]}]\paren*{\loccoord{\allvaraccum}}$ is nonsingular by \cref{lemm:loccoordJacobnonsingaccum}.
        Thus, by the implicit function theorem~\cite[Theorem~C.40]{Lee12IntrotoSmthManibook2ndedn}, 
        there exist a positive scalar $\implconstone > 0$ and a continuously differentiable function $\loccoord{\allvar}\paren*{\implvecone,\implvectwo, \implvecthr}\colon\nbhdone\paren*{\implconstone}\to\setR[\dime]$ such that $\loccoord{\allvar}\paren*{\zerovec[],\zerovec[], \zerovec[]} = \loccoord{\allvaraccum}$ and $\implloccoordbarrKKTvecfld\paren*{\loccoord{\allvar}\paren*{\implvecone,\implvectwo, \implvecthr}, \implvecone, \implvectwo, \implvecthr} = 0$ for any $\paren*{\implvecone, \implvectwo, \implvecthr}\in\nbhdone\paren*{\implconstone}$.
        Notice that the implicit function theorem also ensures the existence and uniqueness of $\loccoord{\allvar}\paren*{\implvecone, \implvectwo, \implvecthr}$ in $\nbhdone\paren*{\implconstone}$.
        The proof of \cref{lemm:solimplloccoordbarrKKTvecfldexist} is complete.

        Next, we consider \cref{lemm:solimplloccoordbarrKKTvecfldorder}.
        By taking $\implconstone > 0$ smaller if necessary, we have
        \isextendedversion{\begin{align}
            \begin{split}\label{eq:loccoordallvarTaylor}
                &\loccoord{\allvar}\paren*{\implvecone[2],\implvectwo[2], \implvecthr[2]} = \loccoord{\allvar}\paren*{\implvecone[1],\implvectwo[1], \implvecthr[1]}\\
                &\quad + \int_{0}^{1} \D\loccoord{\allvar}\paren*{\paren*{\implvecone[1] + \tmenin\paren*{\implvecone[2] - \implvecone[1]}, \implvectwo[1] + \tmenin\paren*{\implvectwo[2] - \implvectwo[1]}, \implvecthr[1] + \tmenin\paren*{\implvecthr[2] - \implvecthr[1]}}} \sbra*{\paren*{\implvecone[2] - \implvecone[1], \implvectwo[2] - \implvectwo[1], \implvecthr[2] - \implvecthr[1]}}\dd{\tmenin}.
            \end{split}
        \end{align}}{\begin{align}
            \begin{split}\label{eq:loccoordallvarTaylor}
                &\loccoord{\allvar}\paren*{\implvecone[2],\implvectwo[2], \implvecthr[2]} = \loccoord{\allvar}\paren*{\implvecone[1],\implvectwo[1], \implvecthr[1]}\\
                &\quad + \int_{0}^{1} \D\loccoord{\allvar}\paren*{\paren*{\implvecone[1] + \tmenin\paren*{\implvecone[2] - \implvecone[1]}, \implvectwo[1] + \tmenin\paren*{\implvectwo[2] - \implvectwo[1]}, \implvecthr[1] + \tmenin\paren*{\implvecthr[2] - \implvecthr[1]}}}\\
                &\qquad\sbra*{\paren*{\implvecone[2] - \implvecone[1], \implvectwo[2] - \implvectwo[1], \implvecthr[2] - \implvecthr[1]}}\dd{\tmenin}.
            \end{split}
        \end{align}}
        for all $\loccoord{\allvar}\paren*{\implvecone[1],\implvectwo[1], \implvecthr[1]}, \loccoord{\allvar}\paren*{\implvecone[2],\implvectwo[2], \implvecthr[2]} \in \nbhdone\paren*{\implconstone}$.
        We first consider the upper bound on the norm of the difference between $\loccoord{\allvar}\paren*{\implvecone[1],\implvectwo[1], \implvecthr[1]}$ and $\loccoord{\allvar}\paren*{\implvecone[2],\implvectwo[2], \implvecthr[2]}$.
        There exists $\consttht > 0$ such that
        \isextendedversion{\begin{equation}
            \begin{aligned}[t]\label{ineq:diffloccoordallvarupperbound}
                &\norm*{\loccoord{\allvar}\paren*{\implvecone[2],\implvectwo[2], \implvecthr[2]} - \loccoord{\allvar}\paren*{\implvecone[1],\implvectwo[1], \implvecthr[1]}}\\
                &\leq \int_{0}^{1} \opnorm{\D\loccoord{\allvar}\paren*{\paren*{\implvecone[1] + \tmenin\paren*{\implvecone[2] - \implvecone[1]}, \implvectwo[1] + \tmenin\paren*{\implvectwo[2] - \implvectwo[1]}, \implvecthr[1] + \tmenin\paren*{\implvecthr[2] - \implvecthr[1]}}}} \norm*{\paren*{\implvecone[2] - \implvecone[1], \implvectwo[2] - \implvectwo[1], \implvecthr[2] - \implvecthr[1]}} \dd{\tmenin}\\
                &\leq \consttht \norm*{\paren*{\implvecone[2] - \implvecone[1], \implvectwo[2] - \implvectwo[1], \implvecthr[2] - \implvecthr[1]}}\\
                &\leq \consttht \paren*{\norm*{\implvecone[2] - \implvecone[1]} + \norm*{\implvectwo[2] - \implvectwo[1]} + \norm*{\implvecthr[2] - \implvecthr[1]}},
            \end{aligned}
        \end{equation}}{
        \begin{align}
            \begin{split}\label{ineq:diffloccoordallvarupperbound}
                &\norm*{\loccoord{\allvar}\paren*{\implvecone[2],\implvectwo[2], \implvecthr[2]} - \loccoord{\allvar}\paren*{\implvecone[1],\implvectwo[1], \implvecthr[1]}}\\
                &\leq \int_{0}^{1} \opnorm{\D\loccoord{\allvar}\paren*{\paren*{\implvecone[1] + \tmenin\paren*{\implvecone[2] - \implvecone[1]}, \implvectwo[1] + \tmenin\paren*{\implvectwo[2] - \implvectwo[1]}, \implvecthr[1] + \tmenin\paren*{\implvecthr[2] - \implvecthr[1]}}}}\\
                &\quad\norm*{\paren*{\implvecone[2] - \implvecone[1], \implvectwo[2] - \implvectwo[1], \implvecthr[2] - \implvecthr[1]}} \dd{\tmenin} \leq \consttht \norm*{\paren*{\implvecone[2] - \implvecone[1], \implvectwo[2] - \implvectwo[1], \implvecthr[2] - \implvecthr[1]}}\\
                &\leq \consttht \paren*{\norm*{\implvecone[2] - \implvecone[1]} + \norm*{\implvectwo[2] - \implvectwo[1]} + \norm*{\implvecthr[2] - \implvecthr[1]}},
            \end{split}
        \end{align}}
        where the second inequality follows from the continuity of $\D\loccoord{\allvar}\paren*{\cdot, \cdot, \cdot}$.
        We next derive the lower bound on $\norm*{\loccoord{\allvar}\paren*{\implvecone[2],\implvectwo[2], \implvecthr[2]} - \loccoord{\allvar}\paren*{\implvecone[1],\implvectwo[1], \implvecthr[1]}}$.
        To this end, we derive auxiliary bounds as follows: by \cref{def:implloccoordbarrKKTvecfld} and the definition of $\loccoord{\allvar}\paren*{\implvecone,\implvectwo, \implvecthr}$, it follows that \isextendedversion{
        \begin{equation}
            \loccoordbarrKKTvecfld[0]\paren*{\loccoord{\allvar}\paren*{\implvecone,\implvectwo, \implvecthr}} = 
            \begin{bmatrix}
                \implvecone\\
                \implvectwo\\
                \implvecthr
            \end{bmatrix}
        \end{equation}}{$\loccoordbarrKKTvecfld[0]\paren*{\loccoord{\allvar}\paren*{\implvecone,\implvectwo, \implvecthr}} = 
            \begin{bmatrix}
                \implvecone\\
                \implvectwo\\
                \implvecthr
            \end{bmatrix}$.}
        Differentiating it yields 
        \begin{align}
            \D\paren*{\loccoordbarrKKTvecfld[0]\circ{\loccoord{\allvar}}}\paren*{\implvecone,\implvectwo, \implvecthr} = \id,
        \end{align}
        which implies
        \begin{align}
            \D\loccoordbarrKKTvecfld[0]\paren*{\loccoord{\allvar}\paren*{\implvecone,\implvectwo,\implvecthr}} \circ \D\loccoord{\allvar}\paren*{\implvecone, \implvectwo,\implvecthr}
            = \id.
        \end{align}
        Since $\D\loccoordbarrKKTvecfld[{0}]\paren*{\cdot}$ is nonsingular at $\loccoord{\allvaraccum}$ and is continuous by definition, $\D\loccoordbarrKKTvecfld[{0}]\paren*{\cdot}$ remains nonsingular around $\loccoord{\allvaraccum}$ and hence
        \begin{align}\label{eq:invderivloccoordbarrKKTvecfld}
            \D\loccoord{\allvar}\paren*{\implvecone, \implvectwo, \implvecthr} = \inv{\D\loccoordbarrKKTvecfld[0]}\paren*{\loccoord{\allvar}\paren*{\implvecone,\implvectwo,\implvecthr}}
        \end{align}
        for any $\paren*{\implvecone, \implvectwo,\implvecthr}\in\nbhdone\paren*{\implconstone}$ by taking $\implconstone > 0$ smaller if necessary.
        Therefore, there exists $\constsvt > 0$ such that 
        \begin{align}
            \begin{split}\label{ineq:diffimplvecsbound}
                &\norm*{\paren*{\implvecone[2] - \implvecone[1], \implvectwo[2] - \implvectwo[1], \implvecthr[2] - \implvecthr[1]}}\\
                &= \norm*{\D\loccoordbarrKKTvecfld[0]\paren*{\loccoord{\allvar}\paren*{\implvecone[1],\implvectwo[1],\implvecthr[1]}} \circ \inv{\D\loccoordbarrKKTvecfld[0]}\paren*{\loccoord{\allvar}\paren*{\implvecone[1], \implvectwo[1], \implvecthr[1]}}\sbra*{\paren*{\implvecone[2] - \implvecone[1], \implvectwo[2] - \implvectwo[1], \implvecthr[2] - \implvecthr[1]}}}\\
                &\leq \opnorm{\D\loccoordbarrKKTvecfld[0]\paren*{\loccoord{\allvar}\paren*{\implvecone[1],\implvectwo[1],\implvecthr[1]}}} \norm*{\inv{\D\loccoordbarrKKTvecfld[0]}\paren*{\loccoord{\allvar}\paren*{\implvecone[1],\implvectwo[1],\implvecthr[1]}}\sbra*{\paren*{\implvecone[2] - \implvecone[1], \implvectwo[2] - \implvectwo[1], \implvecthr[2] - \implvecthr[1]}}}\\
                &\leq \frac{1}{\paren*{1 + \sqrt{3}} \constsvt} 
                \norm*{\D\loccoord{\allvar}\paren*{\implvecone[1], \implvectwo[1], \implvecthr[1]}\sbra*{\paren*{\implvecone[2] - \implvecone[1], \implvectwo[2] - \implvectwo[1], \implvecthr[2] - \implvecthr[1]}}},
            \end{split}
        \end{align}
        where 
        the second inequality follows from \cref{eq:invderivloccoordbarrKKTvecfld} and the boundedness of $\D\loccoordbarrKKTvecfld[0]\paren*{\cdot}$ around $\loccoord{\allvaraccum}$.
        By taking $\implconstone > 0$ smaller again if necessary, we also have
        \begin{align}\label{ineq:diffderivloccoordallvarcont}
            \opnorm{\D\paren*{\loccoord{\allvar}\paren*{\paren*{\implvecone[1] + \tmenin\paren*{\implvecone[2] - \implvecone[1]}, \implvectwo[1] + \tmenin\paren*{\implvectwo[2] - \implvectwo[1]}, \implvecthr[1] + \tmenin\paren*{\implvecthr[2] - \implvecthr[1]}}} - \loccoord{\allvar}\paren*{\implvecone[1], \implvectwo[1], \implvecthr[1]}}} \leq \constsvt
        \end{align}
        for any $\paren*{\implvecone, \implvectwo, \implvecthr}\in\nbhdone\paren*{\implconstone}$ by the continuous differentiability of $\loccoord{\allvar}\paren*{\cdot,\cdot, \cdot}$.
        Now, we derive the lower bound on $\norm*{\loccoord{\allvar}\paren*{\implvecone[2],\implvectwo[2], \implvecthr[2]} - \loccoord{\allvar}\paren*{\implvecone[1],\implvectwo[1], \implvecthr[1]}}$:
        \isextendedversion{\begin{equation}
            \begin{aligned}[t]\label{ineq:diffloccoordallvarlowerbound}
                &\norm*{\loccoord{\allvar}\paren*{\implvecone[2],\implvectwo[2],\implvecthr[2]} - \loccoord{\allvar}\paren*{\implvecone[1],\implvectwo[1]}}\\
                &\geq \norm*{\D\loccoord{\allvar}\paren*{\implvecone[1], \implvectwo[1], \implvecthr[1]}\sbra*{\paren*{\implvecone[2] - \implvecone[1], \implvectwo[2] - \implvectwo[1], \implvecthr[2] - \implvecthr[1]}}}\\
                &\quad - \norm*{\loccoord{\allvar}\paren*{\implvecone[2],\implvectwo[2],\implvecthr[2]} - \loccoord{\allvar}\paren*{\implvecone[1],\implvectwo[1],\implvecthr[1]} - \D\loccoord{\allvar}\paren*{\implvecone[1], \implvectwo[1],\implvecthr[1]}\sbra*{\paren*{\implvecone[2] - \implvecone[1], \implvectwo[2] - \implvectwo[1], \implvecthr[2] - \implvecthr[1]}}}\\
                &\geq \paren*{1 + \sqrt{3}} \constsvt \norm*{\paren*{\implvecone[2] - \implvecone[1], \implvectwo[2] - \implvectwo[1], \implvecthr[2] - \implvecthr[1]}} - \int_{0}^{1} \opnorm{\D\paren*{\loccoord{\allvar}\paren*{\paren*{\implvecone[1] + \tmenin\paren*{\implvecone[2] - \implvecone[1]},\implvectwo[1] + \tmenin\paren*{\implvectwo[2] - \implvectwo[1]}, \\
                &\quad \implvecthr[1] + \tmenin\paren*{\implvecthr[2] - \implvecthr[1]}}} - \loccoord{\allvar}\paren*{\implvecone[1], \implvectwo[1], \implvecthr[1]}}}\norm*{\paren*{\implvecone[2] - \implvecone[1], \implvectwo[2] - \implvectwo[1], \implvecthr[2] - \implvecthr[1]}} \dd{\tmenin}\\
                &\geq \constsvt\sqrt{3}\norm*{\paren*{\implvecone[2] - \implvecone[1], \implvectwo[2] - \implvectwo[1], \implvecthr[2] - \implvecthr[1]}}\\
                &\geq \constsvt \paren*{\norm*{\implvecone[2] - \implvecone[1]} + \norm*{\implvectwo[2] - \implvectwo[1]} + \norm*{\implvecthr[2] - \implvecthr[1]}}
            \end{aligned}
        \end{equation}}{
        \begin{align}
            \begin{split}\label{ineq:diffloccoordallvarlowerbound}
                &\norm*{\loccoord{\allvar}\paren*{\implvecone[2],\implvectwo[2],\implvecthr[2]} - \loccoord{\allvar}\paren*{\implvecone[1],\implvectwo[1]}} \geq \norm*{\D\loccoord{\allvar}\paren*{\implvecone[1], \implvectwo[1], \implvecthr[1]}\sbra*{\paren*{\implvecone[2] - \implvecone[1], \implvectwo[2] - \implvectwo[1], \implvecthr[2] - \implvecthr[1]}}}\\
                &- \norm*{\loccoord{\allvar}\paren*{\implvecone[2],\implvectwo[2],\implvecthr[2]} - \loccoord{\allvar}\paren*{\implvecone[1],\implvectwo[1],\implvecthr[1]} - \D\loccoord{\allvar}\paren*{\implvecone[1], \implvectwo[1],\implvecthr[1]}\sbra*{\paren*{\implvecone[2] - \implvecone[1], \implvectwo[2] - \implvectwo[1], \implvecthr[2] - \implvecthr[1]}}}\\
                &\geq \paren*{1 + \sqrt{3}} \constsvt \norm*{\paren*{\implvecone[2] - \implvecone[1], \implvectwo[2] - \implvectwo[1], \implvecthr[2] - \implvecthr[1]}} \\
                &- \int_{0}^{1} \opnorm{\D\paren*{\loccoord{\allvar}\paren*{\paren*{\implvecone[1] + \tmenin\paren*{\implvecone[2] - \implvecone[1]},\implvectwo[1] + \tmenin\paren*{\implvectwo[2] - \implvectwo[1]}, \implvecthr[1] + \tmenin\paren*{\implvecthr[2] - \implvecthr[1]}}} - \loccoord{\allvar}\paren*{\implvecone[1], \implvectwo[1], \implvecthr[1]}}}\\
                &\norm*{\paren*{\implvecone[2] - \implvecone[1], \implvectwo[2] - \implvectwo[1], \implvecthr[2] - \implvecthr[1]}} \dd{\tmenin} \geq \constsvt\sqrt{3}\norm*{\paren*{\implvecone[2] - \implvecone[1], \implvectwo[2] - \implvectwo[1], \implvecthr[2] - \implvecthr[1]}}\\
                &\geq \constsvt \paren*{\norm*{\implvecone[2] - \implvecone[1]} + \norm*{\implvectwo[2] - \implvectwo[1]} + \norm*{\implvecthr[2] - \implvecthr[1]}}
            \end{split}
        \end{align}}
        for any $\paren*{\implvecone, \implvectwo, \implvecthr}\in\nbhdone\paren*{\implconstone}$, where the second inequality follows from \cref{ineq:diffimplvecsbound,eq:loccoordallvarTaylor}, the third one from \cref{ineq:diffderivloccoordallvarcont}, and the fourth one from the Cauchy-Schwarz inequality with \isextendedversion{
        \begin{equation}
            \trsp{\paren*{\norm*{\implvecone[2] - \implvecone[1]}, \norm*{\implvectwo[2] - \implvectwo[1]}, \norm*{\implvecthr[2] - \implvecthr[1]}}} \text{ and } \trsp{\paren*{1,1,1}}
        \end{equation}
        }{$\trsp{\paren*{\norm*{\implvecone[2] - \implvecone[1]}, \norm*{\implvectwo[2] - \implvectwo[1]}, \norm*{\implvecthr[2] - \implvecthr[1]}}}$ and $\trsp{\paren*{1,1,1}}$} on $\setR[3]$.
        From \cref{ineq:diffloccoordallvarupperbound,ineq:diffloccoordallvarlowerbound}, \eqcref{eq:loccoordbigtheta} holds.
        The proof of \cref{lemm:solimplloccoordbarrKKTvecfldorder} is complete.
        \qed
    \end{proof}

    }{}

\section{Proof of \cref{lemm:diffallvarotriterpaccumbarrparambigO}}\label{appx:proofdiffallvarotriterpaccumbarrparambigO}
    \begin{proof}
        Let $\otriteridx\in\brc*{-1,0}\cup\setN$ be an index where $\allvarotriterp\in\strictfeasirgn\times\setRp[\ineqdime]\times\setR[\eqdime]$ is sufficiently close to $\allvaraccum$ and $\barrparamotriter > 0$ is sufficiently small.         
        We have 
        \begin{align}\label{ineq:diffloccoordallvarotriterpaccum}
            \norm*{\loccoord{\allvarotriterp} - \loccoord{\allvaraccum}} \leq \norm*{\loccoord{\allvarotriterp} - \loccoord{\allvar}\paren*{\zerovec[], \barrparamotriter\onevec, \zerovec[]}} + \norm*{\loccoord{\allvar}\paren*{\zerovec[], \barrparamotriter\onevec, \zerovec[]} - \loccoord{\allvaraccum}}.
        \end{align}
        In the following, we derive the bound on each term of the right-hand side of \cref{ineq:diffloccoordallvarotriterpaccum}.
        For the first term, it holds by \enumicref{lemm:solimplloccoordbarrKKTvecfld}{lemm:solimplloccoordbarrKKTvecfldexist} that the point
        \begin{align}
            \loccoord{\allvar}\paren*{\D\chart[\pt]\paren*{\ptotriterp}\sbra*{\gradstr[\pt]\Lagfun\paren*{\allvarotriterp}}, \Ineqfunmat[]\paren*{\ptotriterp}\ineqLagmultotriterp[], \eqfun[]\paren*{\ptotriterp}}
        \end{align}
        is well-defined and identical to $\loccoord{\allvarotriterp}$
        by the uniqueness of $\loccoord{\allvar}\paren*{\cdot, \cdot, \cdot}$.
        Therefore, there exist $\constfot[1], \constfot[2] \in \setRpp[]$ such that, for any $\otriteridx\in\setNz$ where $\allvarotriterp$ is sufficiently close to $\allvaraccum$ and $\barrparamotriter > 0$ is sufficiently small,
        \isextendedversion{
        \begin{equation}
            \begin{aligned}[t]\label{ineq:diffloccoordallvarotriterpaccumfirst}
                &\norm*{\loccoord{\allvarotriterp} - \loccoord{\allvar}\paren*{\zerovec[], \barrparamotriter\onevec, \zerovec[]}}\\
                &= \norm*{\loccoord{\allvar}\paren*{\D\chart[\pt]\paren*{\ptotriterp}\sbra*{\gradstr[\pt]\Lagfun\paren*{\allvarotriterp}}, \Ineqfunmat[]\paren*{\ptotriterp}\ineqLagmultotriterp[], \eqfun[]\paren*{\ptotriterp}} - \loccoord{\allvar}\paren*{\zerovec[], \barrparamotriter\onevec,\zerovec[]}}\\
                & \leq \constfot[1] \paren*{\norm*{\D\chart[\pt]\paren*{\ptotriterp}\sbra*{\gradstr[\ptotriterp]\Lagfun\paren*{\allvarotriterp}}} + \norm*{\Ineqfunmat[]\paren*{\ptotriterp}\ineqLagmultotriterp[]- \barrparamotriter\onevec} + \norm*{\eqfun[]\paren*{\ptotriterp}}}\\
                & \leq \constfot[1] \paren*{\opnorm{\D\chart[\pt]\paren*{\ptotriterp}}\Riemnorm[\pt]{\gradstr[\pt]\Lagfun\paren*{\allvarotriterp}} + \norm*{\Ineqfunmat[]\paren*{\ptotriterp}\ineqLagmultotriterp[]- \barrparamotriter\onevec} + \norm*{\eqfun[]\paren*{\ptotriterp}}}\\
                &\leq \constfot[1] \paren*{\opnorm{\D\chart[\pt]\paren*{\ptotriterp}} \forcingfungradLag\paren*{\barrparamotriter} + \forcingfuncompl\paren*{\barrparamotriter} + \forcingfuneqvio\paren*{\barrparamotriter}} \leq \constfot[2] \barrparamotriter,
            \end{aligned}
        \end{equation}}
        {\begin{align}
            \begin{split}\label{ineq:diffloccoordallvarotriterpaccumfirst}
                &\norm*{\loccoord{\allvarotriterp} - \loccoord{\allvar}\paren*{\zerovec[], \barrparamotriter\onevec, \zerovec[]}} = \norm*{\loccoord{\allvar}\paren*{\D\chart[\pt]\paren*{\ptotriterp}\sbra*{\gradstr[\pt]\Lagfun\paren*{\allvarotriterp}}, \Ineqfunmat[]\paren*{\ptotriterp}\ineqLagmultotriterp[], \eqfun[]\paren*{\ptotriterp}}\\
                &- \loccoord{\allvar}\paren*{\zerovec[], \barrparamotriter\onevec,\zerovec[]}} \leq \constfot[1] \paren*{\norm*{\D\chart[\pt]\paren*{\ptotriterp}\sbra*{\gradstr[\ptotriterp]\Lagfun\paren*{\allvarotriterp}}} + \norm*{\Ineqfunmat[]\paren*{\ptotriterp}\ineqLagmultotriterp[]- \barrparamotriter\onevec}\\
                &+ \norm*{\eqfun[]\paren*{\ptotriterp}}} \leq \constfot[1] \paren*{\opnorm{\D\chart[\pt]\paren*{\ptotriterp}}\Riemnorm[\pt]{\gradstr[\pt]\Lagfun\paren*{\allvarotriterp}} + \norm*{\Ineqfunmat[]\paren*{\ptotriterp}\ineqLagmultotriterp[]- \barrparamotriter\onevec} \\
                &+ \norm*{\eqfun[]\paren*{\ptotriterp}}} \leq \constfot[1] \paren*{\opnorm{\D\chart[\pt]\paren*{\ptotriterp}} \forcingfungradLag\paren*{\barrparamotriter} + \forcingfuncompl\paren*{\barrparamotriter} + \forcingfuneqvio\paren*{\barrparamotriter}} \leq \constfot[2] \barrparamotriter,
            \end{split}
        \end{align}}
        where the first inequality follows from \enumicref{lemm:solimplloccoordbarrKKTvecfld}{lemm:solimplloccoordbarrKKTvecfldorder}, the third one from \cref{eq:stopcondKKT,eq:stopcondbarrcompl,eq:stopcondeqvio}, and the last one from the boundedness of $\D\chart[\pt]\paren*{\cdot}$ around $\ptaccum$ and \isextendedversion{\cref{ineq:forcingfuncomplbounded,ineq:forcingfungradLagbounded,ineq:forcingfuneqviobounded}}{\cref{assu:forcingfuncomplgradLagbounded}}.
        Moreover, since $\loccoord{\allvar}\paren*{\zerovec[], \zerovec[], \zerovec[]} = \loccoord{\allvaraccum}$ holds by definition, the second term in the right-hand side of \cref{ineq:diffloccoordallvarotriterpaccum} can be bounded as 
        \begin{align}\label{ineq:diffloccoordallvarotriterpaccumsecond}
            \norm*{\loccoord{\allvar}\paren*{\zerovec[], \barrparamotriter\onevec, \zerovec[]} - \loccoord{\allvaraccum}} = \norm*{\loccoord{\allvar}\paren*{\zerovec[], \barrparamotriter\onevec, \zerovec[]} - \loccoord{\allvar}\paren*{\zerovec[], \zerovec[], \zerovec[]}} 
            \leq \constfot[3] \barrparamotriter
        \end{align}
        for some $ \constfot[3] > 0$, where the inequality follows from \enumicref{lemm:solimplloccoordbarrKKTvecfld}{lemm:solimplloccoordbarrKKTvecfldorder} again.
        Therefore, combining \cref{ineq:diffloccoordallvarotriterpaccum} with \cref{ineq:diffloccoordallvarotriterpaccumfirst,ineq:diffloccoordallvarotriterpaccumsecond} yields 
        \begin{align}\label{ineq:diffloccoordallvarotriterpaccumbarrpramotriter}
             \norm*{\loccoord{\allvarotriterp} - \loccoord{\allvaraccum}} \leq \paren*{\constfot[2] + \constfot[3]}\barrparamotriter.
        \end{align}
        
        In the following, we regard the manipulation $\ineqLagmult[] + \dirineqLagmult[]$ and $\eqLagmult[] + \direqLagmult[]$ as retractions on $\setR[\ineqdime]$ and $\setR[\eqdime]$, respectively; that is, we define $\retr[{\ineqLagmult[]}]\paren*{\dirineqLagmult[]}\coloneqq\ineqLagmult[] + \dirineqLagmult[]$ and $\retr[{\eqLagmult[]}]\paren*{\direqLagmult[]}\coloneqq\eqLagmult[] + \direqLagmult[]$.
        We also define the retraction on the product manifold $\mani\times\setR[\ineqdime]\times\setR[\eqdime]$ as $\retr[\allvar]\paren*{\dirallvar[]}=\paren*{\retr[\pt]\paren*{\dirpt}, \retr[{\ineqLagmult[]}]\paren*{\dirineqLagmult[]}, \retr[{\eqLagmult[]}]\paren*{\direqLagmult[]}}$.
        We have
        \isextendedversion{
        \begin{align}
            \begin{split}\label{ineq:distallvarotriterpaccum}
                &\Riemdist{\allvar[\otriteridx + 2]}{\allvaraccum} \leq \Riemdist{\allvar[\otriteridx + 2]}{\allvarotriterp} + \Riemdist{\allvarotriterp}{\allvaraccum}\\
                &= \Riemdist{\retr[\allvarotriterp]\paren*{\dirNewtonallvarbarrparamotriterp}}{\retr[\allvarotriterp]\paren*{\zerovec[\allvarotriterp]}} + \Riemdist{\allvarotriterp}{\allvaraccum}\\
                &\leq \constsix[2] \Riemnorm[{\allvarotriterp}]{\dirNewtonallvarbarrparamotriter} + \constsix[3] \norm*{\loccoord{\allvarotriterp} - \loccoord{\allvaraccum}} \leq \paren*{\constsix[2]\constdirNewton + \constsix[3]\paren*{\constfot[2] + \constfot[3]}}\barrparamotriter
            \end{split}
        \end{align}}
        {\begin{align}
            \begin{split}\label{ineq:distallvarotriterpaccum}
                &\Riemdist{\allvar[\otriteridx + 2]}{\allvaraccum} \leq \Riemdist{\allvar[\otriteridx + 2]}{\allvarotriterp} + \Riemdist{\allvarotriterp}{\allvaraccum}\\
                &= \Riemdist{\retr[\allvarotriterp]\paren*{\dirNewtonallvarbarrparamotriterp}}{\retr[\allvarotriterp]\paren*{\zerovec[\allvarotriterp]}} + \Riemdist{\allvarotriterp}{\allvaraccum}\\ 
                &\leq \constsix[2] \Riemnorm[{\allvarotriterp}]{\dirNewtonallvarbarrparamotriter} + \constsix[3] \norm*{\loccoord{\allvarotriterp} - \loccoord{\allvaraccum}} \leq \paren*{\constsix[2]\constdirNewton + \constsix[3]\paren*{\constfot[2] + \constfot[3]}}\barrparamotriter
            \end{split}
        \end{align}}
        for some $\constsix[2], \constsix[3] \in\setRp[]$, where the equality follows from \isextendedversion{\cref{coro:Newtonstopconds,eq:retrzero}}{\cref{coro:Newtonstopconds} and the definition of the retraction}, the second inequality from \isextendedversion{\cref{lemm:retrtanvecequiv}}{\cite[Lemma~2]{HuangAbsilGallivan2015RiemSymRankOneTRM}} and \isextendedversion{\cref{lemm:EucliRiemdistequiv}}{\cite[Lemma~14.1]{GallivanQiAbsil2012RiemDennisMoreCond}}, and the third one from \cref{lemm:dirNewtonallvarbarrparambound} and \cref{ineq:diffloccoordallvarotriterpaccumbarrpramotriter}.
        Since, without loss of generality, the right-hand side can be made arbitrarily small, \eqcref{ineq:distallvarotriterpaccum} implies that the point $\allvar[\otriteridx + 2]$ remains in a sufficiently small neighborhood of $\allvaraccum$.
        Therefore, using \cref{assu:usedirexactsolution} and the argument above inductively, we conclude that \isextendedversion{\eqcref{ineq:diffloccoordallvarotriterpaccumbigO}}{\cref{lemm:diffallvarotriterpaccumbarrparambigO}} holds.
        \qed
    \end{proof}

\section{Proof of \cref{lemm:diffallvarotriterpaccumbarrparambigtheta}}\label{appx:proofdiffallvarotriterpaccumbarrparambigtheta}
    \begin{proof}
        Recall that, from \cref{lemm:diffallvarotriterpaccumbarrparambigO}, the point $\allvarotriterp$ can be made arbitrarily close to $\allvaraccum$ by taking $\barrparamotriter > 0$ sufficiently small.
        Let $\loccoordtildestd\coloneqq\sbra*{\trsp{0}, \trsp{\onevec[\ineqdime]}, \trsp{0}}\in\setR[\dime]\times\setR[\ineqdime]\times\setR[\eqdime]$.
        Then, we have
        \begin{align}
            \begin{split}
                &\loccoordbarrKKTvecfld\paren*{\loccoord{\allvarotriterp}, \barrparamotriter} = \loccoordbarrKKTvecfld\paren*{\loccoord{\allvaraccum}, \zerovec[]} + \D\loccoordbarrKKTvecfld\paren*{\cdot, \barrparamotriter}\paren*{\loccoord{\allvaraccum}}\sbra*{\loccoord{\allvarotriterp} - \loccoord{\allvaraccum}} + \D\loccoordbarrKKTvecfld\paren*{\loccoord{\allvaraccum}, \cdot}\paren*{0}\sbra*{\barrparamotriter - 0} \\
                &\quad + \int_{0}^{1} \D[]\paren*{\loccoordbarrKKTvecfld\paren*{\paren*{\loccoord{\allvaraccum}, 0} + \tmeten\paren*{\loccoord{\allvarotriterp}, \barrparamotriter}} - \loccoordbarrKKTvecfld\paren*{\loccoord{\allvaraccum}, 0}}\sbra*{\paren*{\loccoord{\allvarotriterp} - \loccoord{\allvaraccum}, \barrparamotriter- 0}}\dd{\tmeten}\\
                &=\D\loccoordbarrKKTvecfld\paren*{\cdot, \barrparamotriter}\paren*{\loccoord{\allvaraccum}}\sbra*{\loccoord{\allvarotriterp} - \loccoord{\allvaraccum}} -\barrparamotriter\loccoordtildestd
                + \residone,
            \end{split}
        \end{align}
        where the second equality follows from $\loccoordbarrKKTvecfld\paren*{\loccoord{\allvaraccum}, \zerovec[]} = \zerovec[]$ and
        \isextendedversion{
        \begin{align}\label{def:residone}
            \residone \coloneqq \int_{0}^{1} \D[]\paren*{\loccoordbarrKKTvecfld\paren*{\paren*{\loccoord{\allvaraccum}, 0} + \tmeten\paren*{\loccoord{\allvarotriterp}, \barrparamotriter}} - \loccoordbarrKKTvecfld\paren*{\loccoord{\allvaraccum}, 0}}\sbra*{\paren*{\loccoord{\allvarotriterp} - \loccoord{\allvaraccum}, \barrparamotriter- 0}}\dd{\tmeten}.
        \end{align}}{$\residone \coloneqq \int_{0}^{1} \D[]\paren*{\loccoordbarrKKTvecfld\paren*{\paren*{\loccoord{\allvaraccum}, 0} + \tmeten\paren*{\loccoord{\allvarotriterp}, \barrparamotriter}} - \loccoordbarrKKTvecfld\paren*{\loccoord{\allvaraccum}, 0}}\sbra*{\paren*{\loccoord{\allvarotriterp} - \loccoord{\allvaraccum}, \barrparamotriter- 0}}\dd{\tmeten}$.}
        This implies 
        \isextendedversion{\begin{align}\label{eq:loccoordbarrKKTvecfldTaylor}
            &\loccoord{\allvarotriterp} - \loccoord{\allvaraccum}\\
            &= \barrparamotriter\inv{\D\loccoordbarrKKTvecfld\paren*{\cdot, \barrparamotriter}}\paren*{\loccoord{\allvaraccum}}\sbra*{\loccoordtildestd}
            + \inv{\D\loccoordbarrKKTvecfld\paren*{\cdot, \barrparamotriter}}\paren*{\loccoord{\allvaraccum}}\sbra*{\loccoordbarrKKTvecfld\paren*{\loccoord{\allvarotriterp}, \barrparamotriter}} - \inv{\D\loccoordbarrKKTvecfld\paren*{\cdot, \barrparamotriter}}\paren*{\loccoord{\allvaraccum}}\sbra*{\residone}
        \end{align}}{\begin{align}\label{eq:loccoordbarrKKTvecfldTaylor}
            \loccoord{\allvarotriterp} - \loccoord{\allvaraccum} = \barrparamotriter\inv{\D\loccoordbarrKKTvecfld\paren*{\cdot, \barrparamotriter}}\paren*{\loccoord{\allvaraccum}}\sbra*{\loccoordtildestd}
            + \inv{\D\loccoordbarrKKTvecfld\paren*{\cdot, \barrparamotriter}}\paren*{\loccoord{\allvaraccum}}\sbra*{\loccoordbarrKKTvecfld\paren*{\loccoord{\allvarotriterp}, \barrparamotriter}} - \inv{\D\loccoordbarrKKTvecfld\paren*{\cdot, \barrparamotriter}}\paren*{\loccoord{\allvaraccum}}\sbra*{\residone}
        \end{align}}
        since $\D\loccoordbarrKKTvecfld\paren*{\cdot, \barrparamotriter}\paren*{\loccoord{\allvaraccum}}$ is nonsingular by \cref{lemm:loccoordJacobnonsingaccum}.
        We derive bounds on terms of \cref{eq:loccoordbarrKKTvecfldTaylor} as follows:
        Let 
        \isextendedversion{\begin{align}\label{def:allvarnonzeroconst}
            \allvarnonzeroconst\paren*{\zerovec[]}\coloneqq\inv{\D\loccoordbarrKKTvecfld\paren*{\cdot, \barrparamotriter}}\paren*{\loccoord{\allvaraccum}}
            \loccoordtildestd
            \in\setR[\dime+\ineqdime+\eqdime]\backslash\brc*{0}.
        \end{align}}{$\allvarnonzeroconst\paren*{\zerovec[]}\coloneqq\inv{\D\loccoordbarrKKTvecfld\paren*{\cdot, \barrparamotriter}}\paren*{\loccoord{\allvaraccum}}\sbra*{\loccoordtildestd}\in\setR[\dime+\ineqdime+\eqdime]\backslash\brc*{0}$.}
        Since $\D\loccoordbarrKKTvecfld\paren*{\cdot, \barrparamotriter}\paren*{\loccoord{\allvaraccum}}$ is independent of $\barrparamotriter$ by \cref{lemm:loccoordJacobnonsingaccum}, again, we also define 
        \isextendedversion{
        \begin{align}
            \constsit\coloneqq\opnorm{\inv{\D\loccoordbarrKKTvecfld\paren*{\cdot, \barrparamotriter}}\paren*{\loccoord{\allvaraccum}}} > 0
        \end{align}}{$\constsit\coloneqq\opnorm{\inv{\D\loccoordbarrKKTvecfld\paren*{\cdot, \barrparamotriter}}\paren*{\loccoord{\allvaraccum}}} > 0$} regardless of the barrier parameter.
        Since $\opnorm{\D\chart[\pt]\paren*{\cdot}}$ is bounded around $\ptaccum$ due to the continuous differentiability of $\chart[\pt]$, we obtain
        \isextendedversion{
        \begin{align}
            &\norm*{\inv{\D\loccoordbarrKKTvecfld\paren*{\cdot, \barrparamotriter}}\paren*{\loccoord{\allvaraccum}}\loccoordbarrKKTvecfld\paren*{\loccoord{\allvarotriterp}, \barrparamotriter}} \leq \constsit \norm*{\loccoordbarrKKTvecfld\paren*{\loccoord{\allvarotriterp}, \barrparamotriter}}\notag\\
            &\leq \constsit \paren*{\opnorm{\D\chart[\pt]\paren*{\ptotriterp}}
            \Riemnorm[\ptotriterp]{\gradstr[\pt]\Lagfun\paren*{\allvarotriterp}} + \norm*{\Ineqfunmat[]\paren*{\ptotriterp}\ineqLagmultotriterp[] - \barrparamotriter\onevec} + \norm*{\eqfun[]\paren*{\ptotriterp}}}\notag\\
            &\leq \frac{\barrparamotriter}{2}\norm*{\allvarnonzeroconst\paren*{0}}\label{ineq:DFFoteriterbound}
        \end{align}}{
        \begin{align}
            \begin{split}\label{ineq:DFFoteriterbound}
                &\norm*{\inv{\D\loccoordbarrKKTvecfld\paren*{\cdot, \barrparamotriter}}\paren*{\loccoord{\allvaraccum}}\loccoordbarrKKTvecfld\paren*{\loccoord{\allvarotriterp}, \barrparamotriter}} \leq \constsit \norm*{\loccoordbarrKKTvecfld\paren*{\loccoord{\allvarotriterp}, \barrparamotriter}} \leq \constsit \paren*{\opnorm{\D\chart[\pt]\paren*{\ptotriterp}}\\
                &\Riemnorm[\ptotriterp]{\gradstr[\pt]\Lagfun\paren*{\allvarotriterp}} + \norm*{\Ineqfunmat[]\paren*{\ptotriterp}\ineqLagmultotriterp[] - \barrparamotriter\onevec} + \norm*{\eqfun[]\paren*{\ptotriterp}}} \leq \frac{\barrparamotriter}{2}\norm*{\allvarnonzeroconst\paren*{0}}
            \end{split}
        \end{align}}
        for any $\otriteridx\in\setNz$ satisfying that $\barrparamotriter$ is sufficiently small, where the second inequality follows from \cref{def:loccoordbarrKKTvecfld} and the third one from \cref{lemm:barrKKTvecflddirNewtonbound} with sufficiently small coefficients $\constgradLagNewton, \constcomplNewton, \consteqvioNewton \in \setRpp[]$.
        For the third term of \cref{eq:loccoordbarrKKTvecfldTaylor}, we have
        \begin{align}
            \begin{split}\label{ineq:residonesmallO}
                &\norm*{\residone} \leq \int_{0}^{1} \opnorm{\D[]\paren*{\loccoordbarrKKTvecfld\paren*{\paren*{\loccoord{\allvaraccum} + \tmeten\loccoord{\allvarotriterp}, 0 + \tmeten\barrparamotriter}} - \loccoordbarrKKTvecfld\paren*{\loccoord{\allvaraccum}, 0}}} \norm*{\paren*{\loccoord{\allvarotriterp} - \loccoord{\allvaraccum}, \barrparamotriter}} \dd{\tmeten}\\
                &=\smallO[\norm*{\paren*{\loccoord{\allvarotriterp} - \loccoord{\allvaraccum}, \barrparamotriter}}] = \smallO[\norm*{\loccoord{\allvarotriterp} - \loccoord{\allvaraccum}} + \barrparamotriter] = \smallO[\barrparamotriter],
            \end{split}  
        \end{align}
        where the first inequality follows from \isextendedversion{\cref{def:residone}}{the definition of $\residone$}, the first equality from the continuous differentiability of $\loccoordbarrKKTvecfld$, and the last equality from \isextendedversion{\cref{ineq:diffloccoordallvarotriterpaccumbigO}}{\cref{lemm:diffallvarotriterpaccumbarrparambigO}}.
        Therefore, for any $\otriteridx\in\setNz$ satisfying that $\barrparamotriter$ is sufficiently small, we have
        \isextendedversion{
        \begin{equation}
            \begin{aligned}[t]\label{ineq:diffallvarotriterpallvaraccumlow}
                &\norm*{\loccoord{\allvarotriterp} - \loccoord{\allvaraccum}}\\
                &= \norm*{\barrparamotriter\inv{\D\loccoordbarrKKTvecfld\paren*{\cdot, \barrparamotriter}}\paren*{\loccoord{\allvaraccum}}
                \sbra*{\loccoordtildestd}
                + \inv{\D\loccoordbarrKKTvecfld\paren*{\cdot, \barrparamotriter}}\paren*{\loccoord{\allvaraccum}}\sbra*{\loccoordbarrKKTvecfld\paren*{\loccoord{\allvarotriterp}, \barrparamotriter}}
                - \inv{\D\loccoordbarrKKTvecfld\paren*{\cdot, \barrparamotriter}}\paren*{\loccoord{\allvaraccum}}\sbra*{\residone}}\\
                &\geq \norm*{\barrparamotriter\inv{\D\loccoordbarrKKTvecfld\paren*{\cdot, \barrparamotriter}}\paren*{\loccoord{\allvaraccum}}
                \sbra*{\loccoordtildestd}
                } 
                - \norm*{\inv{\D\loccoordbarrKKTvecfld\paren*{\cdot, \barrparamotriter}}\paren*{\loccoord{\allvaraccum}}\sbra*{\loccoordbarrKKTvecfld\paren*{\loccoord{\allvarotriterp}, \barrparamotriter}}} - \norm*{\inv{\D\loccoordbarrKKTvecfld\paren*{\cdot, \barrparamotriter}}\paren*{\loccoord{\allvaraccum}}\sbra*{\residone}}\\
                &\geq \barrparamotriter\norm*{\allvarnonzeroconst\paren*{\zerovec[]}} - \frac{\barrparamotriter}{2}\norm*{\allvarnonzeroconst\paren*{0}} - \constsit\norm{\residone}\\
                &= \frac{\barrparamotriter}{2}\norm*{\allvarnonzeroconst\paren*{0}} - \constsit\norm{\residone} > 0,
            \end{aligned}
        \end{equation}}
        {\begin{equation}
            \begin{aligned}[t]\label{ineq:diffallvarotriterpallvaraccumlow}
                &\norm*{\loccoord{\allvarotriterp} - \loccoord{\allvaraccum}}\\
                &= \norm*{\barrparamotriter\inv{\D\loccoordbarrKKTvecfld\paren*{\cdot, \barrparamotriter}}\paren*{\loccoord{\allvaraccum}}
                \sbra*{\loccoordtildestd}
                + \inv{\D\loccoordbarrKKTvecfld\paren*{\cdot, \barrparamotriter}}\paren*{\loccoord{\allvaraccum}}\sbra*{\loccoordbarrKKTvecfld\paren*{\loccoord{\allvarotriterp}, \barrparamotriter}} - \inv{\D\loccoordbarrKKTvecfld\paren*{\cdot, \barrparamotriter}}\paren*{\loccoord{\allvaraccum}}\sbra*{\residone}}\\
                &\geq \barrparamotriter\norm*{\inv{\D\loccoordbarrKKTvecfld\paren*{\cdot, \barrparamotriter}}\paren*{\loccoord{\allvaraccum}}
                \sbra*{\loccoordtildestd}
                }
                - \norm*{\inv{\D\loccoordbarrKKTvecfld\paren*{\cdot, \barrparamotriter}}\paren*{\loccoord{\allvaraccum}} \sbra*{\loccoordbarrKKTvecfld\paren*{\loccoord{\allvarotriterp}, \barrparamotriter}}} - \norm*{\inv{\D\loccoordbarrKKTvecfld\paren*{\cdot, \barrparamotriter}}\paren*{\loccoord{\allvaraccum}}\sbra*{\residone}}\\
                &\geq \barrparamotriter\norm*{\allvarnonzeroconst\paren*{\zerovec[]}} - \frac{\barrparamotriter}{2}\norm*{\allvarnonzeroconst\paren*{0}} - \constsit\norm{\residone} = \frac{\barrparamotriter}{2}\norm*{\allvarnonzeroconst\paren*{0}} - \constsit\norm{\residone} > 0,
            \end{aligned}
        \end{equation}}
        where the equality follows from \cref{eq:loccoordbarrKKTvecfldTaylor} and the nonsingularity of $\D\loccoordbarrKKTvecfld\paren*{\cdot, \barrparamotriter}$ at $\loccoord{\allvaraccum}$, the second inequality from \isextendedversion{\cref{def:allvarnonzeroconst} and \cref{ineq:DFFoteriterbound}}{the definition of $\allvarnonzeroconst\paren*{0}$ and \cref{ineq:DFFoteriterbound}}, and the last one from \cref{ineq:residonesmallO}.
        Hence, from \isextendedversion{\cref{def:allvarnonzeroconst,ineq:residonesmallO,ineq:diffallvarotriterpallvaraccumlow}}{the definition of $\allvarnonzeroconst\paren*{0}$, again, \cref{ineq:residonesmallO}, and \cref{ineq:diffallvarotriterpallvaraccumlow}}, we have $\norm*{\loccoord{\allvarotriterp} - \loccoord{\allvaraccum}} = \bigomega[\barrparamotriter]$, which, together with \cref{lemm:diffallvarotriterpaccumbarrparambigO}, completes the proof.
        \qed
    \end{proof}

\section{Proof of \cref{theo:locnearquadconv}}\label{appx:prooflocnearquadconv}

    \begin{proof}
        By \cref{lemm:diffallvarotriterpaccumbarrparambigtheta}, we have $\frac{\norm*{\loccoord{\allvarotriterp} - \loccoord{\allvaraccum}}}{\norm*{\loccoord{\allvarotriter} - \loccoord{\allvaraccum}}} = \bigtheta[\frac{\barrparamotriter}{\barrparam[\otriteridx-1]}]$, which, together with \cref{assu:barrparamdecreaselowodr}, implies the superlinear convergence, that is,
        \isextendedversion{
        \begin{align}
            \lim_{\otriteridx\to\infty}\frac{\norm*{\loccoord{\allvarotriterp} - \loccoord{\allvaraccum}}}{\norm*{\loccoord{\allvarotriter} - \loccoord{\allvaraccum}}} = 0.
        \end{align}}
        {$lim_{\otriteridx\to\infty}\frac{\norm*{\loccoord{\allvarotriterp} - \loccoord{\allvaraccum}}}{\norm*{\loccoord{\allvarotriter} - \loccoord{\allvaraccum}}} = 0$.}
        If we employ the update \cref{eq:barrparamupdaterule} with any $\constnin\in\paren*{0,1}$, it follows that
        \isextendedversion{
        \begin{align}
            \frac{\norm*{\loccoord{\allvarotriterp} - \loccoord{\allvaraccum}}}{\norm*{\loccoord{\allvarotriter} - \loccoord{\allvaraccum}}^{1+\constnin}} = \bigtheta[\frac{\barrparamotriterp}{\barrparamotriter^{1+\constnin}}] = \bigtheta[1],
        \end{align}}
        {$\frac{\norm*{\loccoord{\allvarotriterp} - \loccoord{\allvaraccum}}}{\norm*{\loccoord{\allvarotriter} - \loccoord{\allvaraccum}}^{1+\constnin}} = \bigtheta[\frac{\barrparamotriterp}{\barrparamotriter^{1+\constnin}}] = \bigtheta[1]$,}
        which is actually the local near-quadratic convergence.
        The proof is complete.
        \qed
    \end{proof}


\section{Proof of \cref{lemm:trsquadposdifardallvaraccum}}\label{appx:prooftrsquadposdifardallvaraccum}
        \begin{proof}
            We first derive auxiliary bounds to prove the positive definiteness of $\trsquad$ around $\paren*{\ptaccum, \ineqLagmultaccum[]}$.
            First, from \cite[Theorem~3]{Debreu1952DefiSemidefiQuadForm}, \isextendedversion{the \SOSC{}~\cref{eq:secondordersufficientcond}}{the \SOSC{}} implies the existence of $\Debreuconst, \Debreueps \in \setRpp[]$ such that, for any $\tanvecone[\ptaccum]\in\tanspc[\ptaccum]\mani\backslash\brc*{\zerovec[\ptaccum]}$, 
            \isextendedversion{
            \begin{equation}
                \begin{aligned}[t]\label{ineq:DebreuSOSCbound}
                    &\metr[\ptaccum]{\paren*{\Hess[]\objfun\paren*{\ptaccum} - \sum_{\ineqidx\in\ineqset}\ineqLagmultaccum[\ineqidx]\Hess[]\ineqfun[\ineqidx]\paren*{\ptaccum}}
                    \sbra*{\tanvecone[\ptaccum]}}{\tanvecone[\ptaccum]} + \sum_{\ineqidx\in\activeineqset\paren*{\ptaccum}} \Debreuconst \metr[\ptaccum]{\gradstr\ineqfun[\ineqidx]\paren*{\ptaccum}}{\tanvecone[\ptaccum]}^{2}\\
                    &\geq \Debreueps \Riemnorm[\ptaccum]{\tanvecone[\ptaccum]}^{2}.
                \end{aligned}
            \end{equation}}{
            \begin{align}
                \begin{split}\label{ineq:DebreuSOSCbound}
                    &\metr[\ptaccum]{\paren*{\Hess[]\objfun\paren*{\ptaccum} - \sum_{\ineqidx\in\ineqset}\ineqLagmultaccum[\ineqidx]\Hess[]\ineqfun[\ineqidx]\paren*{\ptaccum}}
                    \sbra*{\tanvecone[\ptaccum]}}{\tanvecone[\ptaccum]} + \sum_{\ineqidx\in\activeineqset\paren*{\ptaccum}} \Debreuconst \metr[\ptaccum]{\gradstr\ineqfun[\ineqidx]\paren*{\ptaccum}}{\tanvecone[\ptaccum]}^{2}\\
                    &\geq \Debreueps \Riemnorm[\ptaccum]{\tanvecone[\ptaccum]}^{2}.
                \end{split}
            \end{align}}
            Second, \isextendedversion{\eqscref{ineq:HessobjfunLipschitzbound,ineq:HessineqfunLipschitzbound} imply}{\cref{assu:HessobjineqfunLipschitz} implies} that there exists $\Lipschitzconstfou[\Lagfun] > 0$ such that 
            \begin{align}
                &\opnorm{\Hess[]\objfun\paren*{\pt} - \sum_{\ineqidx\in\ineqset}\ineqLagmult[\ineqidx]\Hess[]\ineqfun[\ineqidx]\paren*{\pt} - \partxp[]{\pt}{\ptaccum}\circ\paren*{\Hess[]\objfun\paren*{\ptaccum} - \sum_{\ineqidx\in\ineqset}\ineqLagmultaccum[\ineqidx]\Hess[]\ineqfun[\ineqidx]\paren*{\ptaccum}}\circ\partxp[]{\ptaccum}{\pt}}\notag\\
                &\leq \paren*{\Lipschitzconstfou[\objfun] + \sum_{\ineqidx\in\ineqset}\ineqLagmult[\ineqidx] \Lipschitzconstfou[{\ineqfun[\ineqidx]}]} \Riemdist{\pt[]}{\ptaccum} \leq \Lipschitzconstfou[\Lagfun]\Riemdist{\pt[]}{\ptaccum}\label{ineq:HessLagfunLipschitzbound}
            \end{align}
            for any $\paren*{\pt, \ineqLagmult[]}\in\mani\times\setR[\ineqdime]$ sufficiently close to $\paren*{\ptaccum, \ineqLagmultaccum[]}$.
            Third,
            \isextendedversion{for each $\ineqidx\in\ineqset$, since $\opnorm{\Hess\ineqfun[\ineqidx]\paren*{\cdot}}$ is bounded around $\ptaccum$ by the twice continuous differentiability of $\brc*{\ineqfun[\ineqidx]}_{\ineqidx\in\ineqset}$, there exists $\Lipschitzconstfou[{\ineqfun[\ineqidx]}] > 0$ such that
            \begin{align}\label{ineq:diffpartxpgradLipschitzbound}
                \Riemnorm[\pt]{\gradstr\ineqfun[\ineqidx]\paren*{\pt}-\partxp[]{\pt}{\ptaccum}\sbra*{\gradstr\ineqfun[\ineqidx]\paren*{\ptaccum}}}  \leq \Lipschitzconstfou[{\ineqfun[\ineqidx]}] \Riemdist{\pt[]}{\ptaccum}
            \end{align}
            for any $\pt\in\mani$ sufficiently close to $\ptaccum$ by \cite[Corollaries~10.47, 10.48]{Boumal23IntroOptimSmthMani}. 
            In addition, since the parallel transport is isometric and the functions $\brc*{\ineqfun[\ineqidx]}_{\ineqidx\in\ineqset}$ are continuously differentiable, there exists $\tholdvalthr > 0$ such that, for any $\pt\in\mani$ sufficiently close to $\ptaccum$, 
            \begin{align}
                \begin{split}\label{ineq:gradsboundardptaccum}
                    &\Riemnorm[\pt]{\gradstr\ineqfun[\ineqidx]\paren*{\pt}} + \Riemnorm[\pt]{\partxp[]{\pt}{\ptaccum}\sbra*{\gradstr\ineqfun[\ineqidx]\paren*{\ptaccum}}}\\
                    &= \Riemnorm[\pt]{\gradstr\ineqfun[\ineqidx]\paren*{\pt}} + \Riemnorm[\ptaccum]{\gradstr\ineqfun[\ineqidx]\paren*{\ptaccum}} \leq \tholdvalthr.
                \end{split}
            \end{align}}{
            there exist positive scalars $\brc*{\Lipschitzconstfou[{\ineqfun[\ineqidx]}]}$ and $\tholdvalthr > 0$ such that
            \begin{align}
                \begin{split}\label{ineq:gradLipschitzPartxpbound}
                    &\Riemnorm[\pt]{\gradstr\ineqfun[\ineqidx]\paren*{\pt}-\partxp[]{\pt}{\ptaccum}\sbra*{\gradstr\ineqfun[\ineqidx]\paren*{\ptaccum}}}  \leq \Lipschitzconstfou[{\ineqfun[\ineqidx]}] \Riemdist{\pt[]}{\ptaccum},\\
                    &\Riemnorm[\pt]{\gradstr\ineqfun[\ineqidx]\paren*{\pt}} + \Riemnorm[\pt]{\partxp[]{\pt}{\ptaccum}\sbra*{\gradstr\ineqfun[\ineqidx]\paren*{\ptaccum}}} = \Riemnorm[\pt]{\gradstr\ineqfun[\ineqidx]\paren*{\pt}} + \Riemnorm[\ptaccum]{\gradstr\ineqfun[\ineqidx]\paren*{\ptaccum}} \leq \tholdvalthr,
                \end{split}
            \end{align}
            for each $\ineqidx\in\ineqset$ and any $\pt\in\mani$ sufficiently close to $\ptaccum$, where the inequality in the first bound follows from Lipschitz continuity of the gradient~\cite[Corollaries~10.47, 10.48]{Boumal23IntroOptimSmthMani} with the boundedness of $\opnorm{\Hess\ineqfun[\ineqidx]\paren*{\cdot}}$ around $\ptaccum$, and the inequality in the second bound from the isometry of the parallel transport and the boundedness of $\brc*{\ineqfun[\ineqidx]}_{\ineqidx\in\ineqset}$.}
            Now, we prove the positive definiteness of $\trsquad\paren*{\pt, \ineqLagmult[]}$.
            Let $\tanvecone[\pt]\in\tanspc[\pt]\mani\backslash\brc*{\zerovec[\pt]}$ be any nonzero vector.
            It follows that
            \isextendedversion{
            \begin{equation}
                \begin{aligned}[t]\label{ineq:trsquadboundone}
                    &\metr[\pt]{\trsquad\paren*{\pt, \ineqLagmult[]}\sbra*{\tanvecone[\pt]}}{\tanvecone[\pt]}\\
                    &= \metr[\pt]{\paren*{\Hess[\pt]\objfun\paren*{\pt} - \sum_{\ineqidx\in\ineqset}\ineqLagmult[\ineqidx]\Hess[\pt]\ineqfun[\ineqidx]\paren*{\pt}}\sbra*{\tanvecone[\pt]}}{\tanvecone[\pt]} + \sum_{\ineqidx\in\ineqset}\frac{\ineqLagmult[\ineqidx]}{\ineqfun[\ineqidx]\paren*{\pt}}\metr[\pt]{\gradstr\ineqfun[\ineqidx]\paren*{\pt}}{\tanvecone[\pt]}^{2}\\
                    &\geq \metr[\pt]{\paren*{\Hess[\pt]\objfun\paren*{\pt} - \sum_{\ineqidx\in\ineqset}\ineqLagmult[\ineqidx]\Hess[\pt]\ineqfun[\ineqidx]\paren*{\pt}}\sbra*{\tanvecone[\pt]}}{\tanvecone[\pt]} + \sum_{\ineqidx\in\activeineqset\paren*{\ptaccum}} \Debreuconst \metr[\ptaccum]{\gradstr\ineqfun[\ineqidx]\paren*{\ptaccum}}{\partxp[]{\ptaccum}{\pt}\sbra*{\tanvecone[\pt]}}^{2}\\
                    &\quad - \sum_{\ineqidx\in\activeineqset\paren*{\ptaccum}} \Debreuconst \metr[\ptaccum]{\gradstr\ineqfun[\ineqidx]\paren*{\ptaccum}}{\partxp[]{\ptaccum}{\pt}\sbra*{\tanvecone[\pt]}}^{2} + \sum_{\ineqidx\in\activeineqset\paren*{\ptaccum}}\frac{\ineqLagmult[\ineqidx]}{\ineqfun[\ineqidx]\paren*{\pt}}\metr[\pt]{\gradstr\ineqfun[\ineqidx]\paren*{\pt}}{\tanvecone[\pt]}^{2},
                \end{aligned}
            \end{equation}}{
            \begin{align}
                &\metr[\pt]{\trsquad\paren*{\pt, \ineqLagmult[]}\sbra*{\tanvecone[\pt]}}{\tanvecone[\pt]}\notag\\
                &= \metr[\pt]{\paren*{\Hess[\pt]\objfun\paren*{\pt} - \sum_{\ineqidx\in\ineqset}\ineqLagmult[\ineqidx]\Hess[\pt]\ineqfun[\ineqidx]\paren*{\pt}}\sbra*{\tanvecone[\pt]}}{\tanvecone[\pt]} + \sum_{\ineqidx\in\ineqset}\frac{\ineqLagmult[\ineqidx]}{\ineqfun[\ineqidx]\paren*{\pt}}\metr[\pt]{\gradstr\ineqfun[\ineqidx]\paren*{\pt}}{\tanvecone[\pt]}^{2}\notag\\
                &\geq \metr[\pt]{\paren*{\Hess[\pt]\objfun\paren*{\pt} - \sum_{\ineqidx\in\ineqset}\ineqLagmult[\ineqidx]\Hess[\pt]\ineqfun[\ineqidx]\paren*{\pt}}\sbra*{\tanvecone[\pt]}}{\tanvecone[\pt]} + \sum_{\ineqidx\in\activeineqset\paren*{\ptaccum}} \Debreuconst \metr[\ptaccum]{\gradstr\ineqfun[\ineqidx]\paren*{\ptaccum}}{\partxp[]{\ptaccum}{\pt}\sbra*{\tanvecone[\pt]}}^{2}\notag\\
                &- \sum_{\ineqidx\in\activeineqset\paren*{\ptaccum}} \Debreuconst \metr[\ptaccum]{\gradstr\ineqfun[\ineqidx]\paren*{\ptaccum}}{\partxp[]{\ptaccum}{\pt}\sbra*{\tanvecone[\pt]}}^{2} + \sum_{\ineqidx\in\activeineqset\paren*{\ptaccum}}\frac{\ineqLagmult[\ineqidx]}{\ineqfun[\ineqidx]\paren*{\pt}}\metr[\pt]{\gradstr\ineqfun[\ineqidx]\paren*{\pt}}{\tanvecone[\pt]}^{2},\label{ineq:trsquadboundone}
            \end{align}
            }
            where the inequality follows from $\activeineqset\paren*{\ptaccum}\subseteq\ineqset$.
            We derive the bound on the first two terms as follows: for all $\paren*{\pt, \ineqLagmult[]}\in\strictfeasirgn\times\setRpp[\ineqdime]$ sufficiently close to $\paren*{\ptaccum, \ineqLagmultaccum[]}$, we have
            \isextendedversion{\begin{align}
                &\metr[\pt]{\paren*{\Hess[\pt]\objfun\paren*{\pt} - \sum_{\ineqidx\in\ineqset}\ineqLagmult[\ineqidx]\Hess[\pt]\ineqfun[\ineqidx]\paren*{\pt}}\sbra*{\tanvecone[\pt]}}{\tanvecone[\pt]} + \sum_{\ineqidx\in\activeineqset\paren*{\ptaccum}} \Debreuconst \metr[\ptaccum]{\gradstr\ineqfun[\ineqidx]\paren*{\ptaccum}}{\partxp[]{\ptaccum}{\pt}\sbra*{\tanvecone[\pt]}}^{2}\\
                &= \metr[\pt]{\paren*{\Hess[\pt]\objfun\paren*{\pt} - \sum_{\ineqidx\in\ineqset}\ineqLagmult[\ineqidx]\Hess[\pt]\ineqfun[\ineqidx]\paren*{\pt}}\sbra*{\tanvecone[\pt]}}{\tanvecone[\pt]}\\
                &\quad - \metr[\pt]{\partxp[]{\pt}{\ptaccum} \circ \paren*{\Hess[\pt]\objfun\paren*{\ptaccum} - \sum_{\ineqidx\in\ineqset}\ineqLagmultaccum[\ineqidx]\Hess[\pt]\ineqfun[\ineqidx]\paren*{\ptaccum}} \circ \partxp[]{\ptaccum}{\pt} \sbra*{\tanvecone[\pt]}}{\tanvecone[\pt]}\\
                &\quad + \metr[\pt]{\partxp[]{\pt}{\ptaccum} \circ \paren*{\Hess[\pt]\objfun\paren*{\ptaccum} - \sum_{\ineqidx\in\ineqset}\ineqLagmultaccum[\ineqidx]\Hess[\pt]\ineqfun[\ineqidx]\paren*{\ptaccum}}\circ \partxp[]{\ptaccum}{\pt} \sbra*{\tanvecone[\pt]}}{\tanvecone[\pt]}\\
                &\quad + \sum_{\ineqidx\in\activeineqset\paren*{\ptaccum}}\Debreuconst\metr[\ptaccum]{\gradstr\ineqfun[\ineqidx]\paren*{\ptaccum}}{\partxp[]{\ptaccum}{\pt}\sbra*{\tanvecone[\pt]}}^{2}\\
                &\geq - \opnorm{\paren*{\Hess[\pt]\objfun\paren*{\pt} - \sum_{\ineqidx\in\ineqset}\ineqLagmult[\ineqidx]\Hess[\pt]\ineqfun[\ineqidx]\paren*{\pt}} - \partxp[]{\pt}{\ptaccum} \circ \paren*{\Hess[\pt]\objfun\paren*{\ptaccum} - \sum_{\ineqidx\in\ineqset}\ineqLagmultaccum[\ineqidx]\Hess[\pt]\ineqfun[\ineqidx]\paren*{\ptaccum}}\circ \partxp[]{\ptaccum}{\pt}}\Riemnorm[\pt]{\tanvecone[\pt]}^{2}\\
                &\quad + \metr[\ptaccum]{\paren*{\Hess[\pt]\objfun\paren*{\ptaccum} - \sum_{\ineqidx\in\ineqset}\ineqLagmultaccum[\ineqidx]\Hess[\pt]\ineqfun[\ineqidx]\paren*{\ptaccum}} \circ \partxp[]{\ptaccum}{\pt} \sbra*{\tanvecone[\pt]}}{\partxp[]{\ptaccum}{\pt} \sbra*{\tanvecone[\pt]}}\\
                &\quad + \sum_{\ineqidx\in\activeineqset\paren*{\ptaccum}}\Debreuconst\metr[\ptaccum]{\gradstr\ineqfun[\ineqidx]\paren*{\ptaccum}}{\partxp[]{\ptaccum}{\pt}\sbra*{\tanvecone[\pt]}}^{2}\\
                &\geq -\Lipschitzconstfou[\Lagfun]\Riemdist{\pt}{\ptaccum}\Riemnorm[\pt]{\tanvecone[\pt]}^{2} + \Debreueps[]\Riemnorm[\pt]{\tanvecone[\pt]}^{2}, \label{ineq:trsquadboundtwo}
            \end{align}}
            {\begin{align}
                &\metr[\pt]{\paren*{\Hess[\pt]\objfun\paren*{\pt} - \sum_{\ineqidx\in\ineqset}\ineqLagmult[\ineqidx]\Hess[\pt]\ineqfun[\ineqidx]\paren*{\pt}}\sbra*{\tanvecone[\pt]}}{\tanvecone[\pt]} + \sum_{\ineqidx\in\activeineqset\paren*{\ptaccum}} \Debreuconst \metr[\ptaccum]{\gradstr\ineqfun[\ineqidx]\paren*{\ptaccum}}{\partxp[]{\ptaccum}{\pt}\sbra*{\tanvecone[\pt]}}^{2}\notag\\
                &= \metr[\pt]{\paren*{\Hess[\pt]\objfun\paren*{\pt} - \sum_{\ineqidx\in\ineqset}\ineqLagmult[\ineqidx]\Hess[\pt]\ineqfun[\ineqidx]\paren*{\pt}}\sbra*{\tanvecone[\pt]}}{\tanvecone[\pt]}\notag\\
                &\quad - \metr[\pt]{\partxp[]{\pt}{\ptaccum} \circ \paren*{\Hess[\pt]\objfun\paren*{\ptaccum} - \sum_{\ineqidx\in\ineqset}\ineqLagmultaccum[\ineqidx]\Hess[\pt]\ineqfun[\ineqidx]\paren*{\ptaccum}} \circ \partxp[]{\ptaccum}{\pt} \sbra*{\tanvecone[\pt]}}{\tanvecone[\pt]}\notag\\
                &\quad + \metr[\pt]{\partxp[]{\pt}{\ptaccum} \circ \paren*{\Hess[\pt]\objfun\paren*{\ptaccum} - \sum_{\ineqidx\in\ineqset}\ineqLagmultaccum[\ineqidx]\Hess[\pt]\ineqfun[\ineqidx]\paren*{\ptaccum}}\circ \partxp[]{\ptaccum}{\pt} \sbra*{\tanvecone[\pt]}}{\tanvecone[\pt]}\notag\\
                &+ \sum_{\ineqidx\in\activeineqset\paren*{\ptaccum}}\Debreuconst\metr[\ptaccum]{\gradstr\ineqfun[\ineqidx]\paren*{\ptaccum}}{\partxp[]{\ptaccum}{\pt}\sbra*{\tanvecone[\pt]}}^{2} \geq - \opnorm{\paren*{\Hess[\pt]\objfun\paren*{\pt} - \sum_{\ineqidx\in\ineqset}\ineqLagmult[\ineqidx]\Hess[\pt]\ineqfun[\ineqidx]\paren*{\pt}} \notag\\
                &\quad - \partxp[]{\pt}{\ptaccum} \circ \paren*{\Hess[\pt]\objfun\paren*{\ptaccum} - \sum_{\ineqidx\in\ineqset}\ineqLagmultaccum[\ineqidx]\Hess[\pt]\ineqfun[\ineqidx]\paren*{\ptaccum}}\circ \partxp[]{\ptaccum}{\pt}}\Riemnorm[\pt]{\tanvecone[\pt]}^{2}\notag\\
                &\quad + \metr[\ptaccum]{\paren*{\Hess[\pt]\objfun\paren*{\ptaccum} - \sum_{\ineqidx\in\ineqset}\ineqLagmultaccum[\ineqidx]\Hess[\pt]\ineqfun[\ineqidx]\paren*{\ptaccum}} \circ \partxp[]{\ptaccum}{\pt} \sbra*{\tanvecone[\pt]}}{\partxp[]{\ptaccum}{\pt} \sbra*{\tanvecone[\pt]}}\notag\\
                &\quad + \sum_{\ineqidx\in\activeineqset\paren*{\ptaccum}}\Debreuconst\metr[\ptaccum]{\gradstr\ineqfun[\ineqidx]\paren*{\ptaccum}}{\partxp[]{\ptaccum}{\pt}\sbra*{\tanvecone[\pt]}}^{2} \geq -\Lipschitzconstfou[\Lagfun]\Riemdist{\pt}{\ptaccum}\Riemnorm[\pt]{\tanvecone[\pt]}^{2} + \Debreueps[]\Riemnorm[\pt]{\tanvecone[\pt]}^{2},\label{ineq:trsquadboundtwo}
            \end{align}}
            where the first inequality follows from the adjoint property of the parallel transport and the second one from \cref{ineq:DebreuSOSCbound}, \cref{ineq:HessLagfunLipschitzbound}, and the isometry of the parallel transport.
            As for the last two terms, we have
            \isextendedversion{
            \begin{align}
                &- \sum_{\ineqidx\in\activeineqset\paren*{\ptaccum}} \Debreuconst \metr[\ptaccum]{\gradstr\ineqfun[\ineqidx]\paren*{\ptaccum}}{\partxp[]{\ptaccum}{\pt}\sbra*{\tanvecone[\pt]}}^{2} + \sum_{\ineqidx\in\activeineqset\paren*{\ptaccum}}\frac{\ineqLagmultotriter[\ineqidx]}{\ineqfun[\ineqidx]\paren*{\pt}}\metr[\pt]{\gradstr\ineqfun[\ineqidx]\paren*{\pt}}{\tanvecone[\pt]}^{2}\\
                &= \sum_{\ineqidx\in\activeineqset\paren*{\ptaccum}}\Debreuconst\metr[\pt]{\gradstr\ineqfun[\ineqidx]\paren*{\pt}}{\tanvecone[\pt]}^{2} -\sum_{\ineqidx\in\activeineqset\paren*{\ptaccum}} \Debreuconst \metr[\pt]{\partxp[]{\pt}{\ptaccum}\sbra*{\gradstr\ineqfun[\ineqidx]\paren*{\ptaccum}}}{\tanvecone[\pt]}^{2} \\
                &\quad -\sum_{\ineqidx\in\activeineqset\paren*{\ptaccum}}\Debreuconst\metr[\pt]{\gradstr\ineqfun[\ineqidx]\paren*{\pt}}{\tanvecone[\pt]}^{2} + \sum_{\ineqidx\in\activeineqset\paren*{\ptaccum}}\frac{\ineqLagmultotriter[\ineqidx]}{\ineqfun[\ineqidx]\paren*{\pt}}\metr[\pt]{\gradstr\ineqfun[\ineqidx]\paren*{\pt}}{\tanvecone[\pt]}^{2}\\
                &= \sum_{\ineqidx\in\activeineqset\paren*{\ptaccum}} \Debreuconst \paren*{\metr[\pt]{\gradstr\ineqfun[\ineqidx]\paren*{\pt} + \partxp[]{\pt}{\ptaccum}\sbra*{\gradstr\ineqfun[\ineqidx]\paren*{\ptaccum}}}{\tanvecone[\pt]} \metr[\pt]{\gradstr\ineqfun[\ineqidx]\paren*{\pt} - \partxp[]{\pt}{\ptaccum}\sbra*{\gradstr\ineqfun[\ineqidx]\paren*{\ptaccum}}}{\tanvecone[\pt]}}\\
                &\quad + \sum_{\ineqidx\in\activeineqset\paren*{\ptaccum}}\paren*{\frac{\ineqLagmultotriter[\ineqidx]}{\ineqfun[\ineqidx]\paren*{\pt}}- \Debreuconst}\metr[\pt]{\gradstr\ineqfun[\ineqidx]\paren*{\pt}}{\tanvecone[\pt]}^{2}\\
                &\geq -\sum_{\ineqidx\in\activeineqset\paren*{\ptaccum}}\Debreuconst \paren*{\Riemnorm[\pt]{\gradstr\ineqfun[\ineqidx]\paren*{\pt}} + \Riemnorm[\pt]{\partxp[]{\pt}{\ptaccum}\sbra*{\gradstr\ineqfun[\ineqidx]\paren*{\ptaccum}}}} \Riemnorm[\pt]{\gradstr\ineqfun[\ineqidx]\paren*{\pt}-\partxp[]{\pt}{\ptaccum}\sbra*{\gradstr\ineqfun[\ineqidx]\paren*{\ptaccum}}} \Riemnorm[\pt]{\tanvecone[\pt]}^{2}\\
                &\quad + \sum_{\ineqidx\in\activeineqset\paren*{\ptaccum}}\paren*{\frac{\ineqLagmultotriter[\ineqidx]}{\ineqfun[\ineqidx]\paren*{\pt}} - \Debreuconst}\metr[\pt]{\gradstr\ineqfun[\ineqidx]\paren*{\pt}}{\tanvecone[\pt]}^{2}\\
                &\geq - \sum_{\ineqidx\in\activeineqset\paren*{\ptaccum}}\Debreuconst \tholdvalthr \Lipschitzconstfou[{\ineqfun[\ineqidx]}] \Riemdist{\pt}{\ptaccum}\Riemnorm[\pt]{\tanvecone[\pt]}^{2} + \sum_{\ineqidx\in\activeineqset\paren*{\ptaccum}}\paren*{\frac{\ineqLagmult[\ineqidx]}{\ineqfun[\ineqidx]\paren*{\pt}} - \Debreuconst}\metr[\pt]{\gradstr\ineqfun[\ineqidx]\paren*{\pt}}{\tanvecone[\pt]}^{2},\label{ineq:trsquadboundthr}
            \end{align}}
            {\begin{align}
                \begin{split}\label{ineq:trsquadboundthr}
                    &- \sum_{\ineqidx\in\activeineqset\paren*{\ptaccum}} \Debreuconst \metr[\ptaccum]{\gradstr\ineqfun[\ineqidx]\paren*{\ptaccum}}{\partxp[]{\ptaccum}{\pt}\sbra*{\tanvecone[\pt]}}^{2} + \sum_{\ineqidx\in\activeineqset\paren*{\ptaccum}}\frac{\ineqLagmultotriter[\ineqidx]}{\ineqfun[\ineqidx]\paren*{\pt}}\metr[\pt]{\gradstr\ineqfun[\ineqidx]\paren*{\pt}}{\tanvecone[\pt]}^{2}\\
                    &= \sum_{\ineqidx\in\activeineqset\paren*{\ptaccum}}\Debreuconst\metr[\pt]{\gradstr\ineqfun[\ineqidx]\paren*{\pt}}{\tanvecone[\pt]}^{2} -\sum_{\ineqidx\in\activeineqset\paren*{\ptaccum}} \Debreuconst \metr[\pt]{\partxp[]{\pt}{\ptaccum}\sbra*{\gradstr\ineqfun[\ineqidx]\paren*{\ptaccum}}}{\tanvecone[\pt]}^{2}\\
                    &\quad -\sum_{\ineqidx\in\activeineqset\paren*{\ptaccum}}\Debreuconst\metr[\pt]{\gradstr\ineqfun[\ineqidx]\paren*{\pt}}{\tanvecone[\pt]}^{2} + \sum_{\ineqidx\in\activeineqset\paren*{\ptaccum}}\frac{\ineqLagmultotriter[\ineqidx]}{\ineqfun[\ineqidx]\paren*{\pt}}\metr[\pt]{\gradstr\ineqfun[\ineqidx]\paren*{\pt}}{\tanvecone[\pt]}^{2} = \sum_{\ineqidx\in\activeineqset\paren*{\ptaccum}} \Debreuconst \paren*{\metr[\pt]{\gradstr\ineqfun[\ineqidx]\paren*{\pt} \\
                    &\quad + \partxp[]{\pt}{\ptaccum}\sbra*{\gradstr\ineqfun[\ineqidx]\paren*{\ptaccum}}}{\tanvecone[\pt]} \metr[\pt]{\gradstr\ineqfun[\ineqidx]\paren*{\pt} - \partxp[]{\pt}{\ptaccum}\sbra*{\gradstr\ineqfun[\ineqidx]\paren*{\ptaccum}}}{\tanvecone[\pt]}}\\
                    &\quad + \sum_{\ineqidx\in\activeineqset\paren*{\ptaccum}}\paren*{\frac{\ineqLagmultotriter[\ineqidx]}{\ineqfun[\ineqidx]\paren*{\pt}} - \Debreuconst}\metr[\pt]{\gradstr\ineqfun[\ineqidx]\paren*{\pt}}{\tanvecone[\pt]}^{2} \geq -\sum_{\ineqidx\in\activeineqset\paren*{\ptaccum}}\Debreuconst \paren*{\Riemnorm[\pt]{\gradstr\ineqfun[\ineqidx]\paren*{\pt}}\\
                    &\quad + \Riemnorm[\pt]{\partxp[]{\pt}{\ptaccum}\sbra*{\gradstr\ineqfun[\ineqidx]\paren*{\ptaccum}}}}\Riemnorm[\pt]{\gradstr\ineqfun[\ineqidx]\paren*{\pt}-\partxp[]{\pt}{\ptaccum}\sbra*{\gradstr\ineqfun[\ineqidx]\paren*{\ptaccum}}} \Riemnorm[\pt]{\tanvecone[\pt]}^{2}\\
                    &\quad + \sum_{\ineqidx\in\activeineqset\paren*{\ptaccum}}\paren*{\frac{\ineqLagmultotriter[\ineqidx]}{\ineqfun[\ineqidx]\paren*{\pt}} - \Debreuconst}\metr[\pt]{\gradstr\ineqfun[\ineqidx]\paren*{\pt}}{\tanvecone[\pt]}^{2}\\
                    &\geq - \sum_{\ineqidx\in\activeineqset\paren*{\ptaccum}}\Debreuconst \tholdvalthr \Lipschitzconstfou[{\ineqfun[\ineqidx]}] \Riemdist{\pt}{\ptaccum}\Riemnorm[\pt]{\tanvecone[\pt]}^{2} + \sum_{\ineqidx\in\activeineqset\paren*{\ptaccum}}\paren*{\frac{\ineqLagmult[\ineqidx]}{\ineqfun[\ineqidx]\paren*{\pt}} - \Debreuconst}\metr[\pt]{\gradstr\ineqfun[\ineqidx]\paren*{\pt}}{\tanvecone[\pt]}^{2},
                \end{split}
            \end{align}}
            where the second inequality follows from \isextendedversion{\cref{ineq:diffpartxpgradLipschitzbound,ineq:gradsboundardptaccum}}{\cref{ineq:gradLipschitzPartxpbound}}.
            Therefore, combining \cref{ineq:trsquadboundone} with \cref{ineq:trsquadboundtwo,ineq:trsquadboundthr} yields that, for any $\paren*{\pt, \ineqLagmult[]}\in\strictfeasirgn\times\setRpp[\ineqdime]$ sufficiently close to $\paren*{\ptaccum, \ineqLagmultaccum[]}$,
            \isextendedversion{
            \begin{align}
                &\metr[\pt]{\trsquad\paren*{\pt, \ineqLagmult[]}\sbra*{\tanvecone[\pt]}}{\tanvecone[\pt]} \\
                &\geq \paren*{\Debreueps[] -\paren*{\Lipschitzconstfou[\Lagfun] + \Debreuconst \tholdvalthr \sum_{\ineqidx\in\activeineqset\paren*{\ptaccum}}\Lipschitzconstfou[{\ineqfun[\ineqidx]}]}\Riemdist{\pt}{\ptaccum}}\Riemnorm[\pt]{\tanvecone[\pt]}^{2} + \sum_{\ineqidx\in\activeineqset\paren*{\ptaccum}}\paren*{\frac{\ineqLagmultotriter[\ineqidx]}{\ineqfun[\ineqidx]\paren*{\pt}} - \Debreuconst}\metr[\pt]{\gradstr\ineqfun[\ineqidx]\paren*{\pt}}{\tanvecone[\pt]}^{2}.
            \end{align}}
            {$\metr[\pt]{\trsquad\paren*{\pt, \ineqLagmult[]}\sbra*{\tanvecone[\pt]}}{\tanvecone[\pt]} \geq \paren*{\Debreueps[] -\paren*{\Lipschitzconstfou[\Lagfun] + \Debreuconst \tholdvalthr \sum_{\ineqidx\in\activeineqset\paren*{\ptaccum}} \Lipschitzconstfou[{\ineqfun[\ineqidx]}]}\Riemdist{\pt}{\ptaccum}}\Riemnorm[\pt]{\tanvecone[\pt]}^{2} + \sum_{\ineqidx\in\activeineqset\paren*{\ptaccum}}\paren*{\frac{\ineqLagmultotriter[\ineqidx]}{\ineqfun[\ineqidx]\paren*{\pt}} - \Debreuconst}\metr[\pt]{\gradstr\ineqfun[\ineqidx]\paren*{\pt}}{\tanvecone[\pt]}^{2}$.}
            By re-taking $\paren*{\pt, \ineqLagmult[]}\in\strictfeasirgn\times\setRpp[\ineqdime]$ sufficiently close to $\paren*{\ptaccum, \ineqLagmultaccum[]}$, $\Riemdist{\pt}{\ptaccum}$ and $\brc*{\frac{\ineqLagmultotriter[\ineqidx]}{\ineqfun[\ineqidx]\paren*{\pt}}}_{\ineqidx\in\activeineqset\paren*{\ptaccum}}$ can be made arbitrarily small and large, respectively.
            This implies that, for any $\paren*{\pt, \ineqLagmult[]}\in\strictfeasirgn\times\setRpp[\ineqdime]$ sufficiently close to $\paren*{\ptaccum, \ineqLagmultaccum[]}$, the right-hand side is positive for any $\tanvecone[\pt]\in\tanspc[\pt]\mani\backslash\brc*{\zerovec[\ptaccum]}$, which completes the proof.
            \qed
        \end{proof}

    \isextendedversion{

\section{Proof of \cref{prop:dirNewtonardaccum}}\label{appx:proofdirNewtonardaccum}

    \begin{proof}
        By the definition of $\direxactineqLagmult[]$, we have
        \begin{align}\label{eq:Jacobcompl}
           \IneqLagmultmat[]\coineqgradopr[\pt]\sbra*{\direxactpt} + \Ineqfunmat[]\paren*{\pt}\direxactineqLagmult[] = - \paren*{\Ineqfunmat[]\paren*{\pt}\ineqLagmult[] - \barrparam[]\onevec}.
        \end{align}
        Recall that, by \cref{prop:trsgloboptimiffcond}, there exists $\trsLagmult \geq 0$ such that $\direxactpt$ and $\trsLagmult$ satisfy \cref{eq:trsgloboptimiffcond}.
        From \cref{eq:trsgloboptimiffcondconmpl} and $\Riemnorm[\pt]{\direxactpt} < \trradius[]$, we have $\trsLagmult = 0$.
        Thus, it follows that 
        \begin{align}
            \begin{split}\label{eq:JacobHessLag}
                &\Hess[\pt]\Lagfun\paren*{\allvar}\sbra*{\direxactpt} - \ineqgradopr[\pt]\sbra*{\direxactineqLagmult[]}\\
                &= \Hess[\pt]\Lagfun\paren*{\allvar}\sbra*{\direxactpt} + \ineqgradopr[\pt]\sbra*{\IneqLagmultmat[]\inv{\Ineqfunmat[]\paren*{\pt}}\coineqgradopr\sbra*{\direxactpt}} - \ineqgradopr[\pt]\sbra*{\barrparam[]\inv{\Ineqfunmat[]\paren*{\pt}}\onevec} + \ineqgradopr[\pt]\sbra*{\ineqLagmult[]}\\
                &= - \gradstr[]\objfun\paren*{\pt} + \ineqgradopr[\pt]\sbra*{\ineqLagmult[]} = - \gradstr[\pt]\Lagfun\paren*{\allvar[]},
            \end{split}
        \end{align}
        where the first equality follows from \cref{eq:Jacobcompl}, the second one from \cref{eq:trsquaddef,eq:trslindef,eq:trsgloboptimiffcondquadlineq}, and the third one from \cref{eq:gradLagfundef}. 
        By \cref{eq:Jacobcompl,eq:JacobHessLag}, we have that $\paren*{\direxactpt, \direxactineqLagmult[]}$ is the solution of \cref{eq:RiemNewtoneq}.
        Since the solution is unique due to the nonsingularity of $\Jacobian[\barrKKTvecfld]\paren*{\allvar}$,
        we conclude that $\paren*{\direxactpt, \direxactineqLagmult[]}$ is equivalent to the Newton step~$\paren*{\dirNewtonptbarrparam, \dirNewtonineqLagmultbarrparam[]}\in\tanspc[\pt]\mani\times\setR[\ineqdime]$.       
        \qed
    \end{proof}

\section{Proof of \cref{lemm:dirNewtonwelldefinedbound}}\label{appx:proofdirNewtonwelldefinedbound}

    \begin{proof}
        Since $\Jacobian[\barrKKTvecfld]\paren*{\allvar[]}$ is continuous,
        it follows from \cref{lemm:Jacobnonsingaccum} that, for any $\allvar\in\feasirgn\times\setRp[\ineqdime]$ sufficiently close to $\allvaraccum$,
        $\Jacobian[\barrKKTvecfld]$ is nonsingular and 
        \begin{align}
            \Riemnorm[{\allvar[]}]{\dirNewtonallvarbarrparam} 
            \leq \opnorm{\inv{\Jacobian[\barrKKTvecfld]\paren*{\allvar[]}}} \Riemnorm[{\allvar[]}]{\barrKKTvecfld\paren*{\allvar[];\barrparam[]}} \leq \consteig \paren*{\Riemnorm[{\allvar[]}]{\barrKKTvecfld\paren*{\allvar[]}} + \sqrt{\ineqdime}\barrparam[]}
        \end{align}
        for some $\consteig > 0$.
        Since the point $\allvaraccum$ is the solution of $\barrKKTvecfld\paren*{\allvar[]} = \zerovec[\allvaraccum]$ and $\barrKKTvecfld$ is continuous, the right-hand side can be made arbitrarily small by choosing a sufficiently small $\barrparam[]$.
        The proof is complete.
        \qed
    \end{proof}

    }{}

    \bibliographystyle{spmpsci}
    \bibliography{reference.bib}

\end{document}